\documentclass{article}
\usepackage{array}
\usepackage{amsmath}
\usepackage{amssymb}
\usepackage{amsthm}
\usepackage{amsfonts}
\usepackage{graphicx}
\usepackage[toc,page]{appendix}
\usepackage[margin=0pt,font+=small,labelformat=parens,labelsep=space, skip=6pt,list=false,hypcap=false] {caption}
\usepackage{subfig}
\usepackage{maplestd2e}
\theoremstyle{definition}
\newtheorem{theorem}{Lemma}[section]

\DefineParaStyle{Maple Heading 1}
\DefineParaStyle{Maple Text Output}
\DefineParaStyle{Maple Dash Item}
\DefineParaStyle{Maple Bullet Item}
\DefineParaStyle{Maple Normal}
\DefineParaStyle{Maple Heading 4}
\DefineParaStyle{Maple Heading 3}
\DefineParaStyle{Maple Heading 2}
\DefineParaStyle{Maple Warning}
\DefineParaStyle{Maple Title}
\DefineParaStyle{Maple Error}
\DefineCharStyle{Maple Hyperlink}
\DefineCharStyle{Maple 2D Math}
\DefineCharStyle{Maple Maple Input}
\DefineCharStyle{Maple 2D Output}
\DefineCharStyle{Maple 2D Input}

\usepackage[top=.9in,bottom=.9in]{geometry}
 \numberwithin{equation}{section}
 
\DefineParaStyle{Maple Output}
\DefineCharStyle{2D Comment}
\DefineCharStyle{2D Math}
\DefineCharStyle{2D Output}

\begin{document}

%
\let\labeldefs\iffalse
\let\ifjournal\iftrue 
\pagestyle{plain}
\begin{maplegroup}
\numberwithin{equation}{section}
\begin{flushleft} \vskip 0.3 in 
\iftrue
\centerline{An Integral Equation for Riemann's Zeta Function and its Approximate Solution - Revised} \vskip .3in 
\vskip .2in
\centerline{ Michael Milgram\footnote{mike@geometrics-unlimited.com}}
\centerline{Consulting Physicist, Geometrics Unlimited, Ltd.}
\centerline{Box 1484, Deep River, Ont. Canada. K0J 1P0}
\centerline{}
\centerline{ Dec. 17, 2018}
\centerline{}
Revisions: Dec. 19, 2018 - Substantial revisions to Section 8.\newline
Revisions: Jan. 2, 2019 - Minor revisions, clarifications and corrections throughout;\newline 
\centerline{added Appendix \ref{sec:RomSums} and augmented Section labelled ``More"; revised notation: $\rho\rightarrow t$ and $\Lambda(s)\rightarrow \Upsilon(s)$}.
Revisions: Jan. 7, 2019: Eq. (D.1) corrected;\newline
Revison: July 18, 2019: Extensive reorganization of the text and revisions: Sections 6 and 8 completely rewritten to clarify and to address objections
\centerline{}
\vskip .1in
\fi
\centerline{}
\vskip .1in
MSC classes: 	11M06, 11M26, 11M99, 26A09, 30B40, 30E20, 30C15, 33C47, 33B99, 33F99
\vskip 0.1in
\centerline{\bf Abstract}\vskip .3in

Two identities extracted from the literature are coupled to obtain an integral equation for Riemann's $\xi(s)$ function, and thus $\zeta(s)$ indirectly. The equation has a number of simple properties from which useful derivations flow, the most notable of which relates $\zeta(s)$ anywhere in the critical strip to its values on a line anywhere else in the complex plane. From this, I obtain both an analytic expression for $\zeta(\sigma+i{\it t})$ everywhere inside the asymptotic (${\it t}\rightarrow\infty)$ critical strip, and an approximate solution, within the confines of which the Riemann Hypothesis is shown to be true. The approximate solution predicts a simple, but strong correlation between the real and imaginary components of $\zeta(\sigma+i{\it t})$ for different values of $\sigma$ and equal values of ${\it t}$; this is illustrated in a number of figures. 

\section{Introduction} \label{sec:Intro}
The Riemann Zeta function $\zeta(s)$ is well-known to satisfy a functional equation, and many representations, both integral and series, have been developed over the years for both $\zeta(s)$ and its avatar $\xi(s)$. Additionally, and more significantly, at least two independent contour integral representations are known, \cite[Section (1.4)]{Edwards} and \cite[Eq.(7)]{Milgram} either of which could be utilized as the primary definition, from which many of the properties of $\zeta(s)$ can be derived. In contrast, it has been long-ago proven that $\zeta(s)$ does not satisfy a differential equation of quite general form \cite [(and citations therein)] {Kacinskaite} and I am aware of three integral equations which $\zeta(s)$ does satisfy \cite[Eq.(1.6)]{Patkowski}, \cite[Eq.(3.29)]{Sugiyama}, and \cite[Eq.(1.5) with $\alpha=0$]{AshTFok}. In this work, a series representation of $\xi(s)$ in terms of generalized Integro-exponential functions ($E_s(z)$), is coupled with a contour integral representation of these functions, to obtain a new integral equation satisfied by $\xi(s)$ and, equivalently, $\zeta(s)$. This equation has several remarkable properties:
\begin{itemize}
\item
Through the action of an integral operator, the value of $\zeta(s_1)$, (the ``dependent") anywhere in the complex plane can be determined with respect to the properties of $\zeta(s_2)$, (the ``master") on a line ($\Re(s_2)=constant)$ anywhere else  in the punctured ($s_1\neq s_2$) complex plane;
\item
The transfer function that mediates the aforementioned action is a simple rational polynomial function of $s_1$, and is therefore quite amenable to analysis;
\item
The only $s_1$ dependence resides in the transfer function.
\end{itemize}

In this work, after initially listing some definitions and lemmas (Section \ref{sec:Defs}), I first give a derivation of the integral equation over a bounded region of the $s_2-$plane (Section \ref{sec:LeClair}), and obtain its analytic continuation over the remainder of the punctured ($s_1\neq s_2$) complex plane (Section \ref{sec:AnalCont}). In several subsections of Section \ref{sec:Integrals}, the various representations are then used to obtain a few simple integrals involving $\zeta(s_2)$ corresponding to special values of $s_1$; it is also demonstrated that the simple integrals under consideration are convergent. \newline

Following these preliminaries, the general form of the transfer function is carefully presented in a series of Appendices. Based on a particularly useful property of the transfer function presented in Appendix \ref{sec:Tprops} - it closely mimics a Dirac delta function on the $s_2$ line - a reasonably accurate model for $\zeta(s_1)$ is established, which yields reasonably accurate estimates of the asymptotic nature of $\zeta(s_1)$ inside the critical strip (Section \ref{sec:Asympt}). In consequence, it is shown that, asymptotically 
\begin{equation} \label{Zhalf}
\left| \zeta(\sigma+i{\it t})\right| \sim {\it t}^{1/2-\sigma/2}\times(\log(\log({\it t}))/\log({\it t})+...)\,,
\end{equation}
and, {\bf within the confines of the model}, it is both proven that $\zeta(s_1)\neq 0$ if $\Re(s_1)\neq 1/2$, and explained, in transparent terms (Section \ref{sec:zeros}), why this happens (does not happen?). Graphical examples and comparisons are given, demonstrating the accuracy of a predicted universality between the real and imaginary components of $\xi(\sigma+i{\it t})$ at equal values of ${\it t}$ but different $\sigma$. Section \ref{sec:summary} discusses the requirements to improve the rigour of the model developed here. 

\section{Preamble} \label{sec:Defs}
\subsection{Notation}

Throughout, I use $s=\sigma+i{\it t}$ to denote the independent variable defining $\zeta(s)$. In those instances where dependence on, for example, ${\it t}$ is the focus, I will sometimes just shorten, for example $M(c,s,v)\Rightarrow M({\it t},v)$ for typographical brevity and clarity. Throughout, subscripts 'R' and 'I' refer to the Real and Imaginary components of whatever they are attached to. Much use is made of Riemann's $\xi$ function, defined by
\begin{equation}
\xi(s)\equiv \displaystyle { \left( s-1 \right) {\pi}^{-s/2}}{}\,\zeta \left( s \right) \Gamma \left(1+ s/2 \right)
\label{xidef}
\end{equation}

and

\begin{equation} 
\Upsilon(s)\equiv \zeta \left( s \right)\Gamma \left( s/2 \right)\pi^{-s/2}
\label{Lamdef}
\end{equation}

both of which satisfy
\begin{align} \label{Reflxi}
&\xi(s)=\xi(1-s)\\ \label{ReflLam}
&\Upsilon(s)=\Upsilon(1-s)\,.
\end{align}
\newline
Throughout, $k,N=0,1,2...$.
 
\subsection{Definitions and Lemmas}
\begin{itemize}

\item 
\text {\parbox{.85\textwidth}{{\bf Definition:} polar form }}\newline
\label{polar}
means that a complex function $f(s)$ is written as $\exp(i\theta(s))|f(s)|$ where $\theta(s) \equiv \arg(f(s))$.\newline
\item
The following asymptotic limits \cite[Eq. (5.6.9)] {NIST} - see also \cite[Eq. (4.12.2)]{Titch2} - will be required:

\begin{equation}
\lim_{\substack{{\it t}\rightarrow\infty}}\displaystyle  \left| \Gamma \left( \sigma/2+i{\it t}/2 \right)  \right| \approx \sqrt{2\,\pi} \left( {\sigma}^{2}/4+{{\it t}}^{2}/4 \right) ^{(\sigma-1)/4}
\exp{(-\pi\, {\it t}/4) }\,,
\label{LgammaGen}
\end{equation}
and the special case
\begin{equation}
\lim_{\substack{{\it t}\rightarrow\infty}} \displaystyle \left|\Gamma(1/2\pm i{\it t}/2)\right| \approx \sqrt{2\pi}\,\exp(-\pi{\it t}/4)\,.
\label{LgammaSpec}
\end{equation}
(Remark: I have tested the approximation \eqref{LgammaGen} numerically, and find  that it is remarkably accurate for even modest values of $t$ when $0<\sigma<1\ll{\it t}$.)\newline
\item
From the functional equation of $\zeta(s)$ with \eqref{LgammaGen}, for large values of t
\begin{equation}
\displaystyle  \left| \zeta \left( it \right)  \right| \sim \left| \zeta \left( 1-it \right)  \right| \,\sqrt{\frac {  t}{ {2\pi }}}
\label{FeqId}
\end{equation}


\item \begin{theorem}
\begin{equation}
\arg(\zeta(1/2+i{\it t})) +\arg(\Gamma(1/4+i{\it t}/2))-{\it t}/2\log(\pi)=\pm k\pi
\label{R_S}
\end{equation}
\end{theorem}

Although well-known \cite[Eq. (4.17.2)]{Titch2} with $k=0$, this is sometimes referred to as the Riemann-Siegel identity . A simple derivation follows. \newline

{\bf Proof:}\newline

Expand $\xi(\sigma+it)$ about $\sigma=1/2$ at constant t giving
\begin{equation}
\lim _{\sigma\rightarrow 1/2}\xi \left( \sigma+it \right)  =\displaystyle \xi \left( 1/2+it \right) + \left( {\frac {\rm \partial}{{\rm \partial}{\sigma}}}\xi \left( {\sigma+it} \right)|_{\sigma=1/2} \right)  \left( \sigma-1/2 \right) + ... 
\label{xiExpand}
\end{equation}
and find the imaginary part from \eqref{xidef}, after rewriting $\zeta(s)$ and $\Gamma(s/2)$ in polar form:
\begin{align}
&\Im(\xi \left( 1/2+it \right))\equiv \displaystyle \Im \left( \lim _{\sigma\rightarrow 1/2}\xi \left( \sigma+it \right)  \right)\\ &=-{\frac { \left( 4{t}^{2}+1 \right) \sin \left( \arg(\zeta(\frac{1}{2}+it))+\arg(\Gamma(\frac{1}{4}+it/2))-t/2\,\ln  \left( \pi \right) \right) 
\mbox{} \left| \Gamma \left( \frac{1}{4}+it/2 \right) \zeta \left(\frac{1}{2}+it \right) 
\mbox{} \right| }{8\,{\pi}^{1/4}}}\,.
\label{RsProof}
\end{align}
Because $\xi(\frac{1}{2}+it)$ is known to be real, \eqref{RsProof} vanishes, and \eqref{R_S} follows immediately. {\bf QED} \newline

{\bf Remark:} Apply the identity 
\begin{equation}
\arg(\zeta(\frac{1}{2}+it))=\Im({\log\Gamma(1/2+it))}+N 
\end{equation}
to \eqref{R_S} to obtain Backlund's formula which counts the discontinuities (not zeros)  of  $\arg(\zeta(\frac{1}{2}+it))$ - see \cite[Section 9]{Milgram_Exploring}. This is consistent with \cite[Theorem 1]{Ivic_Note}.\newline

\item
Define
\begin{equation}
\phi(s)\equiv \arg(\zeta(s))+ \arg(\Gamma(s/2)) -(1/2){\it t}\log(\pi)=\arg(\Upsilon(s))\,,
\label{phi}
\end{equation}
\begin{equation}
\theta(v)= \arg(\Gamma(3/2+iv))+\arg(\zeta(1+2iv))-v\ln  \left( \pi \right),
\label{theta}
\end{equation}
\begin{equation}
\Phi(t)\equiv \arg(\zeta(i{\it t}))+ \arg(\Gamma(i{\it t}/2)) -(t/2)\log(\pi)
\label{Phi}
\end{equation}
\begin{equation}
\alpha(t)\equiv \arg{\zeta(1/2+it)}
\label{alpha-def}
\end{equation}
and
\begin{equation}
\beta(t)\equiv \arg{\zeta^{\prime}(1/2+it)}
\label{beta-def}
\end{equation}
but, for brevity, I occasionally abbreviate expressions such as the following: 
\begin{equation}
\cos(\Phi({\it t})):= \cos(\Phi).
\end{equation}

\end{itemize}

\section{LeClair's representation} \label{sec:LeClair}
In a recent work, LeClair \cite[Eq.(15)]{LeClair:2013oda} has obtained the following series representation of Riemann's $\xi$ function 
\begin{align}
\xi(s)= \pi \left( s-1 \right) \sum _{n=1}^{\infty }{n}^{2}{E_{-s/2}} \left(\pi\,{n}^{2} \right)-\pi\,s\sum _{n=1}^{\infty }{n}^{2}{ E_{(s-1)/2}} \left( \pi\,{n}^{2} \right)+4\,\pi\,\sum _{n=1}^{\infty }{n}^{2}{{\rm e}^{-\pi\,{n}^{2}}}
\mbox{} 
\label{Leclair1}
\end{align}
where generically, $E_{s}(z)$ is the (generalized) ``Exponential Integral",  a limiting case of what is elsewhere \cite{Milgram:1985} referred to as the ``Generalized Integro-Exponential Function". Further to the above, it has been shown (\cite[Eq.(3)]{Romik}) that the infinite sum in \eqref{Leclair1} can be readily evaluated:
\begin{equation}
\mapleinline{inert}{2d}{Sum(n^2*exp(-Pi*n^2), n = 1 .. infinity) = (1/8)/(8\,Pi^(3/4)*GAMMA(3/4))}{\[\displaystyle 4\pi\sum _{n=1}^{\infty }{n}^{2}{{\rm e}^{-\pi\,{n}^{2}}}={\frac {\pi^{1/4}}{2\Gamma \left( 3/4 \right) }}\]}\,.
\label{Larry1}
\end{equation}

In his work, LeClair truncates the sum(s) at $N$ terms, refers to the result as an ``approximation" and proceeds to obtain approximations to the location of the zeros on the critical line on that basis. Here, I treat the sums as an infinite series representation of $\xi(s)$, and hence an identity because the series is easily shown to be convergent due to the asymptotic property of $E_s(z)$ ( see \cite[Eq.(2.25)]{Milgram:1985}). Similar, but inequivalent series representations will be found in Paris \cite[Eq.(1.1)] {ParisExp}, Patkowski \cite[Eq.(1.20)] {Patkowski} and elsewhere. From \cite{Milgram:1985}, some useful integral and  contour integral representations of the function $E_{s}(z)$ are
\begin{align}
E_{s}(z)&=z^{s-1}\Gamma(1-s,z) \label{InG}\\
&=\int_{1}^{\infty} v^{-s}\exp(-zv){\rm d}v \label{Edef}
\end{align}
\begin{equation}
\displaystyle {E_{s}} \left(z \right) ={\frac {1}{2\pi i}\int_{c-i\infty }^{c+i\infty }\!{\frac {\Gamma \left( -v \right) {z}^{v}}{s-1
\mbox{}-v}}\,{\rm d}v}
\label{Ei1}
\end{equation}
In \eqref{InG}, $\Gamma(1-s,z)$ is the incomplete Gamma function, and \eqref{Edef} provides the fundamental definition of $E_{s}(z)$. The result \eqref{Ei1} is given in \cite[Eq.(2.6a)]{Milgram:1985}. Here, the integration contour, originally defined to enclose the real axis $v\geq0$ as well as the singularity at $v=s-1$ in a clockwise direction, has been converted into the line $c<\Re(s)-1$, because the integrand vanishes as $v\rightarrow \pm\, i\infty$. This paper investigates the application of \eqref{Ei1} to \eqref{Leclair1}.

\subsection{A novel Integral Representation}

In the following, I focus on the limited range $0\leq \sigma \leq 1$, where $s=\sigma+i{\it t}$, in which case, applying \eqref{Ei1} to the first term in \eqref{Leclair1} yields

\begin{align}\label{EiT1}
\displaystyle \pi\, \left( s-1 \right) \sum _{n=1}^{\infty }{n}^{2}{E_{ -s/2}} \left(\pi\,{n}^{2} \right) =\frac{(s-1)}{2i} 
\int_{c_{{1}}-i\infty }^{c_{{1}}+i\infty }\!{\frac {\Gamma \left( -v \right) {\pi}^{v}\sum _{n=1}^{\infty }{n}^{2+2\,v}}{-s/2-1
\mbox{}-v}}\,{\rm d}v\\ \label{EiT2}
 =\frac{(s-1)}{2i}  
\mbox{}\int_{c_{{1}}-i\infty }^{c_{{1}}+i\infty }\!{\frac {\Gamma \left( -v \right) {\pi}^{v}\zeta \left( -2\,v-2 \right) }{-s/2
\mbox{}-1-v}}\,{\rm d}v
\end{align}
where the interchange of integration and summation in \eqref{EiT1} is justified if $c_1 < -3/2$ since the sum converges under this condition. The general requirement that $c_1<-\sigma/2-1$ imposes the effective constraint $c_1< -3/2$ for $0 \leq \sigma \leq 1$; the contour may be shifted leftwards with impunity since there are no singularities in that direction. Similarly, under the transformation $s\rightarrow 1-s$ we find the following expression for the second term of \eqref{Leclair1}
\begin{align}
\displaystyle -\pi\,s\sum _{n=1}^{\infty }{n}^{2}{E_{(s-1)/2}} \left(\pi\,{n}^{2} \right) =-\frac{s}{2i}\int_{c_{{2}}-i\infty }^{c_{{2}}+i\infty }\!{\frac {\Gamma \left( -v \right) {\pi}^{v}\zeta \left( -2\,v-2 \right) }{-3/2
\mbox{}+s/2-v}}\,{\rm d}v
\label{E1T3}
\end{align}
again valid for $c_2< -3/2$ if $\sigma>0$. In the case that $c_1=c_2=c$ where $c<-3/2$, and focussing on $0\leq\sigma\leq 1$, we find the integral equation
\begin{align} \nonumber
\displaystyle \xi \left( s \right) =&{\frac {i
\mbox{}s \left( s-1/2 \right) }{\pi}\int_{c-i\infty }^{c+i\infty }\!{\frac {\xi \left( v+3 \right) }{ \left( v+3 \right)  \left( 3-s+v \right)  \left( s+2+v \right) }}\,{\rm d}v}+{\frac {i}{2\pi}\int_{c-i\infty }^{c+i\infty }\!{\frac {\xi \left( v+3 \right) }{ \left( v+3 \right)  \left( s+2+v \right) }}\,{\rm d}v}
\\&+{\frac {{\pi^{1/4}}}{2\Gamma \left( 3/4 \right) }}
\label{EqXsiGen}
\end{align}
which, in terms of $\zeta(s)$ can be rewritten
\begin{align} \nonumber
(s-1)\zeta \left( s \right) \Gamma \left( s/2 +1\right) {\pi}^{-s/2}=2\,is \left( s-1/2 \right) \int_{c-i\infty }^{c+i\infty }\!{\frac {\Gamma \left( -v \right) {\pi}^{v}\zeta \left( -2\,v-2 \right) }{ \left( -3+s-2\,v \right)  \left( s+2+2\,v \right) 
}}\,{\rm d}v
\mbox{}\\ 
-i/2\int_{c-i\infty }^{c+i\infty }\!{\frac {\Gamma \left( -v \right) {\pi}^{v}\zeta \left( -2\,v-2 \right) }{s/2+1+v}}\,{\rm d}v+{\frac {\pi^{1/4}}{2\Gamma \left( 3/4 \right) }}\,.
\label{EqcGen}
\end{align}
Notice that with a simple application of the recursion formula for $\Gamma(-v)$, the numerator in \eqref{EqcGen} can be written in terms of $\Upsilon(-2v-2)$ - see \eqref{Lamdef}. Alternatively, under a simple change of variables, with the same condition $c<-3/2$ and the same range of $\sigma$, \eqref{EqcGen} becomes

\begin{align} 
\xi(s) ={\frac {\pi^{1/4}}{2\Gamma \left( 3/4 \right) }}-J(s,c)\,,
\label{EqcGen1}
\end{align}
where
\begin{equation}
\displaystyle J \left( s,c \right) =\int_{-\infty }^{\infty }\!{\pi}^{c-iv}\zeta \left( -2\,c+2\,iv-2 \right)\Gamma \left( -c+iv \right) M(c,s,v)\,{\rm d}v
\label{Jint}
\end{equation}
and
\begin{equation}
M(c,s,v)={\frac {
\mbox{} \left( - \left( 2\,s-1 \right) ^{2}/2+2\,iv-2\, \left( c+5/4 \right)  \right)  }{- \left( 2\,s-1 \right) ^{2}/4+2\,i \left( 2\,iv-4\,\left( c+5/4 \right)  \right) v
\mbox{}+4\, \left( c+5/4 \right) ^{2}}}\,.
\label{ModF}
\end{equation}

This result is worthy of a few comments:
\begin{itemize}
\item Its form is almost (exception: see \cite[Eq.(1.18)]{Patkowski})  unique among representations of $\zeta(s)$. Usually $s$-dependence that is formally embedded inside an integral or series representation of $\zeta(s)$ appears either as an exponent, or buried inside the argument of a transcendental function; here $s$-dependence exists only in the form of a coefficient in a simple rational (polynomial) function. This augers well for further analysis, and leads to some surprising predictions;

\item \eqref{EqXsiGen} and \eqref{EqcGen} present a prescription, in the form of an integral transform, for the value of $\xi(s)$ or $\zeta(s)$ anywhere in the complex strip $0\leq \Re(s)\leq1$, that depends only on its values on a vertical line in the complex plane corresponding to $c<-3/2$. The region in which \eqref{EqcGen1} is valid is labelled ``I" in Figure (\ref{fig:Figb}), and delineated as everything to the left of the point at $v=-1.5$ in Figure (\ref{fig:Figa});

\item In the form \eqref{EqcGen1} and \eqref{Jint} are written, the function $M(c,s,v)$ acts as a transfer function between $\zeta(v)$ on the line $\Re(v)=-2c-2$ and $\zeta(s)$ elsewhere in the complex plane, through the medium of the integral operator \eqref{Jint}. Its properties will be of interest in studying \eqref{EqcGen1} for varying values of the parameters $c$ and $s$.

\item Based on the property \eqref{Reflxi}, it is possible to show \cite[Theorem 1]{GlasserXi} that \eqref{EqcGen} can be cast into a simple closed (Cauchy) contour integral for $\xi(s)$.
\end{itemize}

In the following sections simple choices of $s$ and $c,c_1$, and $c_2$ will be applied to \eqref{EqcGen1}.

\begin{figure}[h] 
\centering
\subfloat [This figure shows the location of the contour of integration (black arrow at $\Re(v)<-3/2$ extending from $-i\infty$ to $i\infty$ as well as the various poles. The fixed pole at $v=-3/2$ corresponds to the singularity belonging to $\zeta(-2v-2)$. The green and red poles correspond to singularities of the integrand that depend on the value of $s$ acting as a parameter. The coloured arrows show the motion of corresponding poles as $\Re(s)$ increases from zero to one if $t>0$. The colours correspond to colours used in Figure \ref{fig:Figb}.]
{
\includegraphics[width=.40\textwidth]{{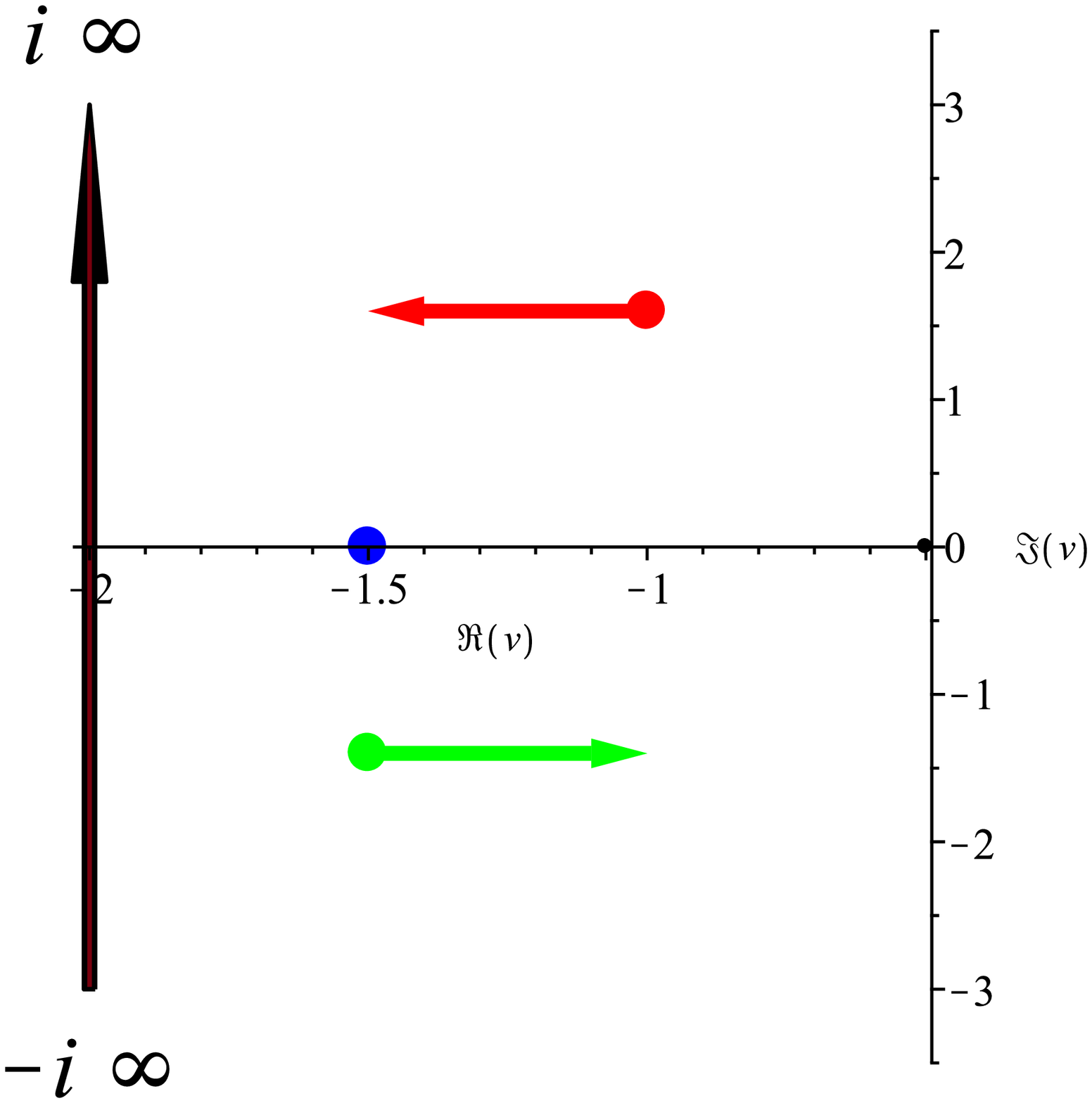}} \label{fig:Figa}
}
\hfill
\subfloat [This figure defines the various regions where different continuations of \eqref{EqcGen1} apply projected onto the ${\it t}=0$ plane. Five distinct regions are bounded by (coded and coloured) lines $c=-3/2$ (blue, solid), $c=-\sigma/2-1$ (red, dash) and $c=\sigma/2-3/2$ (green, dash-dot) when $0\leq \sigma\leq1$. The bounding regions extend vertically out of the plane of the figure when ${\it t}\neq0$ because the contour used in \eqref{Ei1} is chosen to be a straight line, vertical in the complex $v-$plane ($c$ is a constant). The dotted line corresponds to the interesting case $c=-5/4$.]
{
\includegraphics[width=.40\textwidth]{{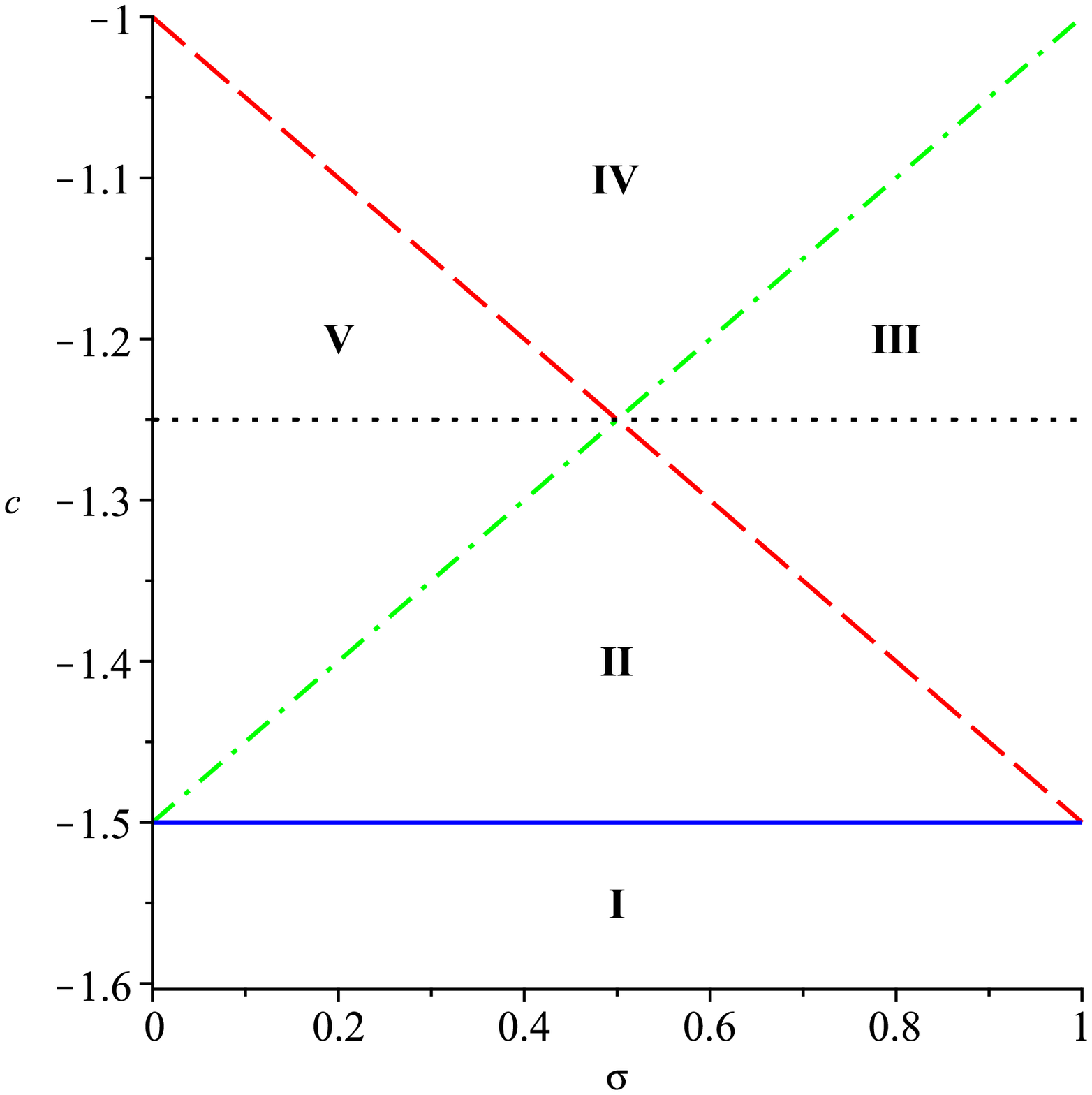}} \label{fig:Figb}
}
\caption{This Figure illustrates the various regions of validity of the analytic continuations developed in Section \ref{sec:AnalCont} both (a) in the complex $v-$plane and (b) in the cross section of a complex plane cut by the real ($\sigma, c)$ plane at ${\it t}=0$.}
\end{figure}

\section{Analytic Continuations} \label{sec:AnalCont}

In the following, let $c_1=c_2=c$. By shifting the contour \eqref{EqcGen} into the region $c>-3/2$, (see Figure (\ref{fig:Figa})) various representations are obtained for each of the regions labelled in Figure (\ref{fig:Figb}), by subtracting ($2 \pi i \times$) the residues  of the singularities that are transited in the complex $v-$plane. To enter the region labelled ``II" in Figure (\ref{fig:Figb}), the contour in Figure (\ref{fig:Figa}) must pass to the right of the fixed pole at $v=-3/2$, whose residue is given by
\begin{equation}
{\rm Residue_{blue}} = 1/(2\pi i)
\label{ResBlue}
\end{equation} 
giving, for region II,
\begin{align} 
\displaystyle \xi(s)={\frac {\pi^{1/4}}{2\,\Gamma \left( 3/4 \right) }}
-1-J(s,c)\,.
\label{Eqn2b}
\end{align}
To enter the region labelled III from region II, (crossing the red line in Figure (\ref{fig:Figb})) the corresponding residue (red in Figure (\ref{fig:Figa})) 
\begin{equation}
{\rm Residue_{red}} =\displaystyle i \xi(s)/(2\pi) 
\label{ResRed}
\end{equation}
must be subtracted. The result, valid for Region III is
\begin{equation} 
0={\frac {\pi^{1/4}}{2\,\Gamma \left( 3/4 \right) }}
-1-J(s,c)\,.
\label{Eqn3b}
\end{equation}
To cross the green line, either from Region II, entering Region V, or from Region III, entering Region IV, the residue
\begin{equation}
\displaystyle {\rm Residue_{green}}= \,{\frac {i s  {\pi}^{(s-3)/2}\zeta \left( 1-s \right) \Gamma \left( -s/2+1/2 \right) 
 \left(s-1 \right) }{{4} }}
\label{ResGreen}
\end{equation} corresponding to the green pole in Figure (\ref{fig:Figa}) must be subtracted. A similar argument applies for a transition from region V to region IV - the residue \eqref{ResRed} must be subtracted. This gives, for Region IV
\begin{align} 
\displaystyle -\xi(s)={\frac {\pi^{1/4}}{2\,\Gamma \left( 3/4 \right) }}
-1-J(s,c)\,.
\label{Eqn4b}
\end{align}
\newline

In particular, this analysis relates the results applicable to region IV to that for region I, including regions where both representations are unconstrained with respect to $0\leq\sigma\leq1$, and equates the result for region V to that for region III. It is emphasized that the regional definitions in Figure (\ref{fig:Figb}) exclude the lines, thereby defining the punctured regions ($s_1\neq s_2$) alluded to in the Introduction. Of particular interest, is the value at the point $(s=1/2, c=-5/4)$, as well as the values along the horizontal (dotted) line $s=\sigma+i{\it t} |_{c=-5/4}$ and a vertical line $c |_{s=1/2+i{\it t}}$. To reiterate, Figure (\ref{fig:Figb}) is intended to illustrate the projection of a (complex) space ($s=\sigma+i{\it t}$) defined by a third axis extending out of the plane of the Figure, corresponding to ${\it t}\neq0$, onto the (real) plane $(\sigma,c)$ where ${\it t}=0$.

\section{Crossing the line at c=-3/2} \label{sec:Integrals}

{\bf Remark:} Several of the identities in this Section are (redundantly) declared to be real and marked $\Re$. This is done as a reminder that unless so-declared, an attempted numerical computer evaluation of that expression will likely fail because the Imaginary part of the identity contains a singularity (e.g. \eqref{Ans0}).  

\subsection{The general case c=-3/2}\label{sec:C3half}

Disregarding the contour integral representations from which they are obtained, \eqref{EqcGen1}, \eqref{Eqn2b} and \eqref{Eqn3b} simply involve several complex functions of two real variables (if ${\it t}=0$), and ought to be amenable to analysis from that viewpoint. Compare \eqref{EqcGen1} and \eqref{Eqn2b} respectively below and above the horizontal line demarking regions I and II of Figure \ref{fig:Figb} by setting $c=-3/2-\eta$ in \eqref{EqcGen1} and $c=-3/2+\eta$ in \eqref{Eqn2b} where $0\leq \eta \approx 0$. Note that $\eta$ is a real variable (there is no cut); since none of the terms in either equation (other than $J(s,c)$) depend on $c$, we are free to choose a convenient value for $c$, subject to the (regional) conditions under which each equation is valid. We have in Region I
\begin{equation}
\xi(s)={\frac {\pi^{1/4}}{2\Gamma \left( 3/4 \right) }}-J(s,-3/2-\eta)
\label{Eqn1regionI}
\end{equation}
and in Region II
\begin{equation} 
\xi(s)={\frac {\pi^{1/4}}{2\,\Gamma \left( 3/4 \right) }}
-1-J(s,-3/2+\eta)\,.
\label{Eqn2regionII}
\end{equation}
Simple comparison between \eqref{Eqn1regionI} and \eqref{Eqn2regionII} suggests that 
\begin{equation}
J(s,-3/2-\eta)-J(s,-3/2+\eta)=1
\label{Eqdiffs}
\end{equation}
and this may be independently verified by a change of variables $v\rightarrow \eta v$ in both integrals, followed by a series expansion about $\eta=0$, leading to
\begin{equation}
J(s,-3/2-\eta)-J(s,-3/2+\eta)\rightarrow \frac{1}{\pi^2}\int_{-\infty}^{\infty} \frac{1}{1+v^2} dv =1\,.
\label{Jint1m2}
\end{equation}

As is often done where discontinuities arise in the theory of functions of a real variable, when $\eta=0$ it is conventional to assign half the discontinuity at that point ($c=-3/2$). This is equivalent to deforming the contour in \eqref{EqcGen} to avoid the fixed pole belonging to $\zeta(-2v-2)$ at $v=-3/2$ by including only half the residue at that point.\newline

{\bf Remark:} This cannot be done for the other poles since they are defined by $s$ (complex), rather then $v=-3/2$ (fixed).

\subsubsection{The case $\sigma=1/2$}

The result from \eqref{Eqn1regionI} and \eqref{Eqn2regionII} after setting $\eta=0$ (see Section \ref{sec:C3half}) when $s=1/2$ is 
\begin{align} \nonumber
\frac{1}{2\pi}\int_{-\infty }^{\infty }\!\Im \left( {\frac {\xi(1+2iv) }
{(2\,iv+1/2)\,v}} \right) \,{\rm d}v&=
\int_{-\infty }^{\infty }\!\Re \left( {\frac {{\pi}^{-3/2-iv}
\zeta \left( 1+2\,iv \right) \Gamma \left( 3/2+iv \right) }
{2\,iv+1/2}} \right) \,{\rm d}v\\ 
&={\frac {\zeta \left( 1/2 \right) {\pi}^{3/4} \sqrt{2}}{8\,\Gamma \left( 3/4 \right) }}+{\frac {\pi^{1/4}}{2\,\Gamma \left( 3/4 \right) }}
\mbox{}-1/2\,.
\label{Ans0}
\end{align} 
This result (compare with \cite[page 204, Lemma $\beta$]{Titch2}) can be verified numerically. The need to specify that only the real part of the integrand in \eqref{Ans0} is to be used, is twofold:
\begin{itemize}
\item
the right-hand side is real and so must be the left-hand side;
\item
ostensibly, the integrand of \eqref{Ans0} appears to be singular at $v=0$ unless one notes that 
\begin{equation}
\zeta(1\pm 2iv)\approx \mp i/(2v)+\gamma+\dots
\label{NearZ}
\end{equation}
and $\gamma$ is Euler's constant. Therefore, the real part of the integrand converges at $v=0$ and the imaginary part, singular at the origin, integrates to zero by anti-symmetry.
\end{itemize}

In general, when $s=1/2+i{\it t}$, the left-hand side of both \eqref{Eqn1regionI} and \eqref{Eqn2regionII} are real, and so the equivalent of \eqref{Ans0} applies. Written in terms of ${\it t}$, \eqref{Ans0} generalizes to

\begin{align} \nonumber
&\xi(1/2+i{\it t})\\ 
&=-{\frac {1}{{2\,\pi}^{3/2}}\int_{-\infty }^{\infty }\!\Re \left({\frac {{\pi}^{-iv}\zeta \left( 1+2\,iv \right) 
\mbox{} \left( 4\,{{\it t}}^{2}+4\,iv+1 \right) \Gamma \left( 3/2+iv \right) 
\mbox{}}{{{\it t}}^{2}+2\,iv-4\,{v}^{2}+1/4}}
 \right) \,{\rm d}v}
-\frac{1}{2}+{\frac {{\pi}^{1/4}}{2\,\Gamma \left( 3/4 \right) }}\,.
\label{AsGenh}
\end{align}

\subsection{The case c=-5/4}
In regions III and V, writing \eqref{Eqn3b} in full gives
\begin{equation}
\displaystyle {\pi}^{c}\int_{-\infty }^{\infty }\!{\frac {{\pi}^{-iv}\zeta \left( -2\,c+2\,iv-2 \right) 
\mbox{}\Gamma \left( -c+iv \right)  \left( 2\,iv-2\,{s}^{2}-2\,c+2\,s-3 \right) }{ \left( -2\,c+2\,iv+s-3 \right) 
\mbox{} \left( -2\,c+2\,iv-s-2 \right) }}\,{\rm d}v={\frac {{\pi}^{1/4}}{2\,\Gamma \left( 3/4 \right) }}-1
\label{GenC}
\end{equation}
or, equivalently
\begin{equation}
\displaystyle \int_{-\infty }^{\infty }\!{\frac { \left( 2\,iv-2\,{s}^{2}-2\,c+2\,s-3 \right)
\mbox{}\xi \left( -2\,c+2\,iv-2 \right) }{ \left( -2\,c+2\,iv+s-3 \right) 
\mbox{} \left( -2\,c+2\,iv-s-2 \right)  \left( -2\,c+2\,iv-3 \right) 
\mbox{}}}\,{\rm d}v=\,{\frac {{\pi}^{1/4}}{2\,\Gamma \left( 3/4 \right) }}-1
\label{GenCXsi}
\end{equation}

Setting $c=-5/4$ in \eqref{GenC} leads to

\begin{equation}
\displaystyle \int_{-\infty }^{\infty }\!{\frac {\zeta \left( 1/2+2\,iv \right) 
\mbox{}{\pi}^{-5/4-iv}\Gamma \left( 5/4+iv \right) 
 \left( iv-\left( 2\,s-1 \right) ^{2}/4 \right) }{ \left( 2\,s-1 \right) ^{2}/8+2\,{v}^{2}}}\,{\rm d}v=1-{\frac {{\pi}^{1/4}}{2\,\Gamma \left( 3/4 \right) }}\,,
\label{cm54}
\end{equation}
a result that is valid for all $s$ (see the dotted line in Figure (\ref{fig:Figb})). Since the left-hand side is a function of $s$ and the right-hand side is not, it must be true that the derivative of the left-hand side with respect to $s$ vanishes. Performing this calculation (the integral is convergent, so the derivative and integral operator can be interchanged) yields

\begin{equation}
\displaystyle \left( s-1/2 \right) \int_{-\infty }^{\infty }\!{\frac {{\pi}^{-5/4-iv} \left( -2-iv+ \left( 2\,s-1 \right) ^{2}/\,4 \right) 
\mbox{}\zeta \left( 1/2+2\,iv \right) \Gamma \left( 5/4+iv \right) 
\mbox{}}{ \left(  \left( 2\,s-1 \right) ^{2}/\,8+2\,{v}^{2} \right) ^{2}}}\,{\rm d}v=0
\label{Deq1}
\end{equation}

Compare \eqref{cm54} with Patkowski \cite{Patkowski}, Theorem 1 .

\subsubsection{$\sigma=1/2$}

Consider \eqref{cm54} in the case $s=1/2$, which gives
\begin{equation}
\displaystyle \int_{-\infty }^{\infty }\!{\frac {i\,{\pi}^{-5/4-iv} }{2\,v}\zeta \left( 1/2+2\,iv \right) 
\Gamma \left( 5/4+iv \right)}\,{\rm d}v=1-{\frac {{\pi}^{1/4}}{2\,\Gamma \left( 3/4 \right) }}\,,
\label{C54&s_half}
\end{equation}
and the integrand apparently has a singularity at $v=0$. However, a simple expansion of the integrand about that point shows that
\begin{equation}
\displaystyle \Re\left(   {\frac {i\,{\pi}^{-5/4-iv}\,}{2v}\zeta \left( 1/2+2\,iv \right) 
\Gamma \left( 5/4+iv \right) } \right) = -\,{\frac {
\zeta \left( 1/2 \right) }{\sqrt{2}\,{\pi}^{1/4}\,\Gamma \left( 3/4 \right) }}+O(v^2)\,,
\label{ReShalf}
\end{equation}

demonstrating that the real part of the integrand in \eqref{C54&s_half} is non-singular at $v=0$. However, near $v=0$ it is similarly shown that the imaginary part of the integrand diverges like $v^{-1}$ and the imaginary part of the integral vanishes by anti-symmetry about $v=0$. \newline

In the case $s=1/2+i{\it t}$, \eqref{cm54} can be written in terms of $\xi(1/2+2iv)$ as
\begin{equation}
\displaystyle {\frac {1}{\pi}\int_{-\infty }^{\infty }\!{\frac {\xi \left( 1/2+2\,iv \right) 
\mbox{} \left( iv+{{\it t}}^{2} \right) }{ \left( -{{\it t}}^{2}/2+2\,{v}^{2} \right)  \left( -1/2+2\,iv \right) 
\mbox{}}}\,{\rm d}v}=1-{\frac {{\pi}^{1/4}}{2\Gamma \left( 3/4 \right) }}
\label{Eq5Ax}
\end{equation}
whose integrand appears to become singular at $v=\pm{\it t}/2$. Writing the integrand in terms of its real and imaginary parts, and noting that $\xi(1/2+2iv)$ is real, we find, for the real part
\begin{equation}
\displaystyle{\frac {4}{\pi
\mbox{}}\int_{-\infty }^{\infty }\!{\frac {\xi \left( 1/2+2\,iv \right) }{16\,{v}^{2}+1}}\,{\rm d}v}=1-{\frac {{\pi
}^{1/4}}{2\,\Gamma \left( 3/4 \right) }}
\label{Eq5AxR}
\end{equation}
and,  for the imaginary part,
\begin{equation}
\displaystyle \int_{-\infty }^{\infty }\!{\frac {v\xi \left( 1/2+2\,iv \right) }{ \left( 4\,{v}^{2}-{{\it t}}^{2} \right) 
\mbox{} \left( 16\,{v}^{2}+1 \right) }}\,{\rm d}v=0\,.
\label{Eq5AxI}
\end{equation}
Notice that in the former case there is no ${\it t}$ dependence, so in \eqref{Eq5Ax}, the variable ${\it t}$ is in reality a free parameter, consistent with the argument applied to obtain \eqref{Deq1}. In the latter case although a singularity exists, the integrand is anti-symmetric about $v=0$ and so the singularities at $\pm v$ cancel.

\section{An asymptotic approximation to $\zeta(\sigma+i{\it t}$)} \label{sec:Asympt}

\subsection{Integral Equations for the real and imaginary parts of $\xi(s)$} \label{sec:ZsApprox}
Consider \eqref{Eqn4b} - region IV - in the general case $s=\sigma+i{\it t}$, written in terms of the transfer function $M(c,s=\sigma+i{\it t},v)$. The case $c=-1$ is of interest since the integral will span the 0-line, which connects to the 1-line by reflection (see \eqref{FeqId}), and, acting in its capacity as a master function, $\zeta(1\pm i{\it t})$ has many well-known properties (e.g. \cite{doi:10.1093/imrn/rnx331}), the most relevant one being that it has no complex zeros. In addition, according to Figure (\ref{fig:Figb}), \eqref{Eqn4b} is valid for all values of $\sigma$ spanning the critical strip $0<\sigma<1$. In this case the transfer function $M_x(-1,\sigma+i{\it t},v)$ and its real and imaginary components, shortened, except where necessary, to $M_R({\it t},v)$ and $M_{I}({\it t},v)$, are respectively:


\begin{align} \label{Mx1}
&\displaystyle {\it M_x}(-1,\sigma+i{\it t},v)\,= \,-{\frac { \left( 2\,\sigma+2\,i{\it t}-1 \right)  \left( \sigma+i{\it t} \right) 
\mbox{}}{ \left( 2\,iv-1+\sigma+i{\it t} \right)  \left( 2\,iv-\sigma-i{\it t} \right) 
\mbox{}}}+ \left( 2\,iv-\sigma-i{\it t} \right) ^{-1}\\ \label{Ms1}
&\displaystyle {\it M_R({\it t},v)} = -\!{\frac {8{v}^{2} \left( {{\it t}}^{2}-{\sigma}^{2}+\sigma \right)\! + 2{\it t} v\left(1 -2\sigma \right)+2\sigma^3(2-\sigma)- \left( 4{{\it t}}^{2}+3 \right) {\sigma}^{2}
+ \left( 4{{\it t}}^{2}+1 \right) \sigma-2\,{{\it t}}^{4}-{{\it t}}^{2}}{ \left( {{\it t}}^{2}+4\,v{\it t}+{\sigma}^{2}+4\,{v}^{2}-2\,\sigma+1 \right)  \left( {{\it t}}^{2}-4\,v{\it t}+{\sigma}^{2}+4\,{v}^{2} \right) 
\mbox{}}}\\ \label{Mc1}
&\displaystyle {\it M_{I}}({\it t},v)= \,-{\frac {8\,{v}^{3}+ \left( -16\,\sigma+8 \right) {\it t}\,{v}^{2}+ \left( -6\,{{\it t}}^{2}+6\,{\sigma}^{2}-6\,\sigma+2 \right) v
\mbox{}+ \left( 2\,\sigma-1 \right) {\it t}}{ \left( {{\it t}}^{2}+4\,v{\it t}+{\sigma}^{2}+4\,{v}^{2}-2\,\sigma+1 \right)  \left( {{\it t}}^{2}-4\,v{\it t}+{\sigma}^{2}+4\,{v}^{2} \right) }}\,.
\end{align}

Furthermore, these functions possess the following symmetry properties: 
\begin{align} \label{McSymm}
M_R(\sigma+i{\it t}, v) =&\;\;\:\; M_R(1-\sigma+i{\it t}, -v)\\ \label{MsSymm}
M_{I}(\sigma+i{\it t}, v) =& -M_{I}(1-\sigma+i{\it t}, -v).
\end{align}
which symmetry also holds true for all values of $c$. A complete description of these functions is given in Appendices \ref{sec:PropsOfMI} and \ref{sec:PropsOfMR} . From \eqref{Eqn4b}, the basic equation in Region IV using $c=-1$ is

\begin{equation}
\displaystyle 1-{\frac {{\pi}^{1/4}}{2\Gamma \left( 3/4 \right) }}-\xi \left( \sigma+i{\it t} \right) =-\frac{i}{\pi}{ {\int_{-\infty }^{\infty }\! \left( M_{{R}} \left( {\it t},v \right) +iM_{{I}} \left( {\it t},v \right)  \right) 
\mbox{}\Upsilon \left( 2\,iv \right)\, v\,{\rm d}v
\mbox{}}{}}\,.
\label{Beq0}
\end{equation}
It is now convenient to rewrite \eqref{Beq0} as an integral over the range $[0\dots\infty)$ by first splitting the integral, setting $v\rightarrow -v$ and combining the two halves. The expression obtained is straightforward, although rather lengthy. A second lengthy equation can be obtained by first replacing $\sigma\rightarrow 1-\sigma$ in \eqref{Beq0} and similarly reducing the range to $[0\dots\infty]$. After splitting $\Upsilon(2iv)$ into its real and imaginary components, these two equations can be added and subtracted, and with the help of  \eqref{McSymm} and \eqref{MsSymm} two new fundamental equations emerge:
\begin{align}
\displaystyle {\frac {{\pi}^{1/4}}{2\Gamma \left( 3/4 \right) }}-1=-\xi_R \left( \sigma+i{\it t} \right) 
+\frac{1}{\pi}{ {\int_{0}^{\infty }\! \left( T_{1}(s,v)\,  \Upsilon_R \left( 2\,iv \right)  + (T_{2}(s,v)-4)\, \Upsilon_{I}\left( 2\,iv \right)  
\right) v\,{\rm d}v}{}}
\label{EqpR}
\end{align}
and


\begin{equation}
\displaystyle  \xi_I \left(s \right)  =\frac{1}{\pi}\,{ {\int_{0}^{\infty }\! \left( T_{{3}} \left( s,v \right) \Upsilon_{{R}} \left( 2\,iv \right) +{ T_4} \left( s,v \right) \Upsilon_{{I}} \left( 2\,iv \right) 
\mbox{} \right) v\,{\rm d}v}{}}\,,
\label{EqmI}
\end{equation}
where
\begin{align}
\displaystyle T_{{1}} \left( s,v \right)\equiv& -M_{{I}} \left( s,v \right) +M_{{I}} \left(s,-v \right) \label{T1Def}
\\
\displaystyle T_{{2}} \left(s,v \right) \equiv& -M_{{R}} \left(s,v \right) -M_{{R}} \left( s,-v \right) 
\mbox{}+4 \label{T2Def}
\\
\displaystyle T_{{3}} \left( s,v \right) \equiv&\hspace{12pt} M_{{R}} \left( s,v \right) -M_{{R}} \left( s,-v \right) \label{T3Def}
\\
\displaystyle T_{{4}} \left( s,v \right) \equiv& -M_{{I}} \left(s,v \right) -M_{{I}} \left( s,-v \right)\,. \label{T4Def}
\end{align}
In this manner, because $\xi(s)$ is self-conjugate, the real and imaginary components of $\xi(s)$ are isolated, and expressed in terms of convergent integrals. Furthermore, subtracting \eqref{Mom1} from \eqref{EqpR} yields a simpler form  

\begin{equation}
\displaystyle \xi_R \left( s \right)   =\frac{1}{\pi} {\int_{0}^{\infty }\! \left( T_{{1}}(s,v)\Upsilon_{{R}} \left( 2\,iv \right) +T_{{2}}(s,v)\Upsilon_{{I}} \left( 2\,iv \right) 
\mbox{} \right) 
\mbox{}v\,{\rm d}v}{}
\label{EqpR2}
\end{equation}

For the case ${\it t}=0,0<\sigma<1$, \eqref{EqpR2} reduces to a relation between $\xi(\sigma)$ on a section of the real line, and $\Upsilon(2iv)$ on the complex line $\Re(v)=0$, specifically
\begin{align} \nonumber
\displaystyle \xi_R \left( \sigma\right) =\frac{1}{\pi}\int_{0}^{\infty }\!\frac {1}{ \left( {\sigma}^{2}+4\,{v}^{2} \right)  \left(  \left( \sigma-1 \right) ^{2}+4\,{v}^{2} \right) } & \left( 4\,{v}^{2} \left( 3\,\sigma\, \left( \sigma-1 \right) +4\,{v}^{2}+1 \right) \Upsilon_R(2iv) \right.\\
&\left.  \hspace{-3cm} +v \left( 64\,{v}^{4}+16\, \left( {\sigma}^{2}-\sigma+1 \right) {v}^{2}-2\,\sigma\, \left( \sigma-1 \right)  \right) \Upsilon_I(2iv) 
\mbox{} \right) \,{\rm d}v\,.
\label{Eq00}
\end{align}

\subsection{The background terms} \label{sec:Backg}

Although the functions $T_{1,2,3,4}(s,v)$ look formidably complicated when written in full, each can be written in a more transparent form when decomposed using partial fractions. Specifically, with $c=-1$:


\begin{equation}
\displaystyle T_{{1}} \left( s,v \right) ={\frac {2\,\sigma\,v-t}{ \left( t-2\,v \right) ^{2}+ \left( \sigma-1 \right) ^{2}}}+{\frac {2\, \left( 1-\sigma \right) v-t}{ \left( t-2\,v \right) ^{2}+{\sigma}^{2}}}
\mbox{}+{\frac {2\,\sigma\,v+t}{ \left( t+2\,v \right) ^{2}+ \left( \sigma-1 \right) ^{2}}}+{\frac {2\, \left( 1-\sigma \right) v+t}{ \left( t+2\,v \right) ^{2}+{\sigma}^{2}}}
\label{T1pa}
\end{equation}
\begin{align} \nonumber
\displaystyle T_2(s,v)= &{\frac {-t \left( t-2\,v \right) + \left( 1-\sigma \right) \sigma}{ \left( t-2\,v \right) ^{2}+{\sigma}^{2}}}+{\frac {-t \left( t-2\,v \right) + \left( 1-\sigma \right) \sigma}{ \left( t-2\,v \right) ^{2}+ \left( 1-\sigma \right) ^{2}}}
\mbox{}+{\frac {-t \left( t+2\,v \right) + \left( 1-\sigma \right) \sigma}{ \left( t+2\,v \right) ^{2}+{\sigma}^{2}}}\\
&+{\frac {-t \left( t+2\,v \right) + \left( 1-\sigma \right) \sigma}{ \left( t+2\,v \right) ^{2}+ \left( 1-\sigma \right) ^{2}}}+4
\label{T2pa}
\end{align}

\begin{align} \nonumber
\displaystyle T_{{3}} \left( s,v \right) =& {\frac {-t \left( t-2\,v \right) -\sigma\, \left( \sigma-1 \right) }{ \left( t-2\,v \right) ^{2}+ \left( \sigma-1 \right) ^{2}}}+{\frac {t \left( t-2\,v \right) +\sigma\, \left( \sigma-1 \right) }{ \left( t-2\,v \right) ^{2}+{\sigma}^{2}}}
\mbox{}+{\frac {t \left( t+2\,v \right) +\sigma\, \left( \sigma-1 \right) }{ \left( t+2\,v \right) ^{2}+ \left( \sigma-1 \right) ^{2}}}\\&{\frac {-t \left( t+2\,v \right) -\sigma\, \left( \sigma-1 \right) }{ \left( t+2\,v \right) ^{2}+{\sigma}^{2}}}
\label{T3pa}
\end{align}
\begin{align}
\displaystyle T_{4}(s,v)=&{\frac {-2\,v \left( \sigma-1 \right) -t}{ \left( t-2\,v \right) ^{2}+{\sigma}^{2}}}+{\frac {-2\,\sigma\,v+t}{ \left( t-2\,v \right) ^{2}+ \left( \sigma-1 \right) ^{2}}}+{\frac {2\,\sigma\,v+t}{ \left( t+2\,v \right) ^{2}+ \left( \sigma-1 \right) ^{2}}}
+{\frac {2\,v \left( \sigma-1 \right) -t}{ \left( t+2\,v \right) ^{2}+{\sigma}^{2}}}
\label{T4pa}
\end{align}
{\bf Remark:} Although not evident when written in the decomposed form, when expressed in factored form, both \eqref{T3pa} and \eqref{T4pa} contain an overall factor $(\sigma-1/2)$, a requirement that $\xi_I(s)$ vanishes when $\sigma=1/2$. Further, it is easily seen that $T_{1,2}(\sigma+it,v)=T_{1,2}(1-\sigma+it,v)$, $T_{3,4}(\sigma+it,v)=-T_{3,4}(1-\sigma+it,v)$, $T_{1,2}(\sigma+it,v)=T_{1,2}(\sigma-it,v)$ and $T_{3,4}(\sigma+it,v)=-T_{3,4}(\sigma-it,v)$.\newline

As written, the first two terms of each of the above explicitly demonstrate the existence of a pole in the complex $v$-plane at $v=(\pm\, t/2\pm i\sigma)$, along with invariance under the symmetry $\sigma\rightarrow 1-\sigma$. Since we are in general interested in the case $t\ggg \{\sigma, 1-\sigma\} >0$, such poles lie relatively close to the positive real $v-$axis as demonstrated by the nature of the peak(s) in Figure \ref{fig:FigMsOver} - see Appendix A. By reducing the integral to the range $[0\dots\infty)$, the influence of poles corresponding to negative values of $v$ has been effectively removed from the integrand ($v>0$). Drawing on these observations, it is suggestive that the functions $T_{1,2,3,4}(s,v)$ be separated into two components - pole terms associated with the first two terms, plus background terms associated with all four terms, in each of the above; it is further suggestive that only the pole terms will contribute to the integrals (see Figures \ref{fig:FigMcOver} and \ref{fig:P1234}), thereby allowing an approximate solution to be obtained. How can this be done?\newline

Consider \eqref{EqpR2}. From \eqref{T1pa} and \eqref{T2pa}, it would appear superficially that both $T_1(s,v)$ and $T_2(s,v)$ vanish as $1/t$ for large values of $t$, which immediately demonstrates a potential inconsistency - the left-hand side of \eqref{EqpR2} vanishes exponentially and so must the right-hand side. Since the remaining integrand factors $\Upsilon_R(2iv)$ and $\Upsilon_I(2iv)$ contain no $t$ dependence, the only possibility is that $1/t$ dependence must vanish from the (highly oscillatory) integral, and somewhere buried in the higher order asymptotic terms will be found a term with exponentially decreasing $t$ dependence. Further, the possibility also exists that cancellations will occur between the two integrand terms containing $T_1(s,v)$ and $T_2(s,v)$ in \eqref{EqpR2}. Any attempt at a numerical evaluation of \eqref{EqpR2} for reasonably large values of $t$ will immediately demonstrate the truth of this prediction in the form of a severe cancellation of significant digits.  An example follows:\newline

Choose $\sigma=1/3$ and $t=50$, chosen sufficiently large to distinguish inverse powers of $t$ from exponentially decreasing terms of the form $\exp(-bt)$, but not so large that multi-digit arithmetic must be used. For this case, calculated with 15 digits, we find
\begin{equation}
\displaystyle  \frac{1}{\pi}{ {\int_{0}^{\infty }\!vT_{{1}} \left( 1/3+50\,i,v \right) 
\mbox{} \Upsilon_R \left( 2\,iv \right)  \,{\rm d}v}}= 0.000810442386190651
\label{Ex1a}
\end{equation}

\begin{equation}
\displaystyle \frac{1}{\pi}{{\int_{0}^{\infty }\!vT_{{2}} \left( 1/3+50\,i,v \right) 
\mbox{} \Upsilon_I \left(2\,iv \right) \,{\rm d}v}{}}=-0.0008104423861872670
\label{Ex1b}
\end{equation}
and the sum $3.3840\times10^{-15}$ is comparable to $\xi(1/3+50i) =3.38361\times10^{-15}$. Clearly, each of the integrals has an absolute value much greater than the sum of the two, and it is only because of a cancellation of digits that the final result can be found. Here, we see that the cancellation of digits has resulted in the loss of 10 digits from the sum; for larger values of $t$, the effect will be more significant, and one might despair of utilizing \eqref{EqpR2} for anything. It will now be shown that these (and other) cancellations can be analysed analytically, and an exponentially decreasing term can be extracted, which is the basis for the approximate model solution.\newline

First of all, for $v\leq t/2$, consider the asymptotic expansion ($t\rightarrow\infty$) of each of the individual terms composing $T_1$, labelled respectively by a second subscript
%
%
\begin{align} \label{Ta_1234p}
\displaystyle {T_{1,1}(s,v)}\,  \rightarrow & \,-\frac{1}{t}+{\frac {2\,\sigma\,v-4\,v}{{t}^{2}}}+{\frac {4\, \left( 2\,\sigma-3 \right) {v}^{2}+ \left( \sigma-1 \right) ^{2}}{{t}^{3}}}+O \left( {t}^{-4} \right) \\ \label{Tb_1234p}
\displaystyle {T_{1,2}(s,v)}\, \rightarrow & \,-\frac{1}{t}+{\frac {2\, \left( 1-\sigma \right) v-4\,v}{{t}^{2}}}+{\frac {-4\, \left( 2\,\sigma+1 \right) {v}^{2}+{\sigma}^{2}}{{t}^{3}}}+O \left( {t}^{-4} \right)\\ \label{Tc_1234p}
\displaystyle {T_{1,3}(s,v)}\,  \rightarrow & \,\:\frac{1}{t}+{\frac {2\, \left( 1-\sigma \right) v-4\,v}{{t}^{2}}}+{\frac {4\, \left( 2\,\sigma+1 \right) {v}^{2}-{\sigma}^{2}}{{t}^{3}}}+O \left( {t}^{-4} \right)\\  \label{Td_1234p} 
\displaystyle {T_{1,4}(s,v)}\, \rightarrow &\,\;\frac{1}{t}+{\frac {2\,\sigma\,v-4\,v}{{t}^{2}}}+{\frac {-4\, \left( 2\,\sigma-3 \right) {v}^{2}- \left( \sigma-1 \right) ^{2}}{{t}^{3}}}+O \left( {t}^{-4} \right)
\end{align}

Clearly, the leading asymptotic dependence ($1/t$) cancels between all of ``pole" and ``background" terms in pairs, demonstrating that the pole terms also contribute significantly to the asymptotic background. This suggests that a numerical evaluation of the $T_1(s,v)$ integral will be challenging as $t$ increases. More importantly, it is apparent that the overall asymptotic dependence of the $T_1(s,v)$ part of the integral will have the leading dependence $1/t^2$ because the multiplicative factor $\Upsilon_R(2iv)$ lacks $t$ dependence as noted. But, as also noted, this is inconsistent with \eqref{EqpR2}, where the left-hand side has an exponentially decreasing asymptotic behaviour in $t$. One explanation for this inconsistency could be that there exists a cancellation between the two terms in the integrand labelled by $T_1$ and $T_2$ to at least order $1/t^2$. This would explain the cancellation observed in the numerical results \eqref{Ex1a} and \eqref{Ex1b}. The following gives the asymptotic expansion of each of the corresponding elements of $T_2(s,v)$ in analogy to \eqref{Ta_1234p} - \eqref{Td_1234p}:

\begin{align}
\displaystyle {T_{2,1}(s,v)}\, \rightarrow & \,-1-2\,{\frac {v}{t}}+{\frac {-4\,{v}^{2}+\sigma}{{t}^{2}}}+{\frac {-8\,{v}^{3}+2\,\sigma\, \left( \sigma+2 \right) v}{{t}^{3}}}+O \left( {t}^{-4} \right) \\
\displaystyle {T_{2,2}(s,v)}\,  \rightarrow & \,-1-2\,{\frac {v}{t}}+{\frac {-4\,{v}^{2}-\sigma+1}{{t}^{2}}}+{\frac {-8\,{v}^{3}+ \left( 2\,{\sigma}^{2}-8\,\sigma+6 \right) v}{{t}^{3}}}+O \left( {t}^{-4} \right) \\
\displaystyle {T_{2,3}(s,v)}\, \rightarrow & \,-1+2\,{\frac {v}{t}}+{\frac {-4\,{v}^{2}+\sigma}{{t}^{2}}}+{\frac {8\,{v}^{3}-2\,\sigma\, \left( \sigma+2 \right) v}{{t}^{3}}}+O \left( {t}^{-4} \right)\\
\displaystyle {T_{2,4}(s,v)}\, \rightarrow & \,-1+2\,{\frac {v}{t}}+{\frac {-4\,{v}^{2}-\sigma+1}{{t}^{2}}}+{\frac {8\,{v}^{3}+ \left( -2\,{\sigma}^{2}+8\,\sigma-6 \right) v}{{t}^{3}}}+O \left( {t}^{-4} \right) 
\label{T2_1234}
\end{align}
Again, we see a cancellation between each of the four terms to order $1/t$ as well as a cancellation of the terms of order $t^0$ with the fifth term in \eqref{T2pa}. And again we find an inconsistency with \eqref{EqpR2} because the leading term of $T_2(s,v)$ is also of order $1/t^2$. Define each term corresponding to $t^{-k}$ in the asymptotic series of $T_1(s,v)$  by
\begin{equation}
\tilde{T_1}(k)\equiv \text{coefficient of } t^{-2k}\;\text{ in  the series } \underset{\underset{v<t}{t\rightarrow\infty}} {\lim}T_{1}(s,v) 
\label{Ttilde}  
\end{equation}
and similarly for $T_2(s,v)$. When the first two asymptotic terms of the factors in the integrand on the right-hand side of \eqref{EqpR2} are written in full we find

\begin{equation} 
\displaystyle {\it T_1({\it s},v)}\, \sim \,-{\frac {12\,v}{{{\it t}}^{2}}}+{\frac {4\,v \left( 9\,{\sigma}^{2}-20\,{v}^{2}-9\,\sigma+4 \right) }{{{\it t}}^{4}}}+O \left( {{\it t}}^{-6} \right) 
\label{T1Asy}
\end{equation}
and
\begin{equation}
\displaystyle {\it T_2({\it s},v)}\, \sim \,{\frac {-16\,{v}^{2}+2}{{{\it t}}^{2}}}+{\frac {-64\,{v}^{4}+ \left( 48\,{\sigma}^{2}-48\,\sigma+48 \right) {v}^{2}-6\,{\sigma}^{2}+6\,\sigma-2}{{{\it t}}^{4}}}
\mbox{}+O \left( {{\it t}}^{-6} \right). 
\label{T2Asy}
\end{equation}
so that
\begin{align}
\tilde{T_1}(1)&=-{\frac {12\,v}{{{\it t}}^{2}}}\\
\tilde{T_1}(2)&={\frac {4\,v \left( 9\,{\sigma}^{2}-20\,{v}^{2}-9\,\sigma+4 \right) }{{{\it t}}^{4}}}\\
\tilde{T_2}(1)&={\frac {-16\,{v}^{2}+2}{{{\it t}}^{2}}}\\
\tilde{T_2}(2)&={\frac {-64\,{v}^{4}+ \left( 48\,{\sigma}^{2}-48\,\sigma+48 \right) {v}^{2}-6\,{\sigma}^{2}+6\,\sigma-2}{{{\it t}}^{4}}}
\end{align}

and clearly there can be no cancellation of terms of equal order in $1/t^2$, or $1/t^4$ between terms containing $T_1(s,v)$ and $T_2(s,v)$. The only remaining possibility is that the integral containing each of these asymptotic forms itself vanishes {\bf identically}, at least to the orders discussed above. In fact, this must be true to all asymptotic orders in $t^{-2k},\; k\geq 0$ because it is impossible for a series of the form $\sum_{k=1}^{\infty}a_k {\it t}^{-2k}$ to add up to an exponential of the form $\exp(-bt)$. Furthermore, if the integrals vanish as predicted, they must vanish for all $\sigma$. Otherwise, \eqref{EqpR2}, which is an exact result cannot be asymptotically and universally true. And somewhere, accompanying all these dominant terms that vanish asymptotically to inverse polynomial order, exist other terms associated with the pole terms, that carry the exponentially vanishing result corresponding to the left hand-side of \eqref{EqpR2} in the asymptotic limit $t\rightarrow \infty$. Numerically, these lesser terms equate to the difference between \eqref{Ex1a} and \eqref{Ex1b}. It will now be proven that the integrals do in fact vanish, as predicted, to at least the first few orders in $t^{-2k}, \, k=1,2,3$, following which a model will be constructed to extract the terms that instead vanish exponentially.\newline


Based on the the first terms of \eqref{T1Asy} and \eqref{T2Asy}, it is predicted that the sum of the integrals involving ${t}^{-2}$ must vanish, that is
\begin{equation}
\displaystyle{\frac {1}{{{\it t}}^{2}}\int_{0}^{\infty } -12\,{v}^{2}\Upsilon_R(2iv) 
 +v \left( -16\,{v}^{2}+2 \right)\Upsilon_I(2iv)\,{\rm d}v}=0
\label{EqRhom2}
\end{equation}

{\bf Proof:} Consider the classical result \cite[Eq.(2.15.6)]{Titch2} relating the inverse Mellin transform of $\Upsilon(c+iv)$ to the Jacobi $\Theta_3$ function:

\begin{equation} 
\displaystyle \int_{-\infty }^{\infty }\!\Upsilon(c+iv)\,{\rm d}v=4\,\pi\,\sum _{n=1}^{\infty }{{\rm e}^{-{n}^{2}x}} = 2\pi(\Theta_3(0,\exp(-\pi x))-1)
\label{R3}
\end{equation}
valid for $c>1$. We are interested in the case $c=0$, so shift the contour by subtracting the residue at $v=(1-c)/i$ and half the residue at $v=ic$ resulting in the identity
\begin{align} \nonumber
\displaystyle \pi\, \left( 2\,\sum _{n=1}^{\infty }{{\rm e}^{-{n}^{2}x}}-{\frac { \sqrt{\pi}}{ \sqrt{x}}}+\frac{1}{2} \right)& =\int_{0}^{\infty }\! \left( {\pi}^{iv}{x}^{-iv}+{\pi}^{-iv}{x}^{iv} \right) \Upsilon_R \left( 2\,iv \right)
\mbox{}\,{\rm d}v\\
&+i\int_{0}^{\infty }\!\left( {\pi}^{iv}{x}^{-iv}-{\pi}^{-iv}{x}^{iv} \right) {\Upsilon_I} \left( 2\,iv \right) 
\mbox{}\,{\rm d}v
\label{Rans0}
\end{align} 
if $c=0$, after converting the integration limits to $(0,\infty)$ and writing $\Upsilon(2iv)$ as the sum of its real and imaginary parts. By an increasing sequence of higher order derivatives, evaluated at $x=\pi$, it is possible to obtain relevant sums using identities given in Romik's paper \cite{Romik}. The first seven such sums are listed in Appendix \ref{sec:RomSums}, from which it is possible to obtain the first seven even moments of $\Upsilon_R(2iv)$ and odd moments of $\Upsilon_I(2iv)$, listed also in Appendix \ref{sec:RomSums}. Substituting the moments \eqref{Mom1}...\eqref{Mom3} into \eqref{EqRhom2} verifies that the prediction \eqref{EqRhom2} is true. QED \newline 

Employing similar logic, a further prediction arises regarding the next term of order $1/t^4$. It too must vanish for all $\sigma$. From \eqref{T1Asy}, \eqref{T2Asy} and \eqref{EqRhom2} we predict

\begin{equation}
\displaystyle \,\frac{1}{{\it t}^{4}}\int_0^{\infty}\left( \Upsilon_{{R}} \left( 2\,iv \right) {v}^{4}+{\frac {v \left( 96\,{v}^{4}-40\,{v}^{2}-1 \right) \Upsilon_{{I}} \left( 2\,iv \right) }{120}}\right)\,dv=0\,.
\label{EqRhom4}
\end{equation}
Note that $\sigma$ dependence has disappeared from \eqref{EqRhom4}, as expected. By substituting \eqref{Mom1}...\eqref{Mom5} into \eqref{EqRhom4} the truth of this prediction can be verified. This sequence of forecasts can be continued to the next recursive level, producing the forecast identity
\begin{equation}
\displaystyle \frac{1}{{\it t}^{6}}\int_{0}^{\infty }\left( 224\,\Upsilon_{{R}} \left( 2\,iv \right) {v}^{6}
- \left( -128\,{v}^{7}+112\,{v}^{5}+{\frac {28\,{v}^{3}}{3}}+v/3 \right)\Upsilon_{{I}} \left( 2\,iv \right)\right)\,{\rm d}v=0\,.
\label{EqRhom6}
\end{equation}
Again, $\sigma$ dependence has cancelled of its own volition. The truth of the prediction \eqref{EqRhom6} can similarly be verified by substituting the set of moments \eqref{Mom1}...\eqref{Mom7} into \eqref{EqRhom6}.\newline

To summarize, we have\newline

{\bf Proposition 1}(uproven but partially demonstrated): The above sequence of predictions can be proven in general for all terms on the right-hand side of \eqref{EqpR2} that asymptotically vanish to inverse polynomial order $t^{-2k}, k>3$ independently of $\sigma$.\newline 

Whatever remains must then be of exponentially decreasing order and must therefore equate to the left-hand side of \eqref{EqpR2}, and the challenge is to somehow isolate such term(s). The cancellation of integrals of asymptotically vanishing inverse polynomial order explains the loss of significant digits observed between \eqref{Ex1a} and \eqref{Ex1b} and suggests that such loss of significant digits will only increase for larger values of $t$, as the interested reader can verify for herself. \newline

Expressed another way, the crucial point is that although the peaks shown in Figure \ref{fig:FigMsOver} originate from the pole terms, these pole terms also contribute to the background terms of inverse polynomial order ($t^{-2k},\; k>1$) that cancel. Nonetheless, the signal associated with the pole terms naturally separates from the background . The separation of signal from background is the essence of the approximate solution to be presented here.\newline

To gain further insight into these issues, consider Figure \ref{fig:J1fig} which shows (the absolute value of) the integrand 

\begin{figure}[h] 
\centering 
\subfloat  [An illustration of the integrand $|J_1(1/3+it,v)|$ at several values of $t=\{12.5,25,50,100,200\}$; the arrows mark the points $v=t/2$, that being the location of the peak of the pole term contained in $T_1(s,v)$ for each of the curves at each individual value of $t$.]
{
\includegraphics[width=.44\textwidth]{{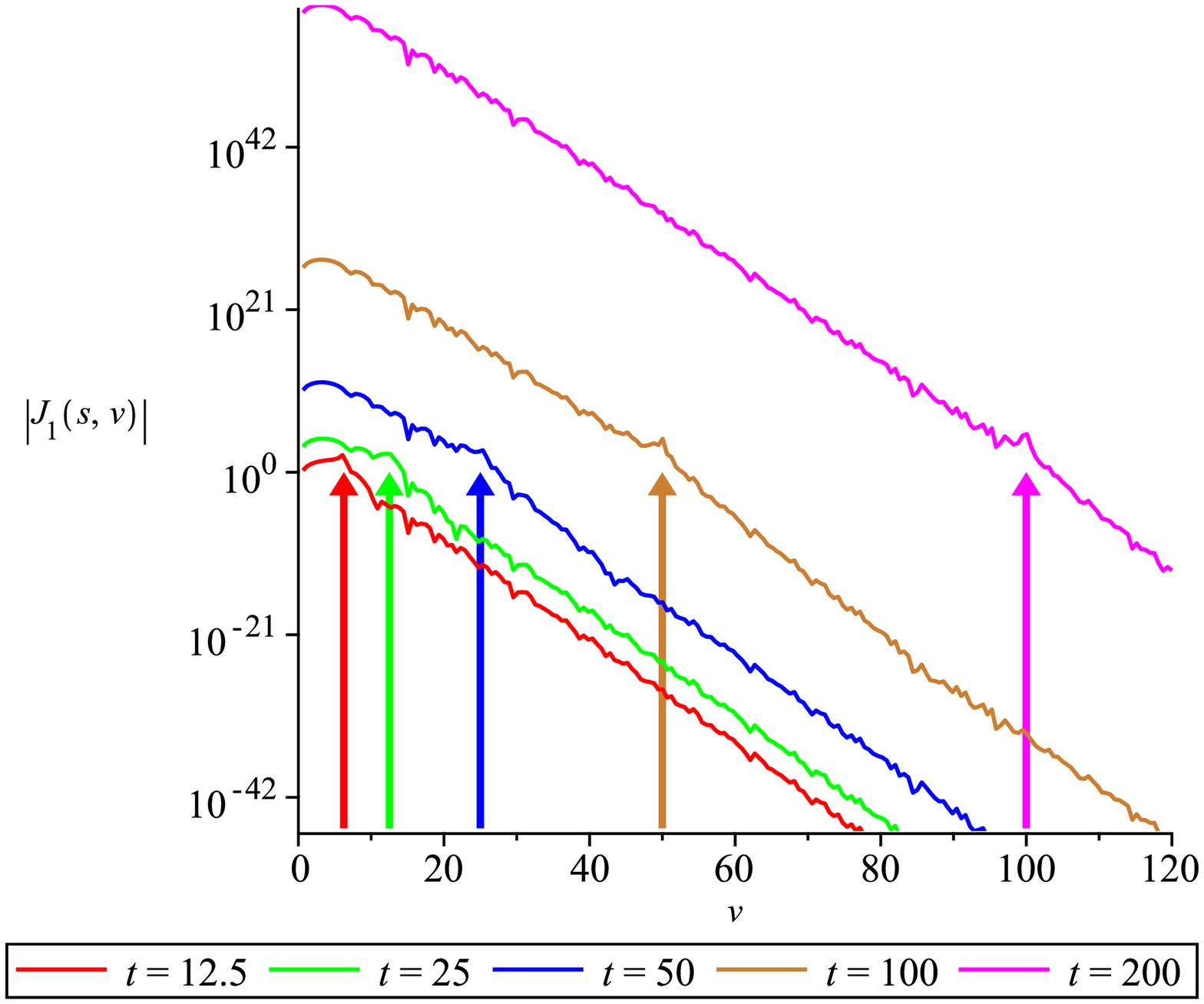}} \label{fig:J1fig} 
}
\hfill                                 
\subfloat [An illustration of the integrand $J_2(1/3+it,v)$ at several values of $t=\{12.5,25,50,100,200\}$; the arrows mark the points $v=t/2$, that being the location of the zero of the pole term contained in $T_2(s,v)$ for each of the curves at each individual value of $t$.]
{
\includegraphics[width=.44\textwidth]{{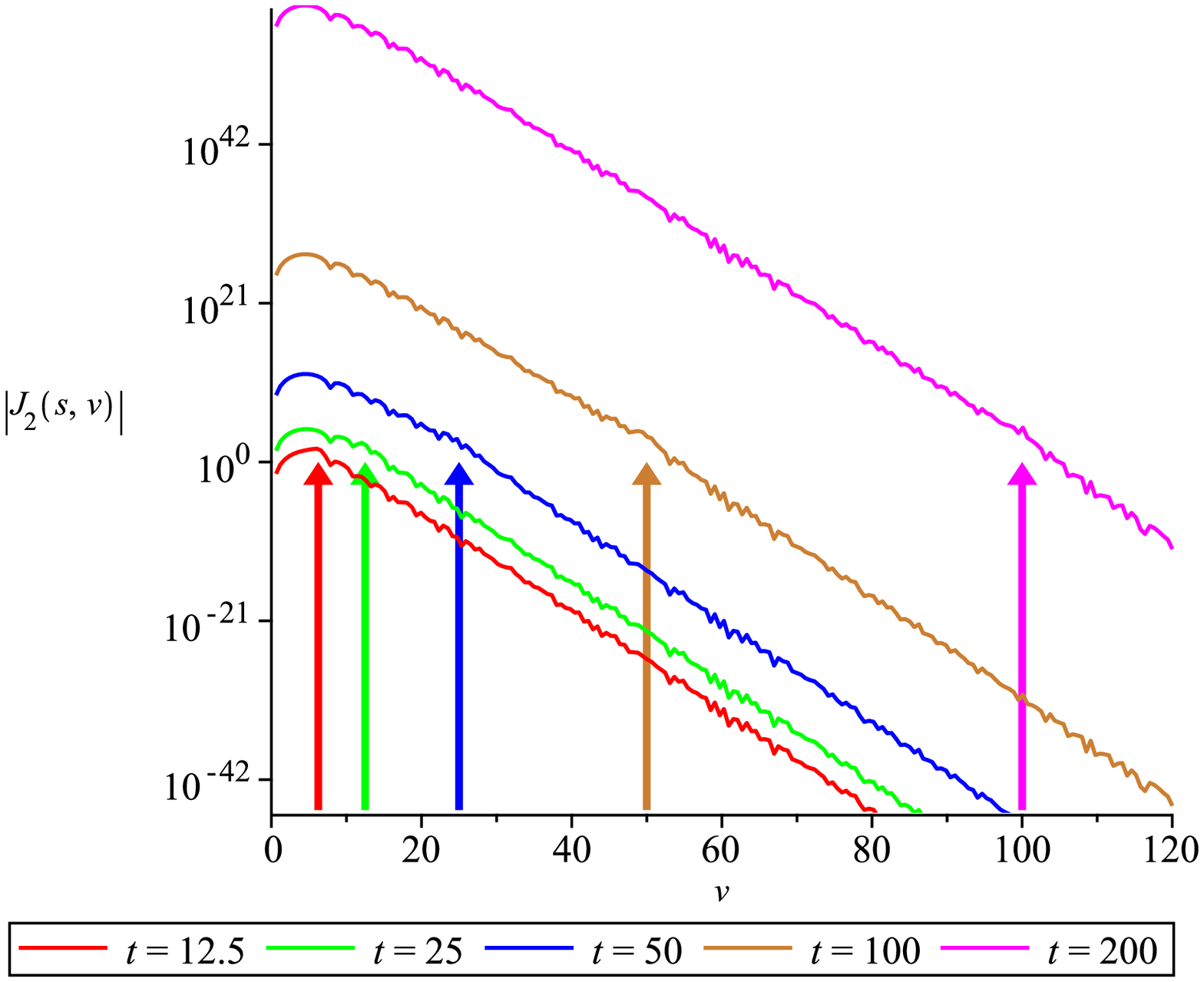}} \label{fig:J2fig} 
}
\caption{This figure illustrates the individual terms composing the integrand of \eqref{EqpR2} as a function of the integration variable $v$, using $\sigma=1/3$ for several values of $t$. In this Figure, each curve has been scaled by a (smooth) factor $\exp(\pi t/4)/t^2$ to separate them for easier viewing - otherwise all of the curves coincide to a good degree of approximation, except near $v\approx 0$.}
\label{fig:J1J2}
\end{figure}

\begin{equation}
J_{1}(s=\sigma+it,v)\equiv T_{1}(s,v) \Upsilon_{{R}} \left( 2\,iv \right)
\label{J1eq}
\end{equation}
as a function of $v$ for several values of $t$ at $\sigma=1/3$. Similarly, Figure \ref{fig:J2fig} shows the same result for the second term in the integrand defined by

\begin{equation}
J_{2}(s=\sigma+it,v)\equiv \,T_{2}(s,v)  \Upsilon_{{I}} \left( 2\,iv \right)\,.
\label{J2eq}
\end{equation} 

Three properties are immediately evident: most of the absolute value of both integrands occurs near $v=0$, both $|J_1(s,v)|$ and $|J_2(s,v)|$ decrease exponentially with increasing $v$ and are of similar magnitude - a necessity if they are to cancel numerically (see \eqref{Ex1a} and \eqref{Ex1b}). A fourth significant property is that the pole term, embedded in $T_1(s,v)$ at $v=t/2$ is clearly visible above the exponentially decreasing background in each curve of Figure \ref{fig:J1fig}, although it lies many orders of magnitude below the numerically significant portion of the integrand. In fact, it appears to be enhanced as $t$ increases, as might be expected because the magnitude of its peak varies as $t$ - see \eqref{T1Min} - whereas the background varies as $t^{-2k},\,k>1$. However, in Figure \ref{fig:J2fig} the zero of $T_2(s,v)$ is not readily apparent because of the resolution of each curve; moreover, each of the curves in Figure \ref{fig:J2fig} is consistent with what appears to be the background as discussed previously. It will now be shown that these observations can be used to extract a signal from the background, and in so doing obtain an approximate solution to \eqref{EqpR2}.

\subsection{The signal} \label{sec:signal}  

Consider Figure \ref{fig:T1abs} which illustrates the peak of the function $|T_1(s,v)|$ over a wide range of $v$ using $s=1/3+10^4i$. Also shown are the cumulative background terms $\sum_{k=1}^N|\tilde{T_1}(k)|$. It can be seen that as N increases, the sum of the background terms give a better and better approximation to the background part of the function $|T_1(s,v)|$ when $v<t/2$. In fact, when the each of the cumulative background terms are subtracted from $|T_1(s,v)|$, that part of the function associated with small values of $v$ decreases (Figure \ref{fig:T1Diff}) as the number of background terms increases. But it is precisely small values of $v$ that contribute most to the integral itself (see Figure \ref{fig:J1J2}). And as we have seen from \eqref{EqRhom2}, \eqref{EqRhom4} and \eqref{EqRhom6}, when coupled with the full part of the integrand associated with $T_2(s,v)$ the integral associated with each of these background terms has either been proven, or predicted, to vanish. What is left? Only the pole term remains, wherein resides the sought after ``signal".\newline

{\bf Remark}: Because absolute values have been used, the corresponding Figures for the function $T_2(s,v)$ are qualitatively identical to Figure \ref{fig:T1absDiff}. 

\begin{figure}[h] 
\centering 
\subfloat  [An illustration of the function $|T_1(1/3+10^4i,v)|$ as well as the cumulative first few background terms $\sum_{k=1}^N|\tilde{T_1}(k)|$ for the first 5 values of $N$ for a wide range of $v$ about $t/2$.]
{
\includegraphics[width=.46\textwidth]{{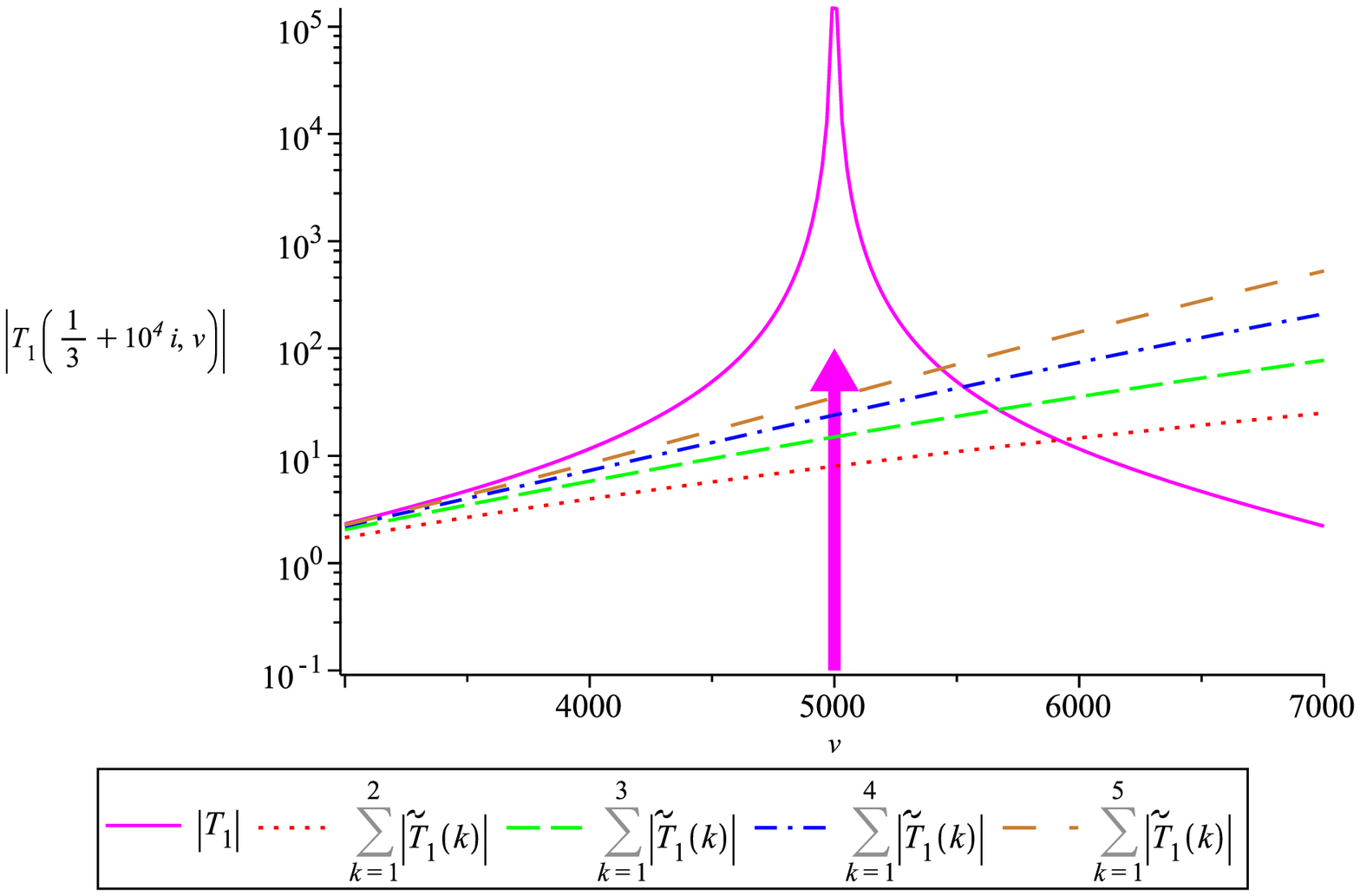}} \label{fig:T1abs} 
}
\hfill                                 
\subfloat [An illustration of the function $|T_1(1/3+10^4i,v)|$ as well as the effect of subtracting the first few cumulative background terms $\sum_{k=1}^N|\tilde{T_1}(k)|$ for the first 5 values of $N$ for a wide range of $v$ about $t/2$.]
{
\includegraphics[width=.46\textwidth]{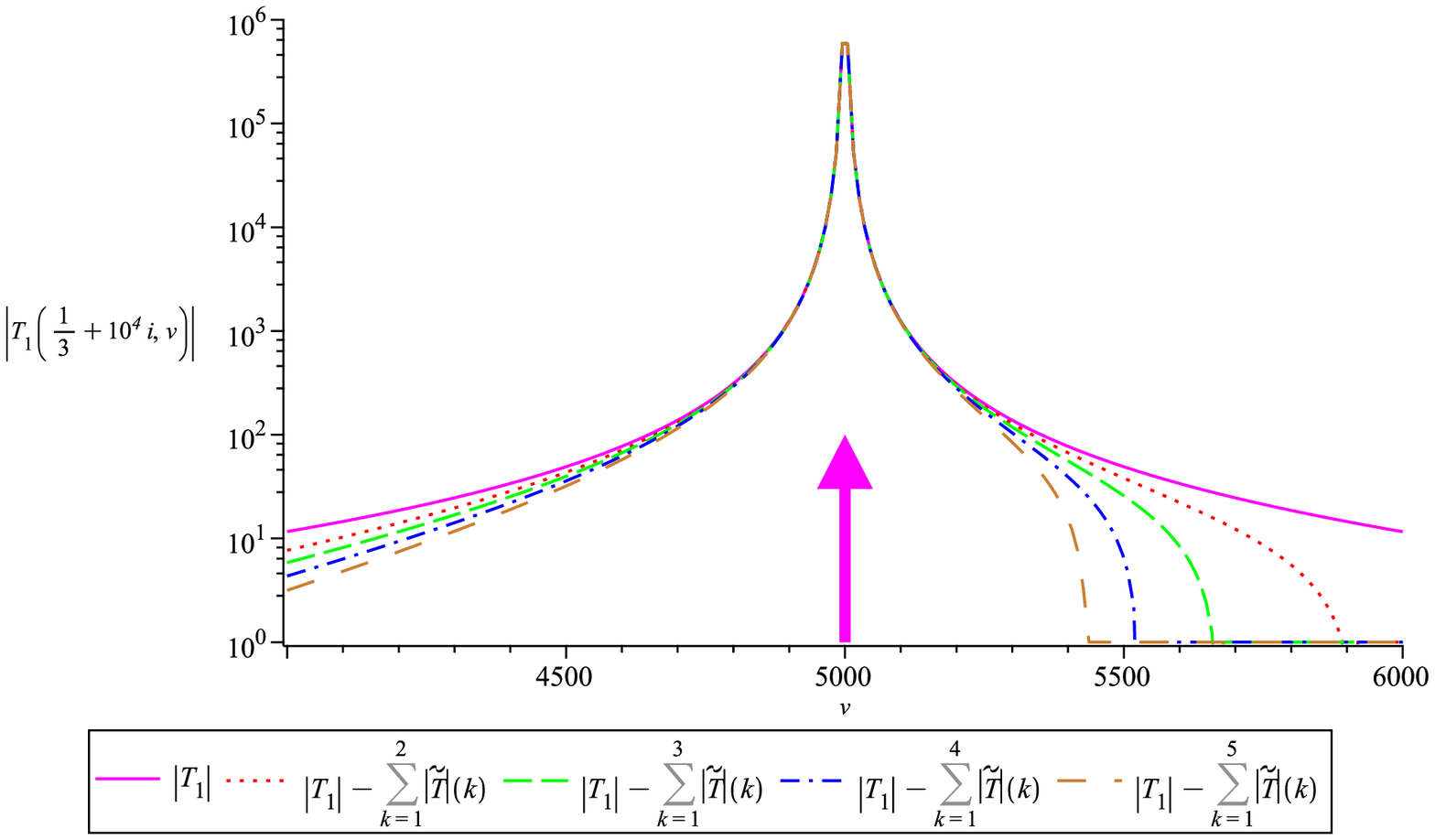} \label{fig:T1Diff} 
}
\caption{This figure illustrates the function $|T_1(1/3+10^4i,v)|$ in comparison to the first few cumulative background terms $\sum_{k=1}^N|\tilde{T_1}(k)|$ for  wide range of $v$ near $t/2$. In both Figures, the vertical arrow marks the location of the peak at $v=t/2$.}
\label{fig:T1absDiff}
\end{figure}\mbox{}\newline

To understand how all this affects the integral itself, consider Figure \ref{fig:J1J2All} which shows both integrands $J_1(s,v)$ and $J_2(s,v)$, again using $s=1/3+it$ with $t=10^4$, in close proximity to the peak region $v=t/2$. Since a logarithmic scale is used, in both Figures the alternating function(s) $J_{1|2}(s,v)$ are shown (dotted) colour coded, positive (red) and negative (blue). The green (dash-dot) line shows the leading asymptotic behaviour of $J_{1|2}(s,v)$ originating from the factor $\Gamma(iv)$ embedded in the factor $\Upsilon(2iv)$, scaled to fit the Figure. This is the background.\newline

Focussing on Figure \ref{fig:J1All}, the peak of the factor $T_1(s,v)$ is indicated by a magenta (dashed) line, also marked by the vertical arrow. The effect of the peak of $T_1(s,v)$ on the integrand $J_1(s,v)$ can be seen as an enhancement of the integrand in the immediate vicinity of the peak. This is indicated by the superimposed solid curve representing $J_1(s,v)$, again colour coded red and blue, but only for those values of $t/2-2<v<t/2+2$. The top horizontal dashed line equates to the ordinate of the point $J_1(s,t/2)$ and the lower line delineates points lying one order of magnitude smaller that the ordinates at the peak. Figure \ref{fig:J2All} is similar to Figure \ref{fig:J1All}, with the exception that instead of a peak at $v=t/2$, the curve representing $J_2(s,v)$ vanishes - see Appendix \ref{sec:Tprops}. 

\begin{figure}[h] 
\centering 
\subfloat  [An illustration of the function $J_1(1/3+10^4i,v)$ accompanied by $T_1(1/3+10^4i,v)$ and the scaled asymptotic shape of $\Re(\Upsilon_1(1/3+10^4i,v))$ very close to $v=t/2$.]
{
\includegraphics[width=.46\textwidth]{{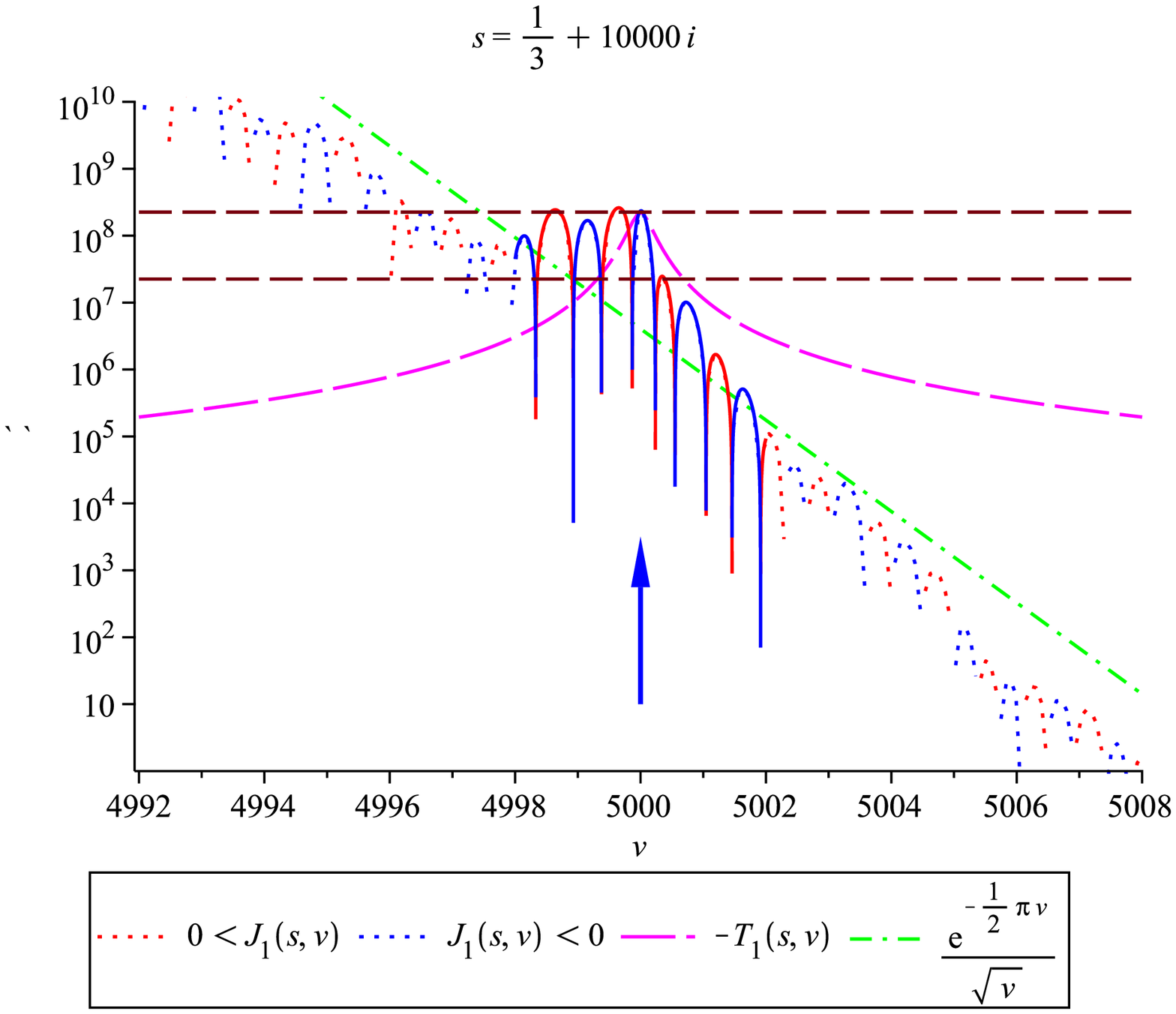}} \label{fig:J1All} 
}
\hfill                                 
\subfloat [An illustration of the function $J_2(1/3+10^4i,v)$ accompanied by $T_2(1/3+10^4i,v)$ and the scaled asymptotic shape of $\Re(\Upsilon_1(1/3+10^4i,v))$ very close to $v=t/2$.]
{
\includegraphics[width=.46\textwidth]{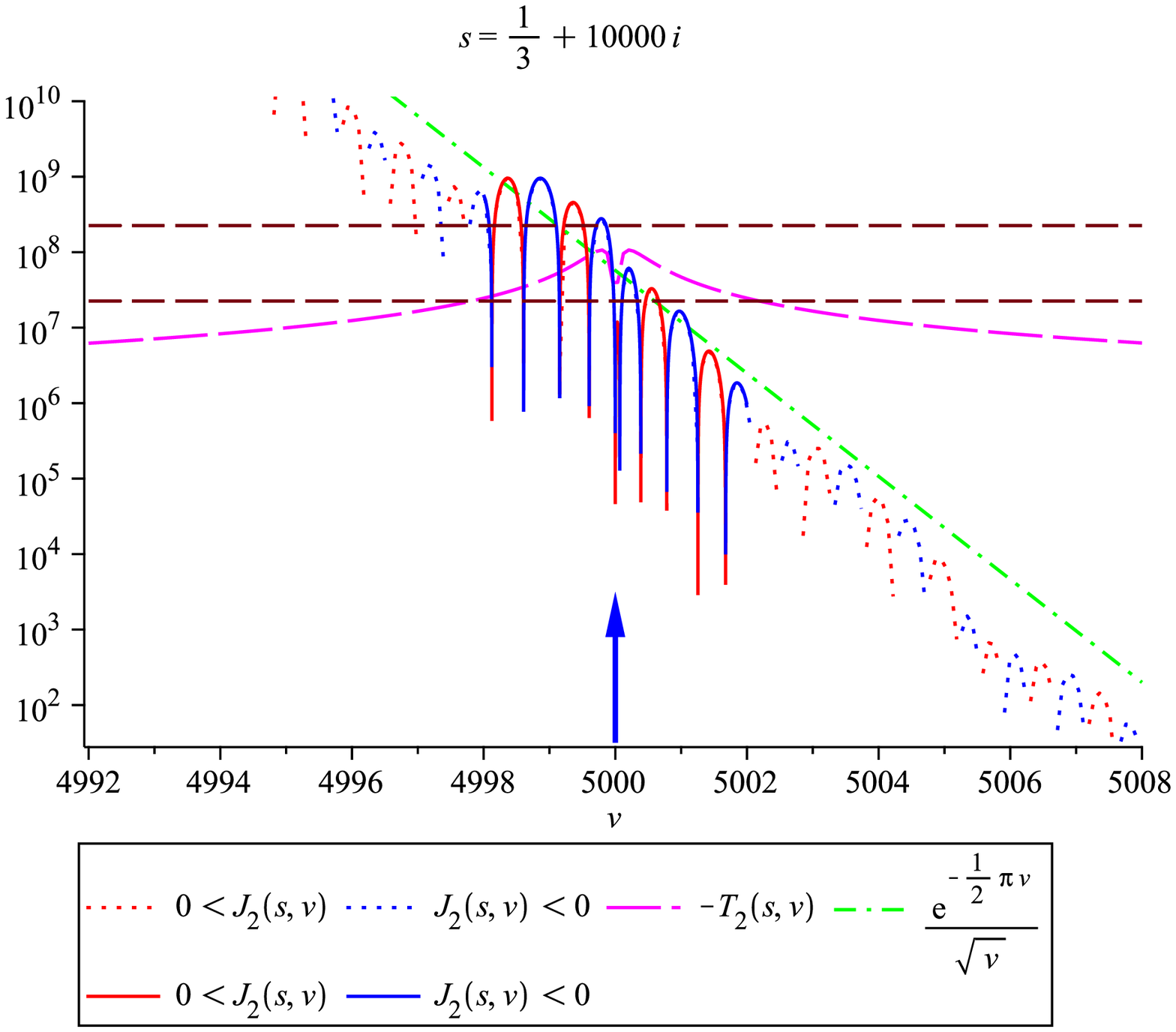} \label{fig:J2All} 
}
\caption{These Figures show $J_1(s,v)$ and $J_2(s,v)$ in close proximity to the peak of $T_1(s,v)$ or zero of $T_2(s,v)$ respectively at $v=t/2$, as well as the shape of the background of both.}
\label{fig:J1J2All}
\end{figure}

\subsection{The Approximate Solution} \label{sec:Approx}

It is Figures \ref{fig:J1All} and \ref{fig:J2All} that hold the key to the approximate solution of \eqref{EqpR2}. Up to this point, all analysis has consisted of observation, (generalized) deduction and a small amount of proof. The crucial point to be extracted from Figure \ref{fig:J1All} is that the peak of $T_1(s,v)$ ``extrudes" the integrand of $J_1(s,v)$ above background in a small region about $v=t/2$. Although the background (dotted part of the curve) represents the vast majority of the magnitude of the integrand over many orders of magnitude, when this background is decomposed into terms that are of inverse polynomial order (asymptotically in $t$), the integration of these terms cancels between the two terms labelled by $T_1(s,v)$ and $T_2(s,v)$ and illustrated in Figures \ref{fig:J1All} and \ref{fig:J2All} - see Section \ref{sec:Backg} and only that part of the curve near $v=t/2$ (shown solid) contributes to the integral. Conversely, with reference to Figure \ref{fig:J2All},  the strong zero of the function $T_2(s,v)$ acts to depress the function $J_2(s,v)$ near $v=t/2$, although the background terms still assume their role of acting to cancel similar background terms shown in Figure \ref{fig:J1All}. All this suggest an approximation that should give a reasonable representation of a solution to \eqref{EqpR2}.\newline

{\bf Approximation 1}\newline

Let
\begin{equation}
T_{j}(s,v)=\tilde{T_j}(s,v)+\overset{o}T_j(s,v)
\label{Approx1}
\end{equation}    
where $\tilde{T_j}(s,v)$ represents the background, $\,\overset{o}T_j(s,v)$ represents the non-background associated with the pole and $j=1,2$. From \eqref{Ttilde} we have

\begin{equation}
\tilde{T_j}(s,v)=\sum_{k=1}^{\infty} \tilde{T}_j(k).
\label{Backgrd}
\end{equation}
When applied to \eqref{EqpR2} this leads to

\begin{align}
J(s,v)\equiv &J_1(s,v)+J_2(s,v)=\left(\tilde{T_1}(s,v)+\overset{o}T_1(s,v)\right)\Upsilon_R(s,v)+\left(\tilde{T_2}(s,v)+\overset{o}T_2(s,v)\right)\Upsilon_I(s,v)\\
=&\left(\tilde{T_1}(s,v)\Upsilon_R(s,v)+\tilde{T_2}(s,v)\Upsilon_I(s,v)\right)+\left(\overset{o}T_1(s,v)\Upsilon_R(s,v)+\overset{o}T_2(s,v)\Upsilon_I(s,v)\right).
\label{ApproxSum}
\end{align}

But, according to Section \ref{sec:Backg}
\begin{equation}
\int_{0}^{\infty }\left(\tilde{T_1}(s,v)\Upsilon_R(s,v)+\tilde{T_2}(s,v)\Upsilon_I(s,v)\right)v\,{\rm{d}}v=0\,;
\label{BackgrdEq0}
\end{equation} 
thus \eqref{EqpR2} becomes
\begin{equation}
\xi_R(s)=\frac{1}{\pi}\int_{0}^{\infty }\overset{o}{T_1}(s,v)\Upsilon_R(s,v)v\,{\rm{d}}v+\frac{1}{\pi}\int_{0}^{\infty }\overset{o}{T_2}(s,v)\Upsilon_I(s,v)v\,{\rm d}v\,.
\label{NewEq}
\end{equation}

Since $\overset{o}{T_j}(s,v)$ represent the non-background (pole) terms this leads to\newline

{\bf Approximation 2}

\begin{equation}
\xi_R(s)=\frac{1}{\pi}\int_{t/2-\delta_1(\sigma,t)}^{t/2+\delta_1(\sigma,t) }\overset{o}{T}_1(s,v)\Upsilon_R(s,v)v\,{\rm{d}}v+\frac{1}{\pi}\int_{t/2-\delta_2(\sigma,t)}^{t/2+\delta_2(\sigma,t)}\overset{o}{T_2}(s,v)\Upsilon_I(s,v)v\,{\rm d}v\,\,,
\label{Approx2}
\end{equation}
where $\delta_{j}(\sigma,t)$ is a (small) function of $O(\sigma^0), O(t^0)$ associated with each of $\overset{o}{T}_j(s,v)$ - see \eqref{V2MinMax}; that is, only that part of the integrand in close proximity to $v=t/2$ contributes to the overall value of the integral, which integral is approximately symmetric about the point $v=t/2$. The remainder cancels. Equivalently, this recognizes that only that section of the curve coloured {\it solid} in Figures \ref{fig:J1All} and \ref{fig:J2All} will contribute to the integral, notwithstanding the fact that the dotted sections ostensibly contribute to the integrand by many more orders of magnitude. This is demonstrated numerically below - see \eqref{NumApprox}.\newline

{\bf Remark:} For clarity, I write $\delta_j(\sigma,t)\equiv \delta_j$ recognizing that $\delta_j(\sigma,t)$ is a simple and low-order function of both $\sigma$ and $t$; the functional dependence is investigated later.\newline

{\bf Approximation 3}\newline

Since each of the $\delta_j$ are much less than any value of $t$ in any range of interest, first rewrite \eqref{Approx2} as
\begin{align} \nonumber
\xi_R(s)=&\displaystyle  \frac{\delta_1}{2\,\pi}\,\int_{-1}^{1}\!\overset{o}{T}_{{1}} \left( s,t/2+\delta_1\,y \right) \Upsilon_{{R}} \left( i \left( 2\,\delta_{{1}}y+t \right)  \right)  \left( 2\,\delta_{{1}}y+t \right) \,{\rm d}y\\
+ &\displaystyle  \frac{\delta_2}{2\,\pi}\,\int_{-1}^{1}\!\overset{o}{T}_{{2}} \left( s,t/2+\delta_2\,y \right) \Upsilon_{{I}} \left( i \left( 2\,\delta_{{2}}y+t \right)  \right)  \left( 2\,\delta_{{2}}y+t \right) \,{\rm d}y
\label{Approx3a}
\end{align}
by applying an obvious change of variables, expand \eqref{Approx3a} about $\delta_j=0$, and consider only those terms of leading asymptotic order in $t$, giving the following to first order in $\delta_{j}$ and asymptotic $t$
%
%
\begin{equation}
\displaystyle \xi_{{R}} \left( s \right) \approx {\frac {{t}^{2}\delta_{1}\,\Upsilon_{{R}} \left( it \right) }{\sigma\, \left( \sigma-1 \right) \pi}}
\mbox{}+{\frac {  t\delta_{2}\,\left( {\sigma}^{2}-\sigma-1 \right)\Upsilon_{{I}} \left( it \right) }{\sigma\, \left( \sigma-1 \right) \pi}}.
\label{Approx3b}
\end{equation}
Notice that the second term in \eqref{Approx3b} is smaller than the first term by one power of $t$, reflecting the fact that $\overset{o}{T}_2(s,v)$ vanishes at $v=t/2$. This means that asymptotically, we expect the second term in \eqref{Approx3b} to vanish, relative to the first term, except if $\Upsilon_R(it)\approx 0$. \newline

{\bf Remark:} By expanding about $\delta_j=0$, we are selecting for values of the integrand at $v={\it t}/2$. This could be interpreted as the introduction of a parameter $w_1$ to model the width of $\overset{o}{T}_1(s,v)$ near its peak - effectively replacing $\overset{o}{T}_1(s,v)$ by inserting a box of equal area around the peak and thereby approximating parts of the integrand near $v={\it t}/2$ which dominate $\overset{o}{T}_1(s,v)$ and thereby, the integral itself.  According to the Mean Value Theorem for Integrals, this can always be done and $w_1\neq0$. The parameter $w_1$ is generalized by $\delta_1$ which captures the area of the product term $\overset{o}{T}_1(s,v)\Upsilon_R(iv)$ and reasonable models for this parameter should be considered in future studies, recognizing that it is a (weak) function of both $\sigma$ and ${\it t}$ (but see Section \ref{sec:Caveats}). Similarly, since $\overset{o}{T}_2(s,v)$ vanishes anti-symmetrically about $v=t/2$, it is reasonable to expect that it will integrate somewhat close to zero. Thus we expect that terms involving $\delta_2$ will be much smaller than those involving $\delta_1$ which appears to be the case for large values of $t$.  Here, I use $\delta_1=\sqrt{\sigma(1-\sigma)}/2$ based on numerical experimentation with the approximation, and later, $\delta_4=\delta_1$. Since the relative width of the peak of $\overset{o}{T}_1(s,v)$ decreases for large values of $t$, this approximation should improve as $t$ increases - see Figure \ref{fig:FigMcOver}.\newline

With reference to the example presented in \eqref{Ex1a} and \eqref{Ex1b}, the two terms on the right-hand side of \eqref{Approx3b} with $\delta_{1}=\delta_{2}=0.2357$ ($=1/(3\sqrt{2})$ give
\begin{equation}
\xi_R(1/3+50i)\approx 3.63105\times10^{-15}-1.01201\times10^{-16}=3.52985\times10^{-15}
\end{equation}
compared to the correct result $\xi_R(1/3+50i)=3.3836\times10^{-15}$, a reasonable approximation. Furthermore, setting $\delta_1(1/3,50)=\delta_2(1/3,50)=0.2357$ in \eqref{Approx2} gives
\begin{align} \nonumber
&\frac{1}{\pi}\int_{24.76429}^{25.2357 }\overset{o}{T}_1(s,v)\Upsilon_R(s,v)v\,{\rm{d}}v+\frac{1}{\pi}\int_{24.76429}^{25.2357}\overset{o}{T_2}(s,v)\Upsilon_I(s,v)v\,{\rm d}v \\= 
&2.3675\times10^{-15}+5.0984\times10^{-16}=2.8773\times 10^{-15}
\label{NumApprox}
\end{align}
where I have used $\overset{o}{T}_{1,2}(s,v)=T_{1,2}(s,v)$ since the latter is dominated by the pole term over the narrow range of integration - see Figure \ref{fig:P1P2}. Note that the second term associated with $\delta_2$ is numerically smaller than the term associated with $\delta_1$ as predicted. Significantly, the order of magnitude of \eqref{NumApprox} is correct, demonstrating that the signal has been extracted, although the numerical estimate is not as good as the approximate one because the functions $\delta_{1,2}$ have not been generally selected for this purpose. In fact, resetting the limits of \eqref{NumApprox} using $\delta_1=\delta_2 = 0.2983$ yields the correct result to four significant figures. Of course there is no reason to expect that $\delta_2=\delta_1$; this becomes less of an issue as $t$ increases, since the relative contribution from the $\overset{o}{T}_{2}(s,v)$ term asymptotically vanishes except at a zero.\newline

{\bf Approximation 4}\newline

To obtain an asymptotic approximation of $\zeta_R(\sigma+i{\it t})$ as ${\it t}\rightarrow\infty$, first write the left-hand side of \eqref{Approx3b} in polar form and employ \eqref{LgammaSpec}, to obtain

\begin{equation}
\displaystyle \xi_{{R}} \left( s \right)\rightarrow -{{\rm e}^{-\pi\,t/4}}{\pi}^{-\sigma/2} \left| \zeta \left( \sigma+it \right)  \right| 
\mbox{} \sqrt{2\pi} \left( t/2 \right) ^{\sigma/2-1/2} \left(  \left(\sigma -1/2 \right) t\sin \left( \phi \right) +{t}^{2}\cos \left( \phi \right) /2 \right)\,. 
\label{Approx3c}
\end{equation}

Now, consider the right-hand side of \eqref{Approx3b} by identifying $\Upsilon_R(i{\it t})$ in terms of its components - see \eqref{Lamdef}, apply the same approximations employed above, reunite both sides of the equation, cancelling common non-zero factors on both sides and eventually arrive at
%
%
\begin{equation}
\displaystyle  \left| \zeta \left( s \right)  \right|  \left(  \left(\sigma -1/2 \right) \sin \left( \phi \right) +t\cos \left( \phi \right)/2  \right) 
\approx {\frac { t\left( \cos \left( \Phi \right)\delta_{1} + \left( {\sigma}^{2}-\sigma-1 \right) \sin \left(\Phi \right)\delta_{{2}}/t  \right)  }{\pi\,\sigma\, \left(1- \sigma \right) }  \left| \zeta \left( it \right)  \right|\left( {\frac {t}{2\pi}} \right)^{-\sigma/2}}
\label{Approx4}
\end{equation}


Starting from \eqref{EqmI}, a similar set of approximations can be applied. As in the model for $\xi_R(s)$, and with reference to Figure \ref{fig:P1234}, the background terms obtained from the decomposition of functions $T_3(s,v)$ and $T_4(s,v)$ integrate to zero exactly as in Section \ref{sec:Backg} with the change that each term contains an additional overall multiplicative factor of $(2\sigma-1)/t$. Again, the dominant contribution to $\xi_I(s)$ arises from the pole term associated with $T_4(s,v)$. Thus the equivalent of \eqref{Approx3b} is

\begin{equation}
\displaystyle \xi_{{I}} \left( s \right) \approx {\frac { 2\,\left( 2\,\sigma-1 \right) t^2 \left(\delta_{{4}}\Upsilon_{{I}} \left( it \right) +\delta_{{3}}/t\Upsilon_{{R}} \left( it \right)  \right) 
\mbox{}}{\pi\,\sigma\, \left( 1-\sigma \right) }}
\label{EqmR2}
\end{equation}

and the equivalent of \eqref{Approx3c} is
\begin{equation}
\displaystyle \xi_I \left( s \right) \rightarrow {{\rm e}^{-\pi\,t/4}}{\pi}^{-\sigma/2} \left| \zeta \left( \sigma+it \right)  \right| 
\mbox{} \sqrt{2\,\pi} \left( t/2 \right) ^{\sigma/2-1/2}\left(  \left(\sigma -1/2 \right) t\cos \left( \phi \right) -{t}^{2}\sin \left( \phi \right)/2  \right) 
\label{EqI4}
\end{equation}
Again, because $0<\sigma<1\ll{\it t}$, follow the exact steps used previously and ultimately arrive at the equivalent of \eqref{Approx4}


\begin{equation}
\displaystyle  \left| \zeta \left( \sigma+it \right)  \right|  \left(  \left(1/2 -\sigma \right) \cos \left( \phi \right) +t\sin \left( \phi \right)  \right) 
\approx{\frac {4 t\left(1/2 -\sigma \right) \left(\sin \left( \Phi \right) \delta_{{4}}+\cos \left( \Phi \right) \delta_{{3}}/t \right) 
}{\pi\,\sigma\, \left( 1-\sigma \right) }  \left| \zeta \left( it \right)  \right| \left( {\frac {t}{2\pi}} \right) ^{-\sigma/2}}
\label{EqI6a}
\end{equation}
\newline

{\bf Remark:} If $\sigma=1/2$ in \eqref{EqI6a}, $\sin(\phi)$ on the left-hand side vanishes. See \eqref{R_S} and Section \ref{sec:SigHalf}.\newline

Further, it is possible to extract approximations for $\zeta_R(\sigma+i{\it t})$ and $\zeta_I(\sigma+i{\it t})$ by solving \eqref{Approx4} and \eqref{EqI6a} simultaneously, using a simple trigonometric expansion (because $\arg(\zeta(\sigma+i{\it t}))$ is embedded within the definition of $\phi$ (see \eqref{phi}), by noting that $\cos(\arg(\zeta(s))) \left| \zeta \left(s \right)  \right|=\zeta_R(s)$ (and similarly for $\zeta_I(s)$). Asymptotically,(because both $\cos(\Phi)$ and $\sin(\Phi)$ cannot both vanish simultaneously) this solution yields a simple representation after applying\newline

{\bf Approximation 5:}
\begin{equation}
\delta_{2}/t=\delta_{3}/t=0, 
\end{equation}
giving
%
%
\begin{equation}
\displaystyle \zeta_{{R}} \left( \sigma+it \right) \approx 2\,\left| \zeta \left( it \right)  \right| {\frac { 
\mbox{} \left( -\cos \left( \omega \right) \cos \left( \Phi \right) \delta_{{1}}-2\, \left( \sigma-1/2 \right) \sin \left( \Phi \right) \sin \left( \omega \right) \delta_{{4}}
\mbox{} \right) }{\sigma\, \left( \sigma-1 \right) \pi} \left( {\frac {t}{2\pi}} \right) ^{-\sigma/2}}
\label{ZReal}
\end{equation}
and
\begin{align}
\displaystyle \zeta_{{I}} \left( \sigma+it \right)\approx 2\left| \zeta \left( it \right)  \right| {\frac { 
\mbox{} \left( -\sin \left( \omega \right) \cos \left( \Phi \right) \delta_{{1}}+2\, \left( \sigma-1/2 \right) \cos \left( \omega \right) \sin \left( \Phi \right) \delta_{{4}}
\mbox{} \right) }{\sigma\, \left( \sigma-1 \right) \pi} \left( {\frac {t}{2\,\pi}} \right) ^{-\sigma/2}}
\label{ZImag}
\end{align}

where
\begin{equation}
\omega=\arg(\Gamma(\sigma/2+i{\it t}/2))-\frac{{\it t}}{2}\log(\pi).
\label{omega}
\end{equation}
Notice that the right-hand sides of both \eqref{ZReal} and \eqref{ZImag} do not contain any components of $\zeta(\sigma+i{\it t})$ and therefore the left-hand sides are expressed in terms of independent functions. Compare with \cite{2012arXiv1201.2633F}.

\begin{figure}[h] 
\centering 
\subfloat  [This is a comparison of the left- and right-hand sides of \eqref{ZReal} for $\sigma=0.7$ near ${\it t}=500$.]
{
\includegraphics[width=.44\textwidth]{{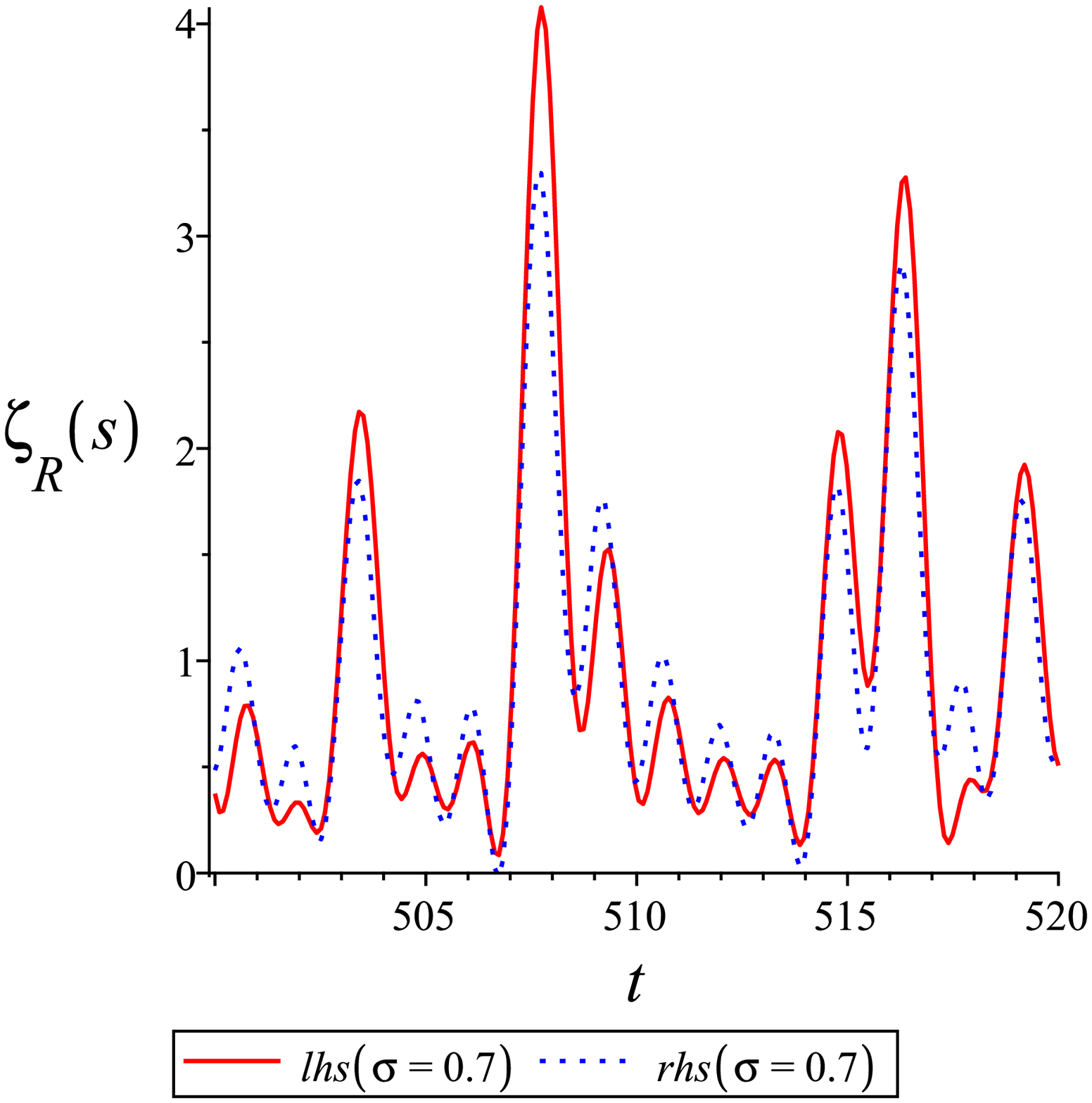}} \label{fig:ZRfig} 
}
\hfill                                 
\subfloat [This is a comparison of the left- and right-hand sides of \eqref{ZImag} for $\sigma=0.7$ near ${\it t}=500$.]
{
\includegraphics[width=.44\textwidth]{{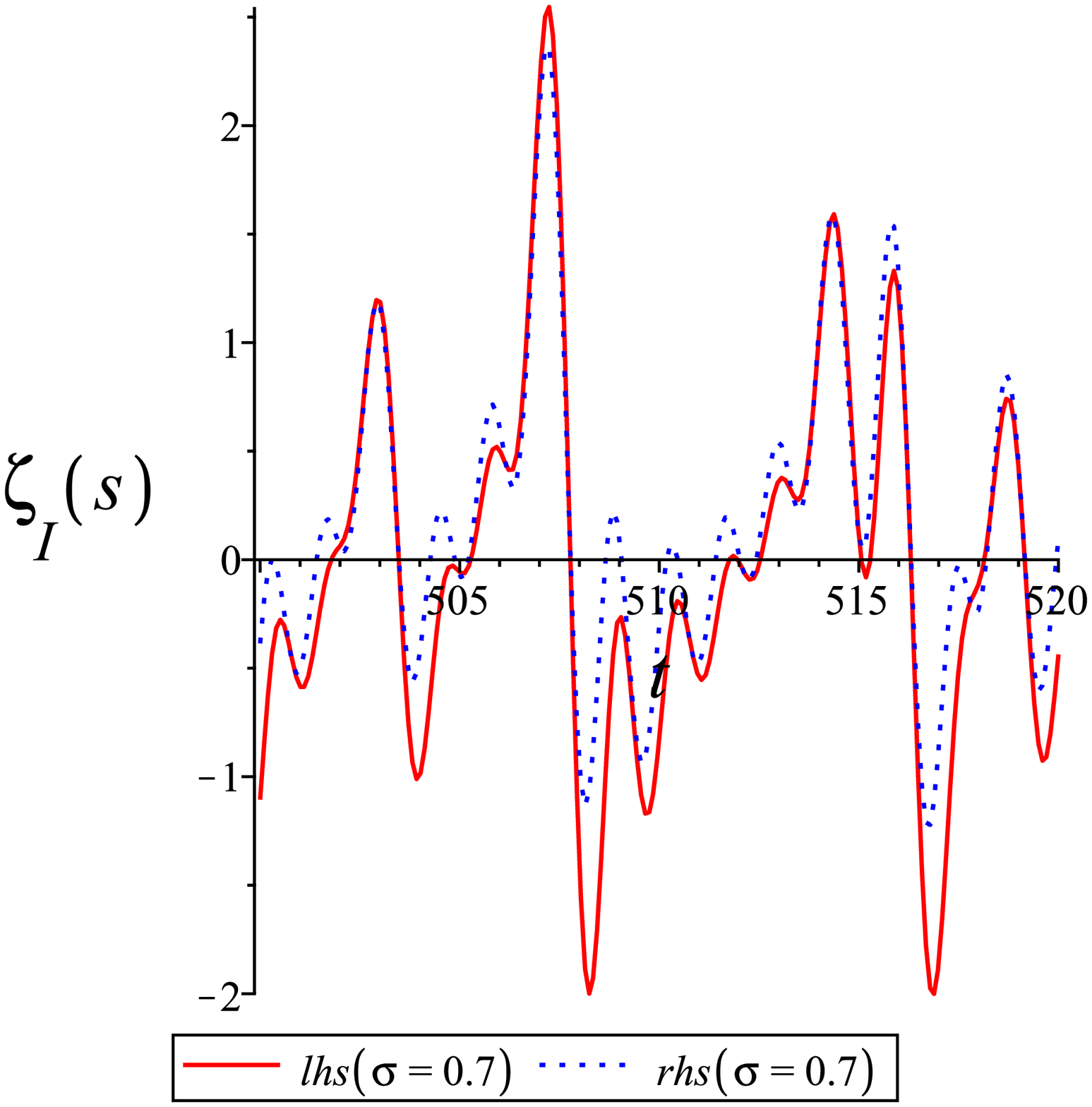}} \label{fig:ZIfig} 
}
\caption{This figure illustrates the approximations \eqref{ZReal} and \eqref{ZImag} to $\zeta_R(\sigma+i{\it t})$ and $\zeta_I(\sigma+i{\it t})$ respectively near ${\it t}=500$ using $\delta_1=\delta_4=\sqrt{\sigma(1-\sigma)}/2$.} 
\label{fig:ZRZI}
\end{figure}

\subsubsection{Numerical Examples} 

\begin{figure}[h] 
\centering
\begin{subfloat} [This is a comparison of $\xi_R(\sigma+i{\it t})$ for two very different values of $\sigma$ near ${\it t}=700$.]
{
\includegraphics[width=.45\textwidth]{{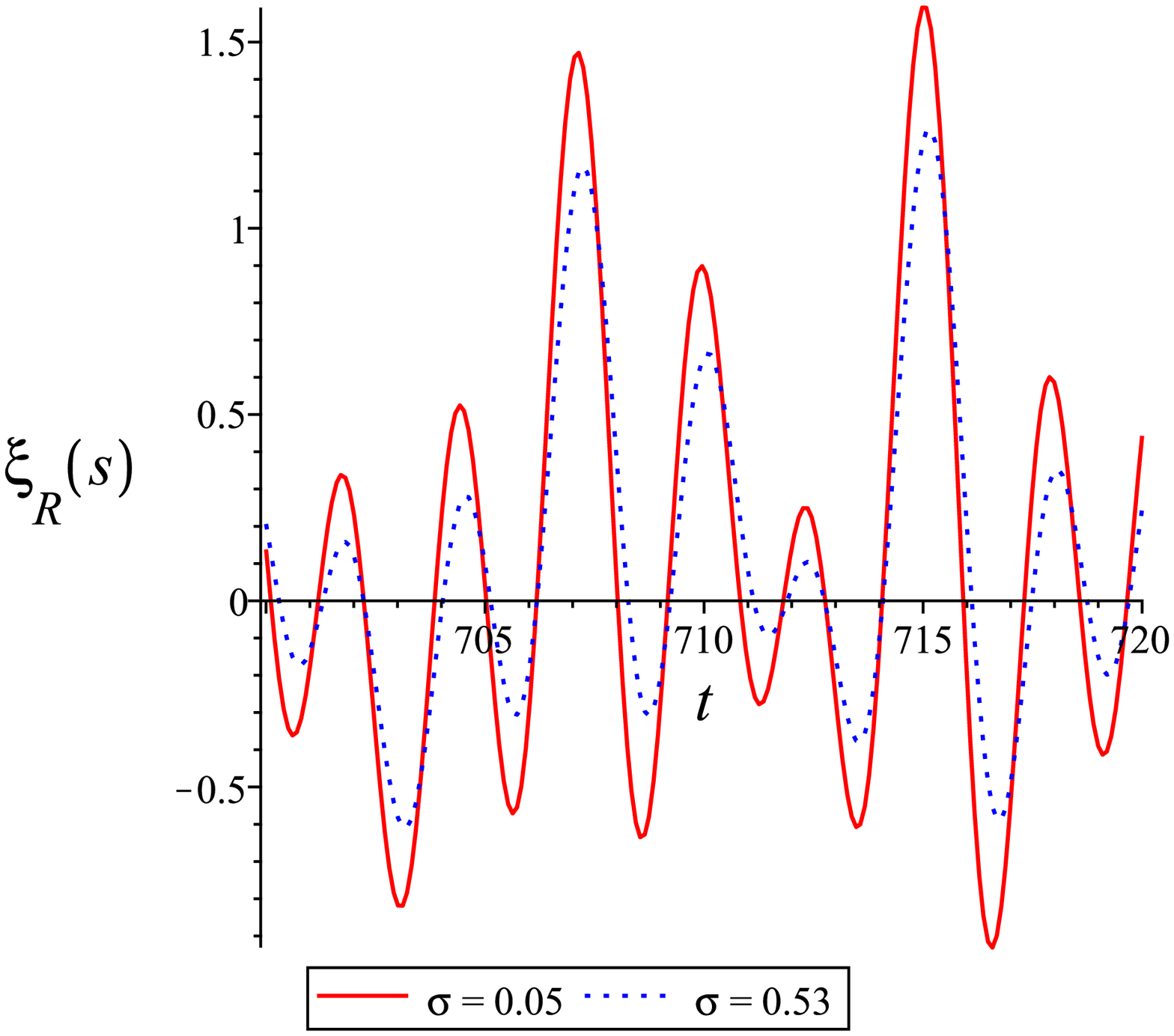}}
\label{fig:TwoXsi}
}                                 
\end{subfloat}
\hfill 
\begin{subfloat}[This is a comparison of the left- and right-hand sides of \eqref{XR-XL} for two very different values of $\sigma$ near ${\it t}=500$.]
{
\includegraphics[width=.45\textwidth]{{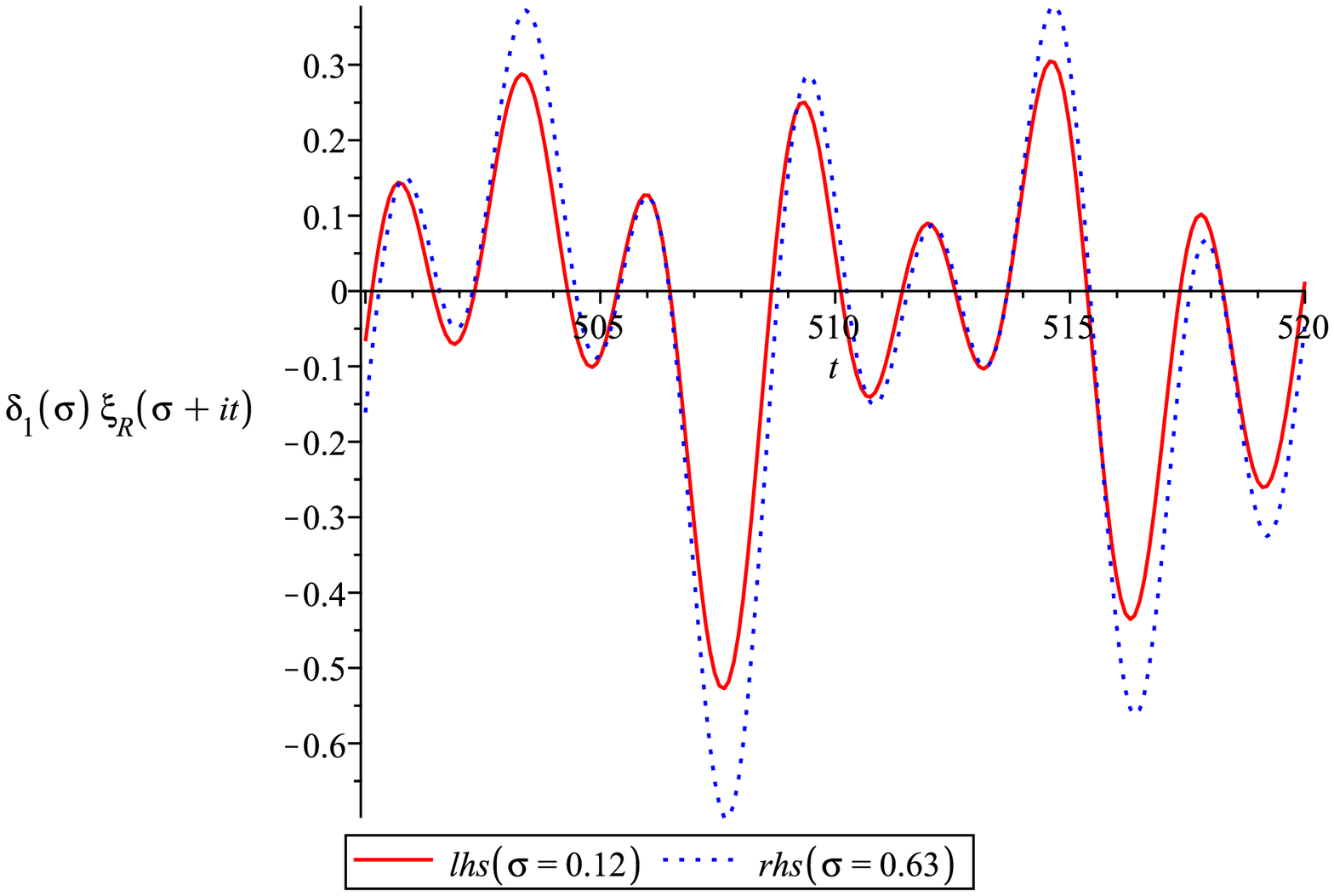}}
\label{fig:Q1Q2}
}
\end{subfloat}
\caption{For several values of $\sigma$ and $t$, this is a comparison of $\xi_R(\sigma_1+i{\it t})$ and $\xi_R(\sigma_2+i{\it t})$ as they exist (left), and as they are predicted to exist by the approximation (right) - see \eqref{XR-XL}. In both cases, the functions have been scaled by a factor $\exp(\pi{\it t}/4)/{\it t}^2$.}
\end{figure}

Figure \ref{fig:ZRZI} illustrates both approximations \eqref{ZReal} and \eqref{ZImag} for an arbitrary choice of $\sigma$. From \eqref{Approx3b}, it is seen that any structure in $\xi_R(\sigma+i{\it t})$ must originate from the function $\Upsilon(i{\it t})$ on the right-hand side, and this function lacks $\sigma$ dependence. Thus \eqref{Approx3b} predicts that the structure of $\xi_R(\sigma+i{\it t})$, particularly the location of the zeros, must be relatively independent of $\sigma$, since, as was noted in Section \ref{sec:ZsApprox}, $\sigma$ dependence only resides in the transfer functions, and to a lesser extent any model developed to approximate $\delta_{1,2}$ as a function of $\sigma$ and/or $t$. Thus, in \eqref{EqpR2} if $\Upsilon_R(iv)$ happens to vanish at $v={\it t}/2$, then the peak associated with $T_1({\it t},v)$ is multiplied by a vanishing factor and contributions to the integral from both functions $T_1({\it t},v={\it t}/2)$ and $T_2({\it t},v={\it t}/2)$ vanish - independently of $\sigma$!. So, the location of the zeros is approximately governed by $\cos(\Phi)=0$. See the discussion surrounding \eqref{CosPhi0} and Figure \ref{fig:CosPhi}. This prediction can be checked - see Figure (\ref{fig:TwoXsi}) which  demonstrates that this property appears to be approximately true, for two randomly chosen, disparate values of $\sigma$.\newline

From this observation, arises a constraint on the approximation. Consider two instances of \eqref{Approx3b} with different choices $\sigma=\sigma_1$ and $\sigma=\sigma_2$. The approximation predicts that 
\begin{equation}
\displaystyle {\frac {\xi_{{R}} \left( \sigma_{{1}}+i{\it t} \right) }{\xi_{{R}} \left( \sigma_{{2}}+i{\it t} \right) }}={\frac {\delta_{{1}} \left( \sigma_{{1}} \right) 
\mbox{}\sigma_{{2}} \left( 1-\sigma_{{2}} \right) }{ \delta_{{1}} \left( \sigma_{{2}} \right)\,\sigma_{{1}} \left( 1-\sigma_{{1}} \right) }}
\label{XiRat}
\end{equation}
since the model suggests that $\delta_1$ is at least a (weak) function of $\sigma$. Now suppose that $\sigma_2=1-\sigma_1$, in which case the left hand side of \eqref{XiRat} becomes unity (see \eqref{xidef}); for consistency, this imposes a constraint on the parameter $\delta_1$, specifically
\begin{equation}
\delta_1(\sigma)=\delta_1(1-\sigma)\,.
\label{w1Strain}
\end{equation}

With this constraint in mind, and consistent with the approximations specified previously, it is possible to obtain some insight into a reasonable model for $\delta_1(\sigma,t)$ by solving \eqref{Approx3b} using the approximation $\delta_{2}/t=0$ - see Approximation 5. Introduce $\sigma,t$ dependence by letting the solution be denoted by $\Delta_1(\sigma,t)$, and consider $\Delta_{1}(\sigma,t)$ over the range $0\leq \sigma\leq 1$, for several choices of $t$. The result is given in Figure \ref{fig:Delta1S}, providing evidence that it is likely that a smooth function exists to capture the $\sigma$ dependence of $\delta_1(\sigma,t)$, although the normalization factor varies considerably as a function of $t$. The somewhat arbitrary and simple model 

\begin{equation}
\delta_1(\sigma,t)=\sqrt{ \sigma(1-\sigma) }/2\,.
\label{w1def}
\end{equation}

used throughout this paper is also shown; it is clear that this model, although a reasonable first approximation, requires future refinement. A similar result can be found for $\Delta_4(\sigma,t)$, the solution of \eqref{EqmR2} - see Figure \ref{fig:Delta4S}. (Remark: To obtain the normalization of the latter Figure, expand the function $\xi(\sigma+it)$ about the point $\sigma=1/2$ and cancel the numerator and denominator factors $(\sigma-1/2)$ that appear in the imaginary part). \newline

\begin{figure}[h] 
\centering
\begin{subfloat} [This Figure shows the solution function $\Delta_1(\sigma,t)$ for several values of $t$ over a range of $\sigma$. Each of the curves has been normalized by $\Delta_1(1/2,t)$ with values shown.]
{
\includegraphics[width=.45\textwidth]{{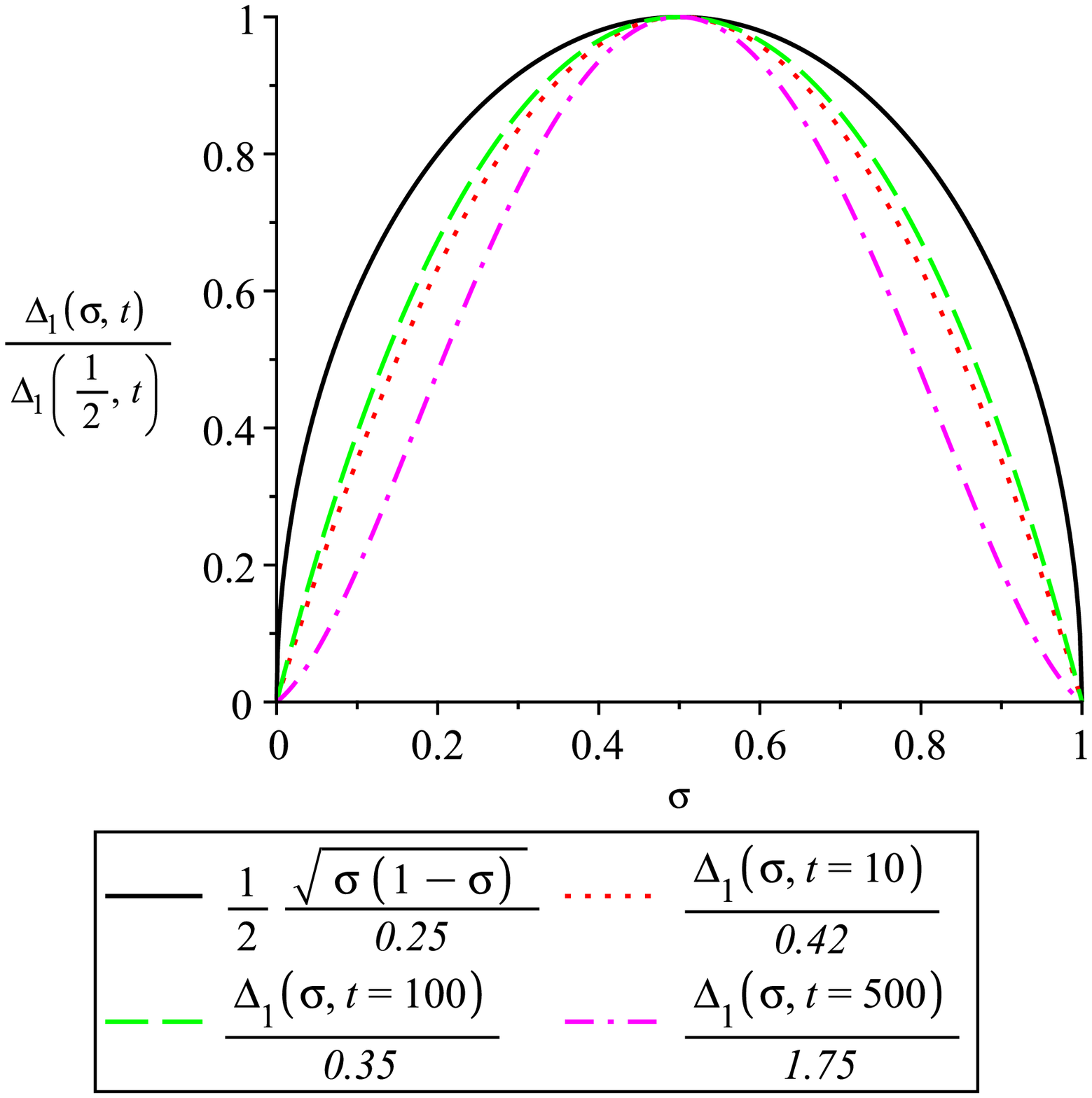}}
\label{fig:Delta1S}
}                                
\end{subfloat}
\hfill 
\begin{subfloat}[This Figure shows the solution function $\Delta_4(\sigma,t)$ for several values of $t$ over a range of $\sigma$. Each of the curves has been normalized by $\Delta_4(1/2,t)$ with values shown.]
{
\includegraphics[width=.45\textwidth]{{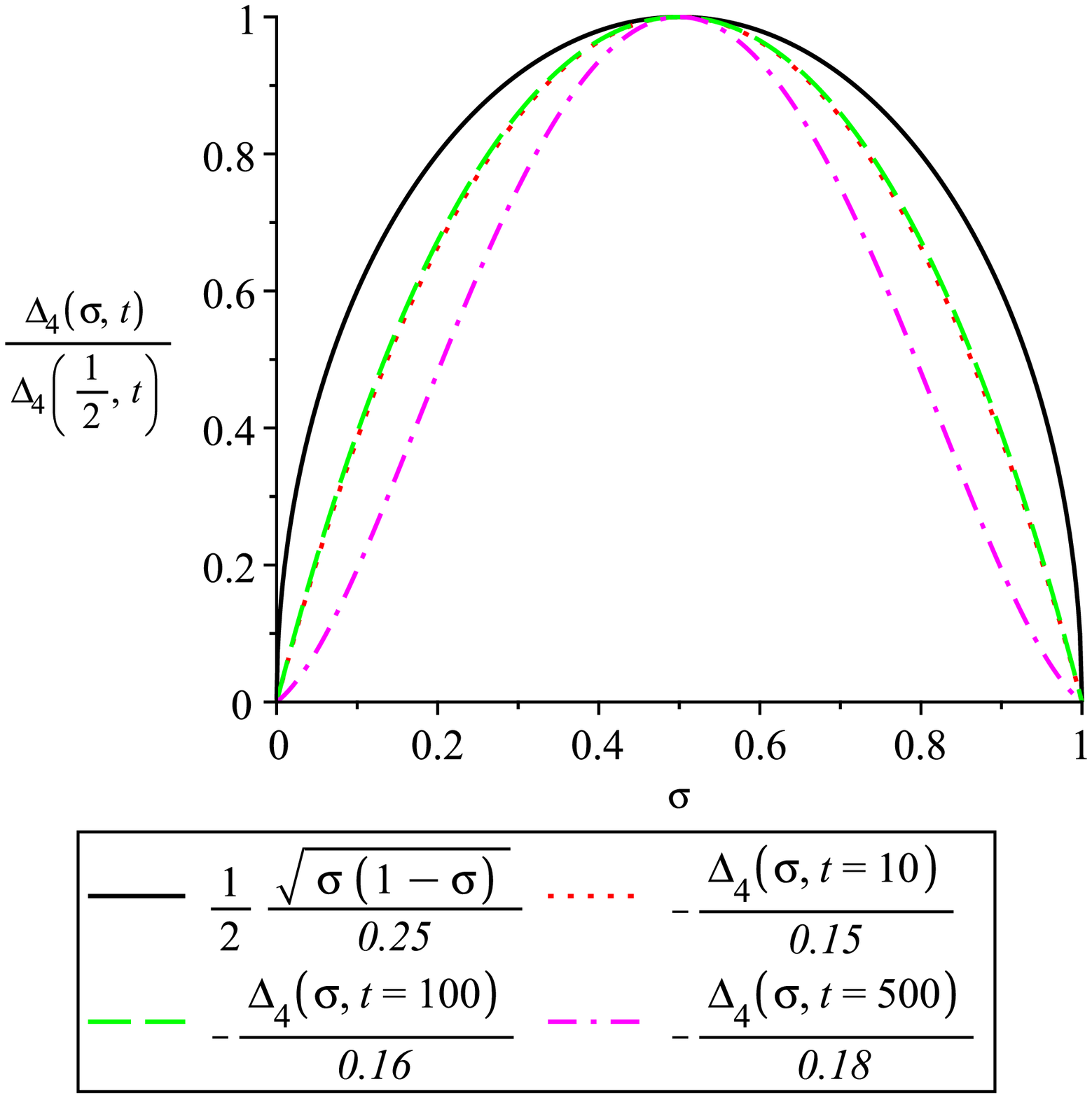}}
\label{fig:Delta4S}
}
\end{subfloat}
\caption{This Figure indicates that a simple model for both $\delta_1(\sigma,t)$ and $\delta_4(\sigma,t)$ exists as a function of $\sigma$, although the dependence on $t$ may not be simple - see Figure (\ref{delta4}) below.}
\label{fig:Delta14S}
\end{figure}

All this suggests that it is of interest to investigate $\xi(\sigma+i{\it t})$ for two different choices of $\sigma$. \eqref{XiRat} and \eqref{w1def} predict that 
\begin{equation}
\displaystyle  \sqrt{\sigma_{{1}} \left( 1-\sigma_{{1}} \right) }\xi_{{R}} \left( \sigma_{{1}}+i{\it t} \right) 
\mbox{}\approx\xi_{{R}} \left( \sigma_{{2}}+i{\it t} \right)  \sqrt{\sigma_{{2}} \left( 1-\sigma_{{2}} \right) }\,,
\label{XR-XL}
\end{equation}
an example of which is shown in Figure (\ref{fig:Q1Q2}). Finally, it is of interest to consider approximations to $|\zeta(\sigma+i{\it t})|$ as given in \eqref{Approx4}. See Figure \ref{fig:TwoSigmas}.

\begin{figure}[ht] 
\centering
\includegraphics[width=.8\textwidth]{{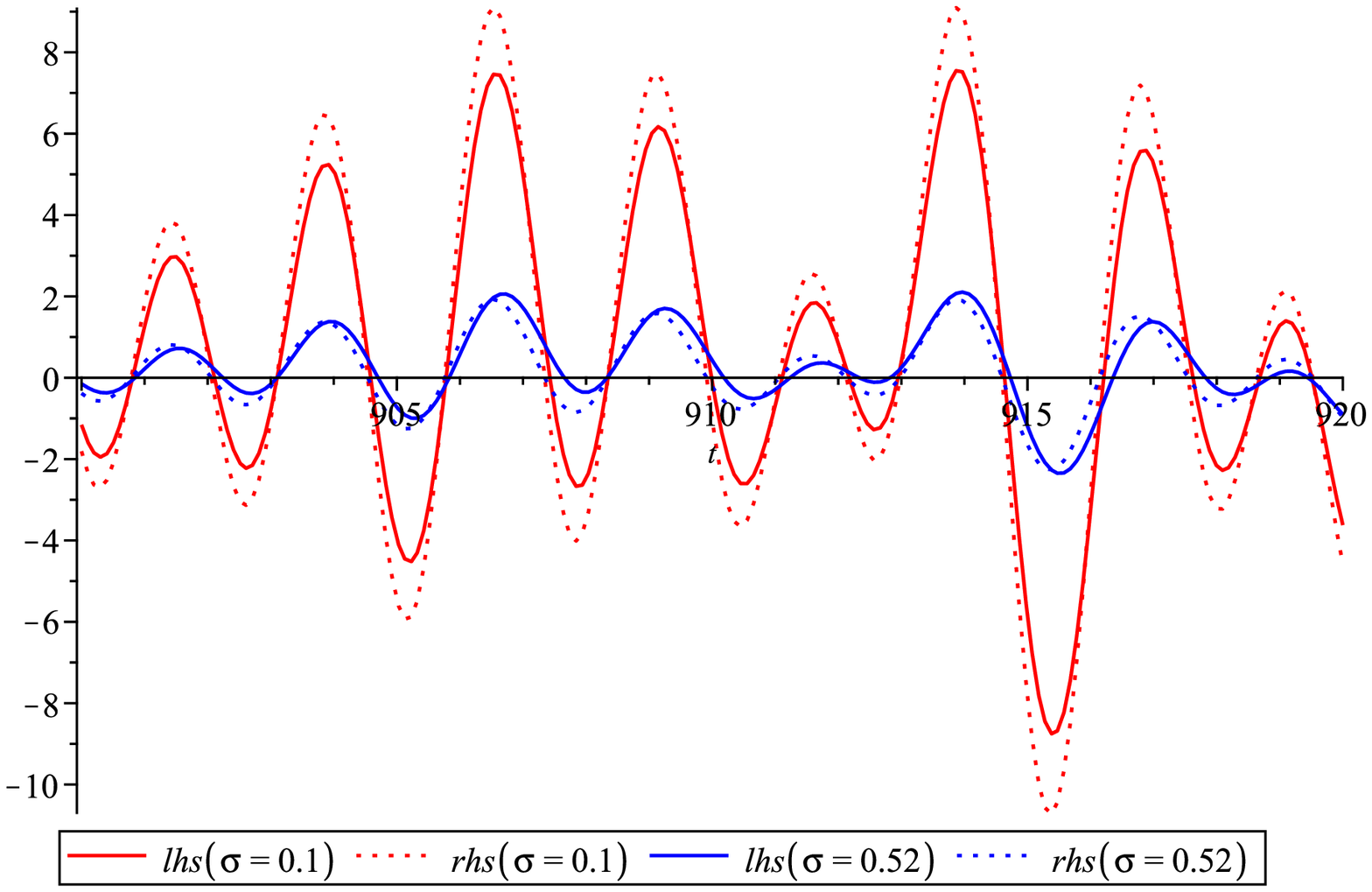}}
\caption{This is a comparison of the left- and right-hand sides of \eqref{Approx4} for two very different values of $\sigma$ near ${\it t}=900$. The parameter $\delta_1=\sqrt{\sigma(1-\sigma)}/2$ in both cases and all curves have been normalized by the factor $1/t$.}
\label{fig:TwoSigmas}
\end{figure}

\subsection{The case $\sigma=1/2$} \label{sec:SigHalf}

From the Real part of $\xi(s)$, \eqref{Approx4} is of particular interest when $\sigma=1/2$. By straight substitution, using \eqref{w1def} and \eqref{FeqId}, we obtain

\begin{equation}
\displaystyle  \left| \zeta \left( 1/2+it \right)  \right| \cos \left( \phi \right)\approx 4\delta_{{1}} \left( {\scriptstyle \frac{1}{2}},t \right)\,{\frac { 
\cos \left( \Phi \right)  \left| \zeta \left( 1-it \right)  \right| {t^{1/4}}{2}^{3/4}}{{\pi}^{5/4}}}\,.
\label{Req1}
\end{equation}

\begin{figure}[h] 
\centering
\includegraphics[width=.5\textwidth]{{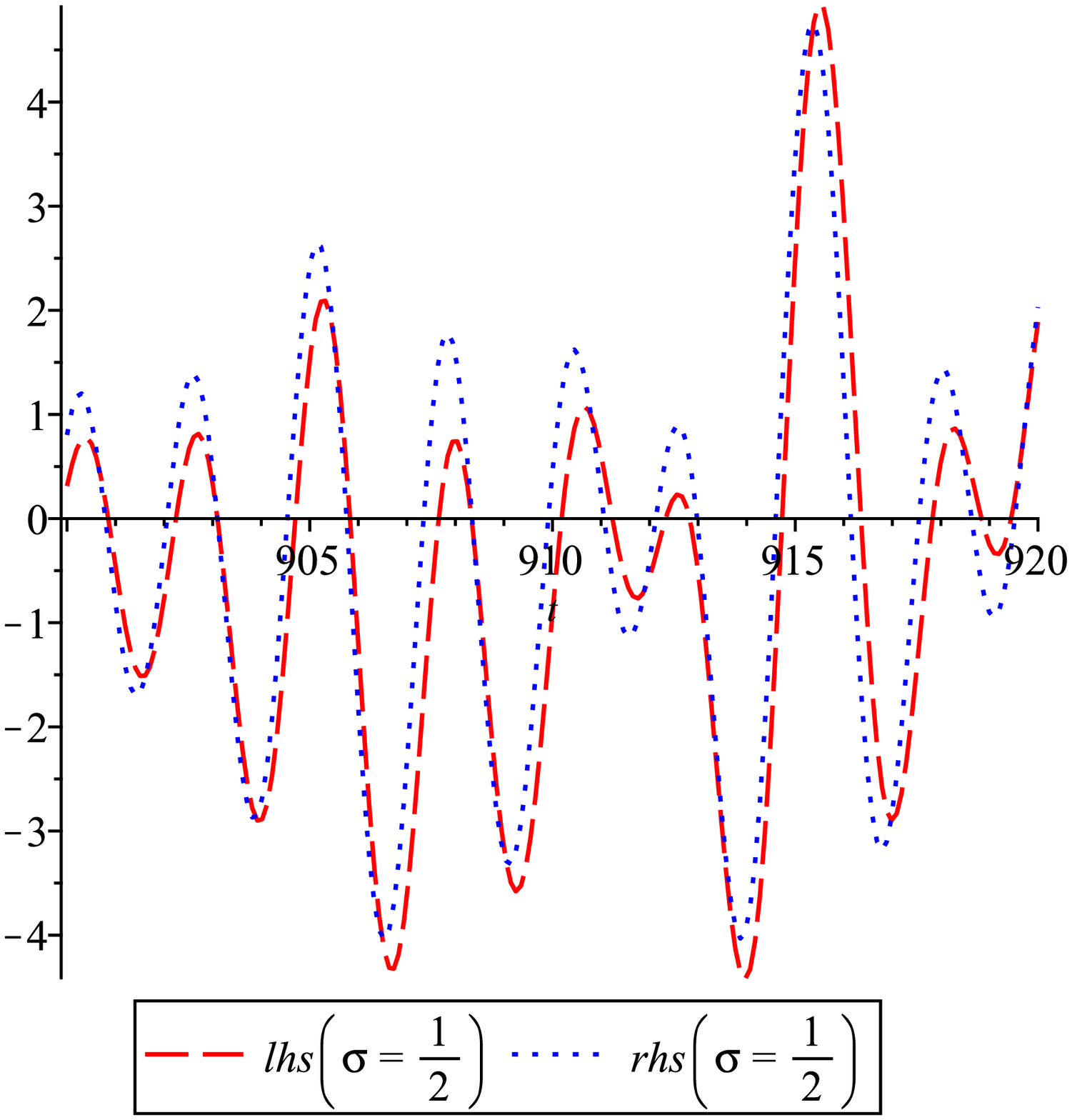}}
\caption{This is a comparison of the left- and right-hand sides of \eqref{Req1} near ${\it t}=900$. The parameter $\delta_1=\sqrt{\sigma(1-\sigma)}/2$ and $\sigma=1/2$.}
\label{fig:ReqHalf}
\end{figure}

This is consistent with \cite[Eq. (5.1.8)]{Titch2}; for a numerical example, see Figure \ref{fig:ReqHalf}. In \eqref{Req1}, note that the factor $\cos(\phi)=(-1)^k$ because of \eqref{phi}, where $k$ counts the number of zeros of $\zeta(1/2+i{\it t})$ relative to ${\it t}=0$. For further discussion, see Section \ref{sec:Lind}.\newline 

From the Imaginary part of $\xi(s)$, in the case of \eqref{EqmR2} with $\sigma=1/2$, both  sides vanish, as they must because it is known that $\xi_I(1/2+i{\it t})=0$. This is also true of \eqref{EqI6a} because in the limit $\sigma=1/2$, we have $\phi=k\pi$ due to \eqref{R_S} and \eqref{phi}. From \eqref{xiExpand} by straightforward differentiation of \eqref{xidef} applied to \eqref{EqmR2}, we obtain  

\begin{align} \nonumber
\lim_{\substack{{\it \sigma}\rightarrow 1/2}} &\frac{\xi_I(1/2+it)}{(\sigma-1/2)}=\Im \left( {\frac {\rm \partial}{{\rm \partial}{\sigma}}}\xi \left( {\sigma+it} \right)|_{\sigma=1/2} \right)  \\ \nonumber
 \displaystyle&=\frac{\left( -1 \right) ^{k} \left| \Gamma \left( \frac{1}{4}+it/2 \right)  \right|}{{\pi}^{1/4}}   \left(  \left( t-\frac{1}{4} ( {t}^{2}+1/4 ) \Im ( \psi ( 1/4+it/2)  )  \right)  \left| \zeta \left( 1/2+it \right)  \right| \right. \\ 
 \nonumber
&\left.+\frac{1}{2} \left( {t}^{2}+1/4 \right) \sin \left( \alpha -\beta  \right)  \left| \zeta^{\prime} \left( 1/2+it \right)  \right| 
\mbox{} \right)\\
&\approx  {\frac {4\, t^2 \left(\delta_{{4}}\Upsilon_{{I}} \left( it \right) +\delta_{{3}}/t\Upsilon_{{R}} \left( it \right)  \right) 
\mbox{}}{\pi\,\sigma\, \left( 1-\sigma \right) }} 
\label{XsiD}
\end{align}
where $\psi(s)$ is the digamma function, $\alpha$ and $\beta$ are defined in \eqref{alpha-def} and \eqref{beta-def} respectively, and asymptotically this becomes


\begin{equation} 
\displaystyle ({2\pi})^{3/4} \left( \frac{1}{2} \psi_{I} \left( 1/4+it/2 \right)    \left| \zeta \left( 1/2+it \right)  \right| 
-\sin ( \alpha -\beta )  \left| \zeta^{\prime} \left( 1/2+it \right)  \right|  \right) 
 \approx -{ {64\,\delta_{4} \left| \zeta \left( it \right)  \right| \sin \left(\Phi  \right)t^{-1/4} 
}{{}}} \,.
\label{Ceq3}
\end{equation}

For this fixed value of $\sigma=1/2$, it is possible to obtain an estimate of $\delta_4(1/2,t)$ by solving \eqref{Ceq3} for $\delta_4$. Denote the solution by $\Delta_4(1/2,t)$ and consider Figure (\ref{delta4}) which shows $\Delta_4(1/2,t)$ for a range of $t$. As a first approximation, choose $\delta_4(1/2,t)=0.16$ as suggested in the Figure. It should be noted that the denominator of the solution can vanish but the numerator does not vanish simultaneously -  these singular points of the solution are apparent in the Figure; however at such points, the parameter $\delta_3(1/2,t)$ must be incorporated into the model; this will prevent the singularity from occurring. Thus the approximate solution in Figure (\ref{delta4}) is valid in between the singular points; the horizontal line $\Delta_4(1/2,t)\approx 0.16$ used in this example has been chosen on that basis.


\begin{figure}[h] 
\centering
\begin{subfloat} [The curve in this Figure shows the solution $\Delta_4(1/2,t)$ of \eqref{Ceq3} over a range of $t$. The horizontal line $\Delta_4(1/2,t)\approx 0.16$ illustrates a first approximation.]
{
\includegraphics[width=.40\textwidth] {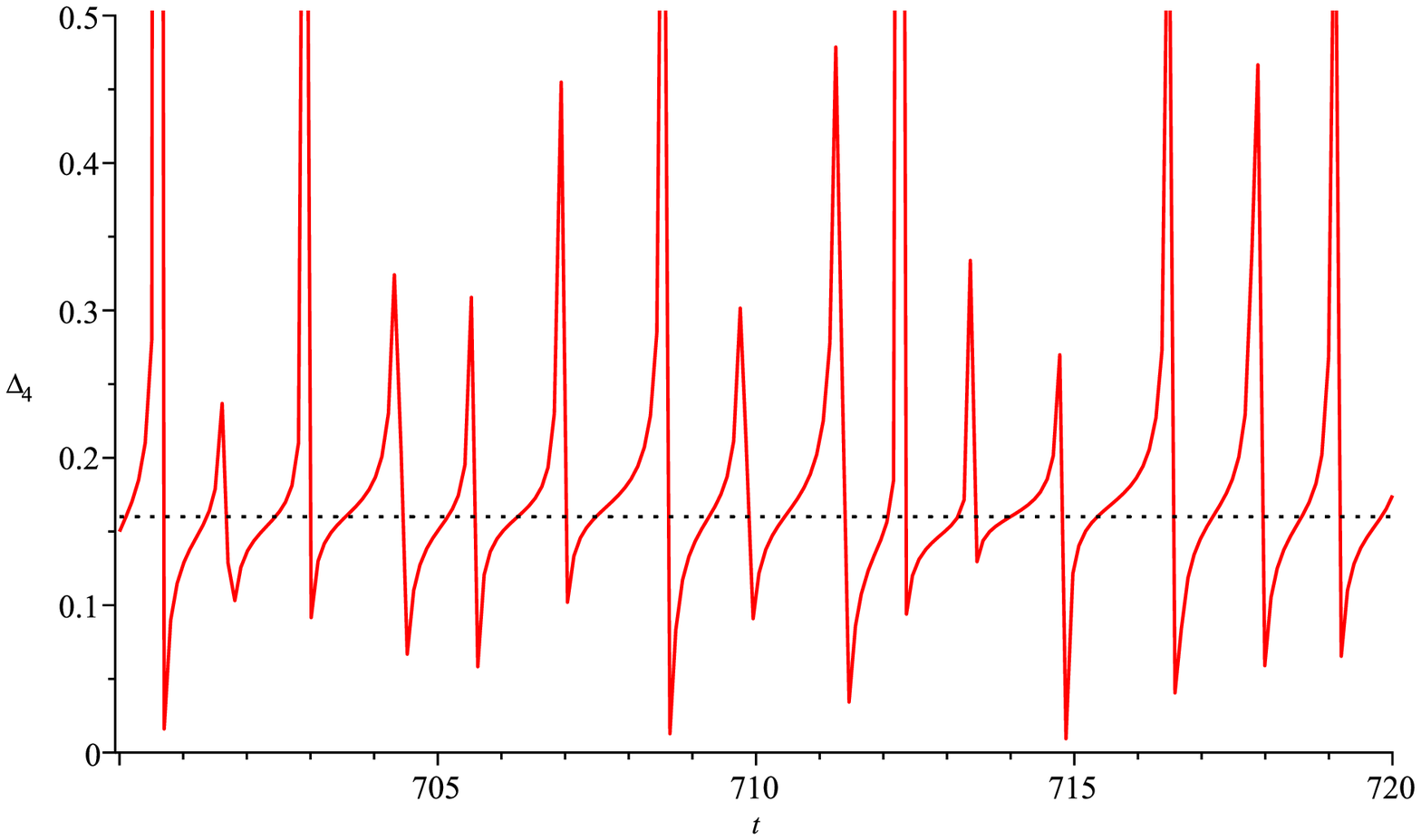}
\label{delta4}
}
\end{subfloat}
\quad \quad
\begin{subfloat} [This is a comparison of the absolute value of the right and left-hand sides of \eqref{XsiD} using $\delta_4(1/2,t)=0.16$.]
{
\includegraphics[width=.45\textwidth]{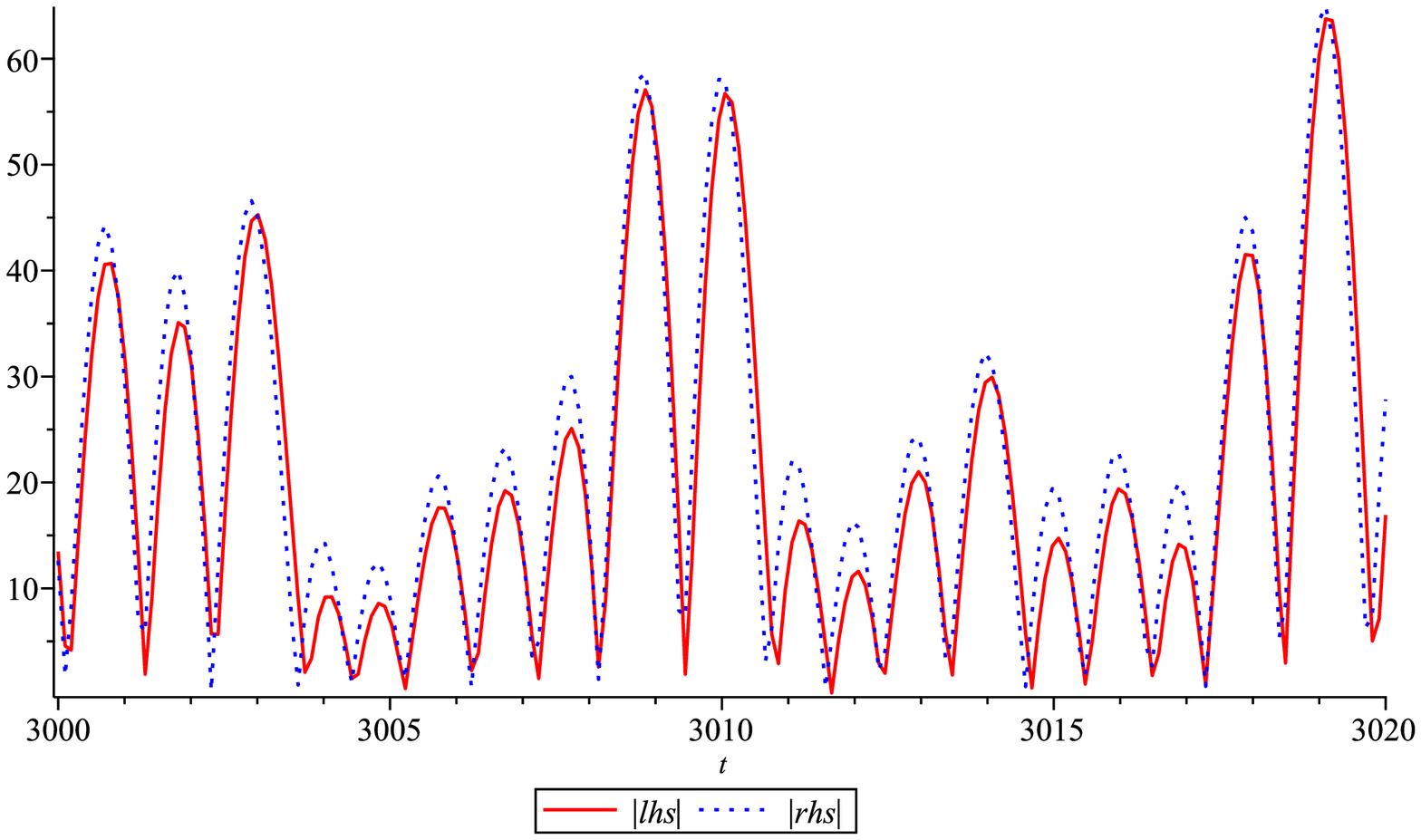}

\label{XiImag}
}
\end{subfloat}
\caption{A test (a) of the solution to \eqref{Ceq3} over the range $700\leq t \leq 720 $ applied (b) to that same equation over a completely different range of $t$.}
\label{fig:delta4}
\end{figure}

\section{Zeros of $\zeta(\sigma+i{\it t})$} \label{sec:zeros}

As noted earlier, $c=-1$ was chosen in \eqref{Beq0} because it would lead to an expression involving a master function $\zeta(1- i{\it t})$ (see \eqref{FeqId}) whose properties are fairly well-known, and whose validity spans the critical strip. In particular, it is known that $|\zeta(1-i{\it t})|\neq0$, (see Ivic, \cite[Theorem (6.1)]{Ivic}), although both $\zeta_R(1-i{\it t})$ and $\zeta_I(1-i{\it t})$ do have many (non-coincident) zeros, and their location is the subject of some interest (e.g. \cite{2011arXiv1112.4910A}). Further, it is well-known that $\xi(s)$ shares the zeros of $\zeta(s)$ and $\Upsilon(1-i{\it t})$ shares the zeros of $\zeta(1-i{\it t})$ so any properties of the zeros belonging to the first of each pair, apply to the second (of each pair).\newline

Therefore consider \eqref{Approx3b}, \eqref{EqmR2} in the company of Approximation 5. The first of these relates $\xi_R(s)$ to $\Upsilon_R(i{\it t})$ and predicts that the zeros of the these functions should coincide, and, in fact the two should be proportional to one another if $\sigma\neq0$ and $\sigma\neq1$ where the model ($c=-1$) breaks down anyway. The second of these relates  $\xi_I(s)$ to $\Upsilon_I(i{\it t})$ and predicts a further coincidence between the zeros of $\xi_I(\sigma+i{\it t}))$ and $\Upsilon_I(i{\it t})$ unless $\sigma=1/2$, in which case $\xi_I(1/2+i{\it t})$ is always zero \cite[Section 6.5]{Edwards}; $\Upsilon_I(i{\it t})$ is therefore independent if $\sigma=1/2$. But, as noted, because $\,|\zeta(1-i{\it t})|\neq0|$ then no zeros of $\Upsilon_R(i{\it t})$ and $\Upsilon_I(i{\it t})$ coincide (see also \eqref{FeqId}). Consequently, neither do the zeros of $\xi_R(\sigma+i{\it t})$ and $\xi_I(\sigma+i{\it t})$. Further, if $\xi(\sigma+i{\it t})=0$ and $\sigma\neq 1/2$, this would require the right-hand sides of \eqref{Approx3b} and \eqref{EqmR2} to vanish simultaneously, a contradiction with the known properties of $|\zeta(1-i{\it t})|$ because neither of the two parameters $\delta_1$ and $\delta_4$ vanish. \newline

Therefore, if $\sigma\neq1/2$, it is impossible for $\xi(\sigma+i{\it t})$ to vanish because it is overly constrained by the properties of $\zeta(1-i{\it t})$; if $\sigma=1/2$, one of the constraints vanishes, and $\xi(1/2+i{\it t})=0$ becomes possible. Thus, within the confines of this approximation, the Riemann Hypothesis (RH) is true. However, none of the above predicts that $\zeta(1/2+i{\it t})=0$ anywhere. That is a separate issue. A comparison between the right- and left-hand sides of \eqref{Approx3b} and \eqref{EqmR2} is given in Figure \ref{fig:XiVsLam}, demonstrating that the discussion here closely approximates reality.\newline

\begin{figure}[h] 
\centering
\begin{subfloat} [This is a comparison of the left- and right-hand sides of \eqref{Approx3b} near ${\it t}=800$. The parameter $\delta_1=\sqrt{\sigma(1-\sigma)}/2$, $\sigma=1/2$ and $\delta_2=0$.]
{
\includegraphics[width=.45\textwidth] {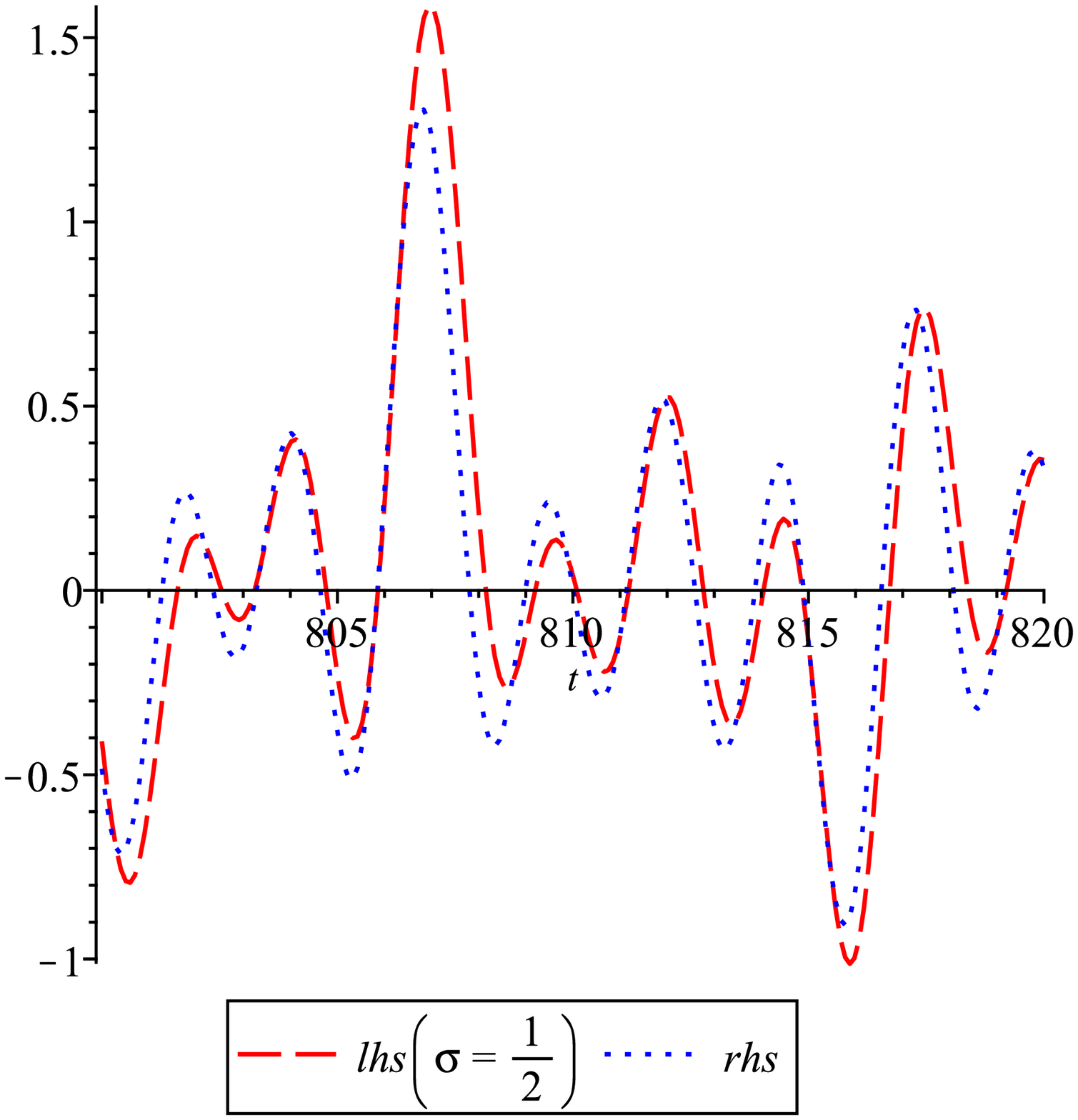}
\label{XiVsLamR}
}
\end{subfloat}
\quad
\begin{subfloat} [This is a comparison of the left- and right-hand sides of \eqref{EqmR2} near ${\it t}=800$. The parameter $\delta_4=\sqrt{\sigma(1-\sigma)}/2$, $\sigma=0.22$ and $\delta_3=0$.]
{
\includegraphics[width=.45\textwidth]{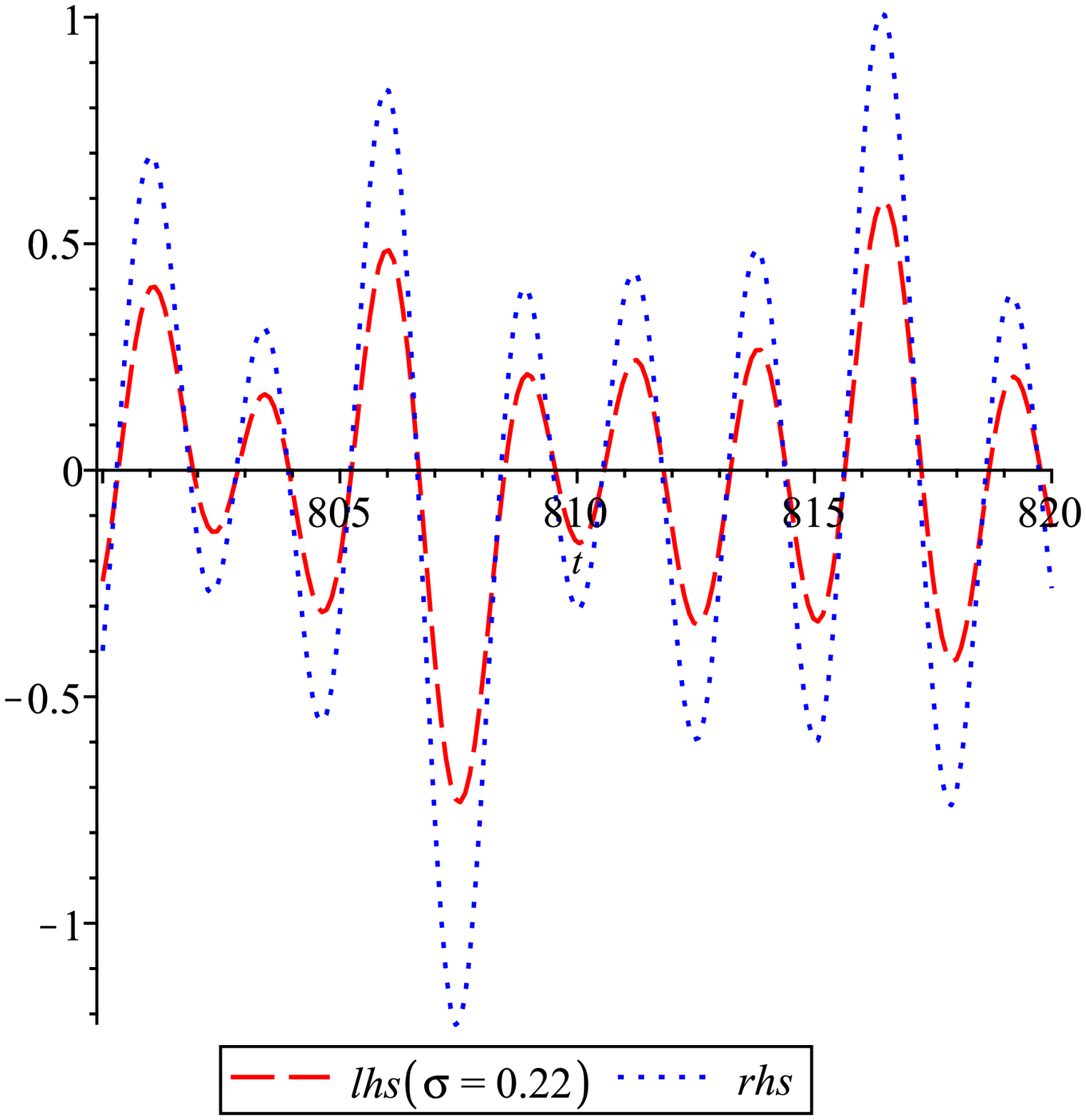}
\label{XiVsLamI}
}
\end{subfloat}
\caption{This Figure gives a comparison of the predicted relationship between $\xi_R(\sigma+i{\it t})$ (left) and $\xi_I(\sigma+i{\it t})$ (right) versus the corresponding real and imaginary components of $\Upsilon(i{\it t})$ at two different choices of $\sigma$ near ${\it t}=800$.  The figures are scaled by $e^{\pi{\it t}/4}/{\it t}^2$.}
\label{fig:XiVsLam}
\end{figure}

It is also worthwhile to consider how \eqref{ZReal} and \eqref{ZImag} allow us to reach the same conclusion. Suppose that $\zeta_R(\sigma+i{\it t})=0$ on the left-hand side of \eqref{ZReal} and simultaneously $\zeta_I(\sigma+i{\it t})=0$ on the left-hand side of \eqref{ZImag}. Since it is known that $|\zeta(1-i{\it t})|\neq 0$ (and therefore $|\zeta(it)|\neq 0)$ it must be that each of the terms in parenthesis in those two equations vanish simultaneously. That is, after removing all non-zero common factors, and noting that the equations are only valid for $0<\sigma<1$ we have
\begin{align} \nonumber
&\displaystyle {\frac { \sqrt{2} \left( -1+2\,\sigma \right) \sin \left( \omega \right) {\it \delta_4}\,\sin \left( \Phi \right) 
}{ \sqrt{1-\sigma}}}+{\frac { \sqrt{2}\cos \left( \Phi \right) \cos \left( \omega \right) {\it \delta_1}}{2(1-\sigma)}}=0\\
&\displaystyle {\frac { \sqrt{2} \left( -1+2\,\sigma \right) \cos \left( \omega \right) {\it \delta_4}\,\sin \left( \Phi \right) 
}{\sqrt{1-\sigma}}}-{\frac { \sqrt{2}\sin \left( \omega \right) {\it \delta_1}\,\cos \left( \Phi \right) }{2 \left( 1-\sigma \right) }}=0
\label{EqZ2}
\end{align}
If $\sigma\neq 1/2$, if $\cos(\Phi)\neq0$, the only solution to \eqref{EqZ2} requires that 
\begin{equation}
\cos^{2}(\omega)+\sin^{2}(\omega)=0
\label{NoGood}
\end{equation}
which is obviously impossible. If $\cos(\Phi)=0$, a second formal solution exists; it requires that either $\sin(\Phi)=0$ also, or that $\sin(\omega)$ and $\cos(\omega)$ vanish simultaneously. Both of these conditions are incompatible. Therefore $\zeta(s)\neq 0$ if $\sigma\neq1/2$. However in the case that $\sigma=1/2$, the only solution to \eqref{EqZ2} is 
\begin{equation}
\cos(\Phi(t))=0
\label{CosPhi0}
\end{equation}
which is entirely possible and consistent with everything else that has been deduced so far - see Figure \ref{fig:CosPhi}.\newline

\begin{figure}[h] 
\centering
\includegraphics[width=.90\textwidth]{{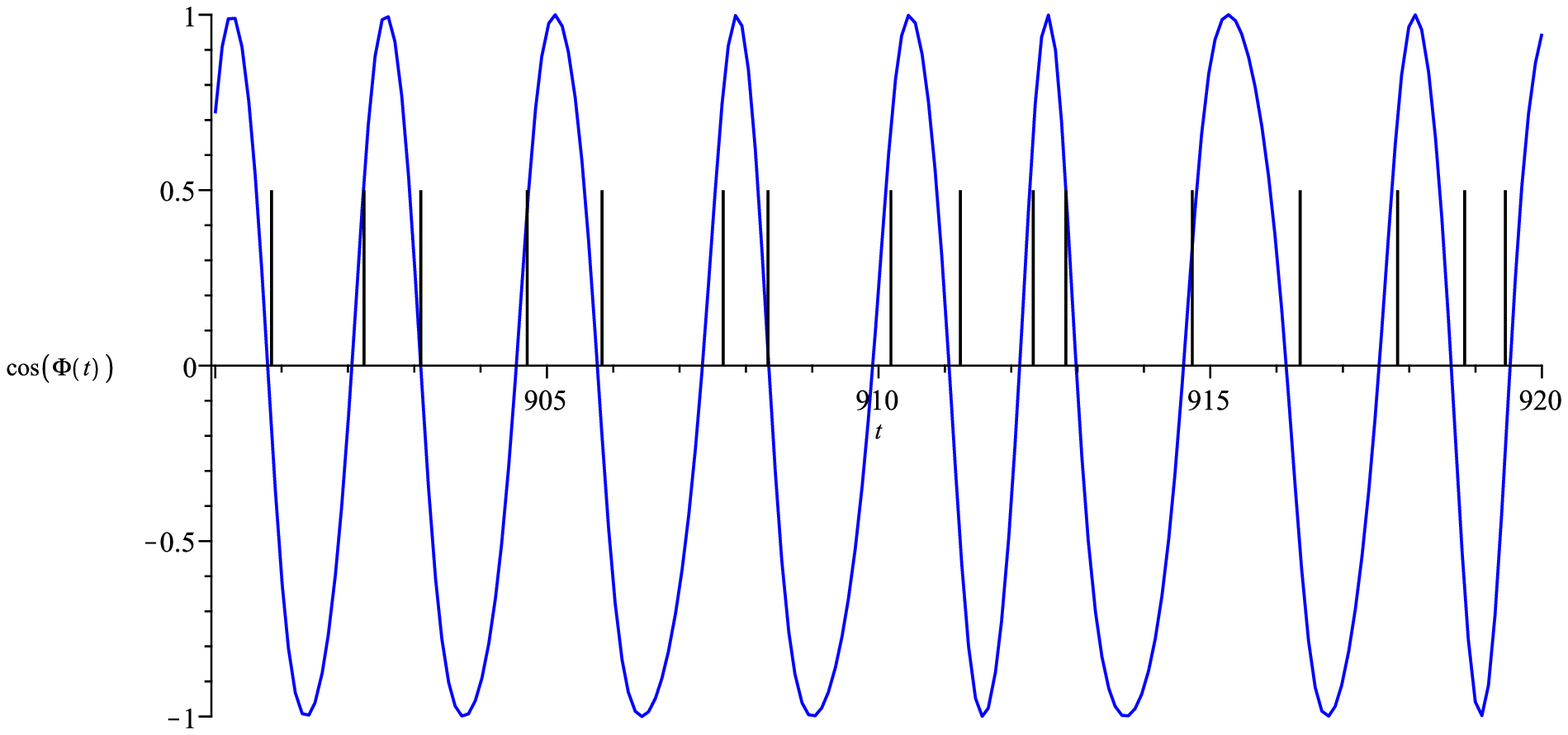}}
\caption{This is a test of \eqref{CosPhi0} when $\sigma=1/2$ near $t=900$. The vertical lines mark the location of the points $t_0$ where $\zeta(1/2+it_0)=0$, theorized to coincide with $\cos(\Phi(t_0))=0$}.
\label{fig:CosPhi}
\end{figure}

\begin{figure}[h] 
\includegraphics[width=.90\textwidth]{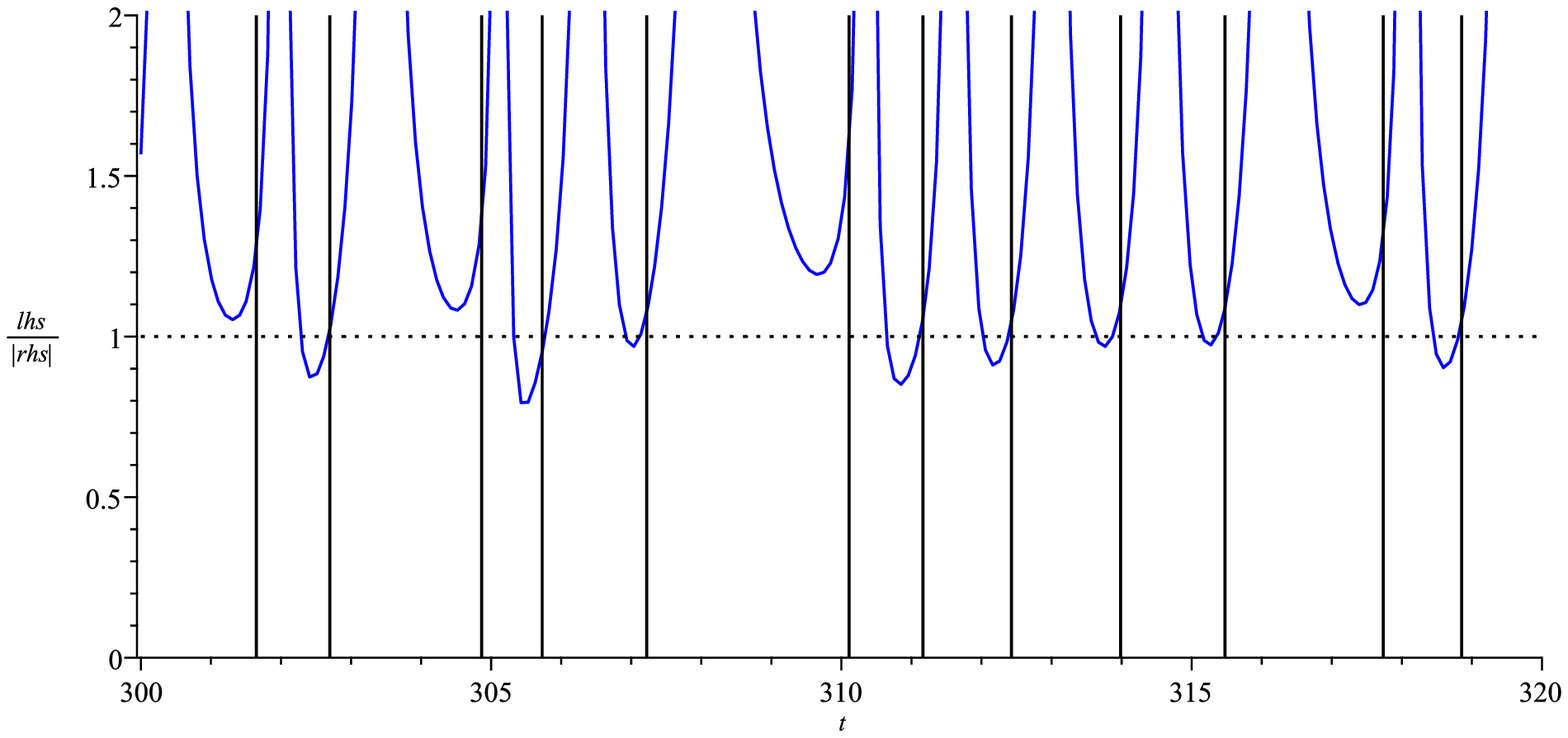}
\caption{Ratio of the left and absolute value of the right-hand sides of \eqref{Zdiff} for general values of $t>0$. Vertical lines mark the location of the zeros.}
\label{fig:Zdiff}
\end{figure}

Notice that if $t=t_0$ defines a zero on the critical line, the term $\left| \zeta \left( 1/2+it_0 \right)  \right|$ on the left-hand side of \eqref{Ceq3} vanishes, and since it is known \cite[Eq. (4.29) and Figure 3]{Milgram_FurtherEx} that at a zero defined by $t=t_0$, the argument functions 
\begin{equation}
 \alpha(t_0) -\beta(t_0) =\pm \pi/2
\label{argZero}
\end{equation}
then the accuracy of the approximation $\delta_4(1/2,t_0)=0.16$ can be tested, because \eqref{Ceq3} becomes

\begin{equation}
\displaystyle  \left| \zeta^{\prime} \left(1/2+it \right)  \right|_{t=t_0} = 2.169739306\, \left| \zeta \left( it \right)  \right| 
\mbox{}\sin \left( \phi(it) \right)( { 2}/{t})^{1/4}\bigg|_{t=t_0}\,.
\label{Zdiff}
\end{equation}

A numerical exploration of \eqref{Zdiff} is given in Figure \ref{fig:Zdiff} where the absolute ratio of the left- and right-hand sides of \eqref{Zdiff} is shown. The vertical lines mark each of the zeros $t=t_0$ for the range of $t$ under consideration, and if \eqref{Zdiff} is true, the solid (blue) curve should cross each vertical line with the value of unity, indicated by the dotted horizontal line.\newline

Finally, it needs to be mentioned that the above discussion did not consider any detailed functional dependence of $\delta_1$ and $\delta_4$. Both these parameters were introduced to model the width of the peaks (i.e. {\it extrusions}) in the products $T_1(s,v)\Upsilon_R(iv)$ and $T_4(s,v)\Upsilon_I(iv)$, and it has been shown in the Appendices that the underlying widths of $T_1(s,v)$ and $T_4(s,v)$ never vanish. Therefore neither do $\delta_1$ and $\delta_4$, and so neither affects the discussion regarding the solution of \eqref{EqZ2} and the existence of the zeros. \newline

\section{Caveats} \label{sec:Caveats}
\subsection{Universality}
An anonymous referee has written: {\it Unfortunately the fundamental flaw of the paper is the premise that (approximately) the
zeta function at $s =\sigma + it$ is related in an elementary way to the same function at
$z =\sigma_0 + it$ (i.e. for the same t but different $\sigma$). This contradicts Voronin's Universality
Theorem, which basically says the zeta function, sufficiently far up the critical strip,
mimics (to within arbitrary accuracy) all non-vanishing analytic functions in simply
connected compact subsets of the strip $1/2 <\sigma <1.$ So while the behaviour of $\zeta(\sigma+it)$
may appear to be quite regular for  $1/2 <\sigma <1$ and fixed $t$ of moderate size (indeed
monotonic for real, imaginary and absolute values), this is cannot be true as one moves
further up the strip. Indeed, one can choose an arbitrary ``wild" nonvanishing analytic
function, and an arbitrary $\epsilon > 0$ , such that the zeta function approximates that function to
within $\epsilon$ for some (probably extremely large!) fixed $t$ and  $1/2 <\sigma <1$ . Therefore,
approximations of the form \eqref{Approx4} and \eqref{Req1} cannot hold uniformly for large t.} \newline

The author comments: The fundamental result of this paper is a proof that $\zeta(s_1)$ is related to the same function $\zeta(s_2)$ integrated over a contour through the intermediary of a specific integral operator. This is a rigorous result, which, as noted earlier, has been shown \cite[Theorem 1]{GlasserXi} to be equivalent to the Cauchy Contour Integral Theorem. In this paper, a specific example has been chosen such that $s_2$ covers the 1-line through the intermediary of the functional equation \eqref{FeqId}. Subject to Proposition 1, the representation \eqref{Approx2} introduces arbitrary functions $\delta_{1,2,3,4}(\sigma,t)$ and the remainder of the paper investigates the observation that if these functions are chosen in an elementary manner, a reasonable representation of the function $\zeta(\sigma+it)$ can be found with reference to the function $\zeta(it)$ (and, by reflection, $\zeta(1-it))$ for $0<\sigma <1$. Surprisingly, the numerical examples presented suggest that the connection between the values of $\zeta(\sigma+it)$ and $\zeta(it)$ are remarkably insensitive to the particular model used for $\delta_{1,4}(\sigma,t)$. This was investigated in a superficial fashion - See Figures \ref{fig:Delta14S} and (\ref{delta4}); a detailed analysis is beyond the scope of this paper.
 Certainly, if at an extremely large value of $t$ the zeta function assumes a ``wild" value, the functions $\delta_{1,2,3,4}(\sigma,t)$ could be chosen appropriately to satisfy whatever wildness arises. However, the purpose of an approximate solution to a complex problem is to convert chaos into order, not the reverse. \newline 
 
Therefore, notice that the referee, perhaps inadvertently, raises an interesting point: in addition to the Universality Theorem, (see \cite[Chapter XI]{Titch2}), the function $\xi(s)$ (and hence $\zeta(s)$) must also obey restraints embedded in the newly derived integral equations \eqref{EqmI} and \eqref{EqpR2}. Consider the consequence(s), taking note that the following results are all rigorously proven and widely accepted:
\begin{itemize} 
\item
The Universality Theorem recognizes that somewhere ``further up the strip", the function $\zeta(s)$ will assume a ``wild" value that approximates an arbitrary function within $\epsilon$. 
\item
Equations \eqref{EqmI} and \eqref{EqpR2} impose a strict relationship between $\xi(s)$ and $\xi(2iv)$ on the $s=0+2iv$ line, and therefore a relationship exists between $\zeta(s)$ and $\zeta(1-2iv)$ through these equations, \eqref{xidef} and the functional equation. 
\item
The function $\zeta(1-2iv)$ appearing in \eqref{EqmI} and \eqref{EqpR2} via the functional equation and its polar-form is bounded by known limits on $|\zeta(1-2iv|$ through well-established theorems. Crudely $|\zeta(1-2iv)|\ll \log^{2/3}(2v)$, as well as other more strict limits given in \cite{doi:10.1093/imrn/rnx331} and \cite{Heap} - see \eqref{ZAs1} and \eqref{ZAs2} below. 
\item
The transfer functions embedded in \eqref{EqmI} and \eqref{EqpR2} are simple polynomials whose maximum $t-$value anywhere on the contour is limited to $O(t^{1})$, and these functions are the only elements in the integrand(s) that carry $t$ dependence.
\item 
The only other significant function appearing in the integrand of \eqref{EqmI} and \eqref{EqpR2} is the product $\Gamma(iv)\pi^{-iv}$ whose smoothness (continuity) and magnitude  properties are well-known. 
\item
None of the elements of the integrand in \eqref{EqmI} and \eqref{EqpR2} contain real poles or divergences and the integrals themselves are obviously convergent at both endpoints.
\item
Consequently, the integrals in \eqref{EqmI} and \eqref{EqpR2} are all limited in magnitude and constrained in their $t-$dependence and so cannot assume an arbitrarily ``wild" value or approximate an arbitrarily ``wild" function. This is a rigorous result contingent on the latter four results cited above, all of which are rigorous and true both in their own right as well as more general instances corresponding to any $c>-1$ (see Section \ref{sec:AnalCont}).
\end{itemize}
All of this suggests that an incompatibility exists among one or more of the above-cited rigorous results. Specifically, it is not apparent how an arbitrarily ``wild" function, whose magnitude exceeds the known maximum value of the product of all the components of the integrands of \eqref{EqmI} and \eqref{EqpR2} integrated over the range $[0,\infty)$, could be approximated ``within $\epsilon$" by $\zeta(\sigma+it)$ which is subject to the integral equation constraint(s) specified in the above. Further, the referee's observation that the approximate solutions \eqref{Approx4} and \eqref{Req1} do not satisfy the Universality Theorem suggests that their precursor - the (rigorously proven) integral equation itself - does not satisfy the Universality Theorem. \ifjournal \else If so, all the above could provide the foundation for a demonstration that RH is unprovable. \fi In summary, the reconciliation of the noted apparent incompatibility suggests a number of fruitful subjects for investigation - again, beyond the scope of this work. 

\subsection{Lindel{\"o}f's Hypothesis} \label{sec:Lind}

An anonymous referee has written: ``{\it \eqref{Zhalf} seems to indicate that $t^{-1/4}\zeta(1/2+it)$ is unbounded as $t\rightarrow \infty$ which we know is not true; see \cite{Bourgain}. In fact the Lindel{\"o}f hypothesis states that $\zeta(1/2+it)=O(t^{\epsilon}$) as $t\rightarrow \infty$ for all $\epsilon>0$, and if this is false, so is the Riemann hypothesis"}.\newline

The author comments:  According to Edwards \cite[Section 9.2]{Edwards}, Lindel{\"o}f's Theorem is based on Lindel{\"o}f's observation that ``there is a constant K such that 
\begin{equation}
|\zeta(\sigma+it)|<Kt^{1/2-\sigma/2}\log\,t
\label{Lobs}
\end{equation}
throughout the half-strip $0\le\Re(s)\le1$." Assuming that $|\zeta(1/2+it)|$ asymptotically obeys a power law, Lindel{\"o}f denotes by $\mu(\sigma)$
the least upper bound of the numbers $A$ such that $|\zeta(\sigma+it)|t^{-A}$ is bounded as $t\rightarrow \infty$,
shows that 
\begin{equation}
0\le\mu(1/2)\le 1/4
\label{muLims} 
\end{equation}
and hypothesizes that $\mu(1/2)=0$. This basically posits that as $t\rightarrow \infty$, $|\zeta(1/2+it)|$ increases at least faster than $O(t^0)$. Several proven refinements of \eqref{muLims} are known, over the years gradually reducing the upper inequality limit from $1/4$ to the most recent \cite{Bourgain} value $13/84=0.1548$, yielding the least upper asymptotic bound
\begin{equation}
\left| \zeta(1/2+it)\right| \ll O(t^{13/84}) \,.
\label{low_Upper}
\end{equation}

By comparison, the equivalent result obtained here (see \eqref{Req1}) is 

\begin{equation}
\displaystyle  \left| \zeta \left( 1/2+it \right)  \right| \cos \left( \phi \right)\approx 4\delta_{{1}} \left( {\scriptstyle \frac{1}{2}},t \right)\,{\frac { 
\cos \left( \Phi(t) \right)  \left| \zeta \left( 1-it \right)  \right| {t^{1/4}}{2}^{3/4}}{{\pi}^{5/4}}},
\label{XiReal}
\end{equation}
where $\cos \left( \phi \right)=(-1)^k$ counts the zeros and $-1\le \cos \left( \Phi(t) \right) \le 1$. For emphasis, $\Phi$ has been written fully as $\Phi(t)$ - see \eqref{Phi}.
With \cite[Conjecture A] {Heap} 
\begin{equation}
\max_{\substack{t\in[1,T]}}|\zeta(1+it)|\sim\exp(\gamma)\log\log\,T
\label{ZAs1}
\end{equation}
or \cite[Theorem 1]{doi:10.1093/imrn/rnx331}

\begin{equation}
\max_{\substack{t\in[\sqrt{T},T]}}|\zeta(1+it)|\ge \exp(\gamma)(\log\log\,T+\log\log\log\,T-C)
\label{ZAs2}
\end{equation}
both of which are valid for sufficiently large $T$, \eqref{XiReal} becomes
\begin{equation}
\displaystyle  \left| \zeta \left( 1/2+it \right)  \right| \cos \left( \phi \right)\approx 4\delta_{{1}} \left( {\scriptstyle \frac{1}{2}},t \right)\,{\frac { 
\cos \left( \Phi(t) \right)  {t^{1/4}}{2}^{3/4}}{{\pi}^{5/4}}}(\log \log (t)+\dots) \,.
\label{zR}
\end{equation}

Notice the similarity between \eqref{Approx4} (with \eqref{FeqId}) and \eqref{Lobs} for general $\sigma$, as well as \eqref{Lobs} and \eqref{zR} when $\sigma=1/2$. Up to this point, the functions $\delta_{1,2,3,4}(\sigma,t)$ have been simply modelled as a function of $\sigma$, although clearly $t$ dependence exists - e.g. see Figure \ref{delta4}. So, considering that $\delta_1(1/2,t)$ is expected to be a weak function of $t$ as noted previously, it is reasonable to anticipate that $t$ dependence residing in $\delta_1(\sigma,t)$ and to a lesser extent $\cos(\Phi(t))$ could affect any (dis)agreement between \eqref{Lobs}, \eqref{low_Upper}, \eqref{zR} and  Lindel{\"o}f's hypothesis.\newline

In contrast to the previous Section, where an arbitrary ``wildness" was needed to generate agreement, here it is noted that the functions $\delta_{1,2}(1/2,t)$ capture the value of the integrand of the integral equation only in close proximity to the point $v=t/2$ in \eqref{NumApprox}. Thus t-dependence of $\delta_{1,3}$ and $\delta_{2,4}$ is sensitive to the density of the real and imaginary zeros (respectively) of the function $\zeta(1-2iv)$ near a particular point on the 1-line. And, as previously noted, such $t-$dependence must persist for any $c>-1$, in which case it is useful to recall \cite{2011arXiv1112.4910A} that $\Re(\zeta(\sigma+it))>0$ if $\sigma>1.192...$ corresponding to $c>-0.904..$. No arbitrariness here. In any case, it is evident that investigations into these observations might bear fruit. Again, this is beyond the scope of this work.  

\section{Summary} \label{sec:summary}

({\bf Remark:} Throughout this work, it should be considered that the Figures presented were all obtained from the computer program Maple \cite{Maple} the accuracy of whose numerical evaluation of functions that do not lie in its mainstream is open to question. I have performed spot check comparisons between Figures from Maple and Mathematica \cite{Math} and find no grounds for concern.)\newline

Although it is known in the literature that $\xi(s)$ satisfies integral equations as cited in Section \ref{sec:Intro}, only the one presented by Patkowski \cite[Theorem (1.1)]{Patkowski} comes close to sharing the desirable properties uncovered here. In particular, because the $s-$dependence is isolated within the transfer function, it became possible to associate properties of $\xi(s)$ anywhere in the critical strip with its properties on the 1-line, which are relatively well known. One surprising prediction is that the zeros of $\xi(\sigma+it)$ and those of $\Upsilon(1-i{\it t})$ (and $\Upsilon(it)$) are closely correlated, as demonstrated by a number of Figures. Although the correlation shown in the Figures is not exact, the prediction is based on a model that leaves several questions unanswered:
\begin{itemize}
\item
What is a more reasonable model to employ for the functions $\delta_1(\sigma,t)$ and $\delta_4(\sigma,t)$?
\item
How will the two terms $\delta_2(\sigma,t)$ and $\delta_3(\sigma,t)$ affect \eqref{Approx3b} and \eqref{EqmR2} once they are incorporated into the model?
\item
What will the numerical estimates look like at values of $t$ significantly greater than those explored here, which are limited by restrictions on simple computational arithmetic?
\end{itemize}
These questions as well as  those raised in Section \ref{sec:Caveats}, implore further investigation.\newline

\section{Acknowledgements}

The author notes the contribution of two anonymous referees, one of whose remarks made it clear that the original version of this paper was poorly worded and misleading. Hence this version has been extensively rewritten to clarify all such issues. This referee also noted a number of errors in the original that have been corrected. A second referee raised a number of interesting issues that merit public discussion. Hence Section \ref{sec:Caveats}.\newline 

I am grateful to Larry Glasser who, early on provided me with a bespoke derivation of \eqref{Larry1}; shortly afterwards he outdid himself by both pointing out Romik's paper \cite{Romik}, thereby negating his invitation to pen a guest Appendix. I also thank Larry, who identified the value of the Mellin transforms, and Vini Anghel for commenting on a preliminary version of this manuscript. \ifjournal The \else Although the author is a graduate of Canada's leading university -``one of the world's greatest", as it styles itself - that university refuses to give its graduates internet access to its research library; thus the       \fi author acknowledges and thanks those authors who take the trouble to make their research work available outside a paywall. Citations to inaccessible papers are withheld from this paper's reference list. \ifjournal \else Again, although Canada is a first world country and does subsidize basic research, because the author does not have any affiliation with an educational institution, he is ineligible for a government research grant. \fi All expenses associated with this work have been borne by myself. 



\bibliographystyle{unsrt}

\bibliography{biblio}

\begin{appendices}

\addcontentsline{toc}{section}{Appendices} 
\section{Appendix: Properties of $M_{I}({\it t},v)$} \label{sec:PropsOfMI}


\begin{figure}[htbp] 
\centering
\subfloat 
{
\includegraphics[width=.49\textwidth] {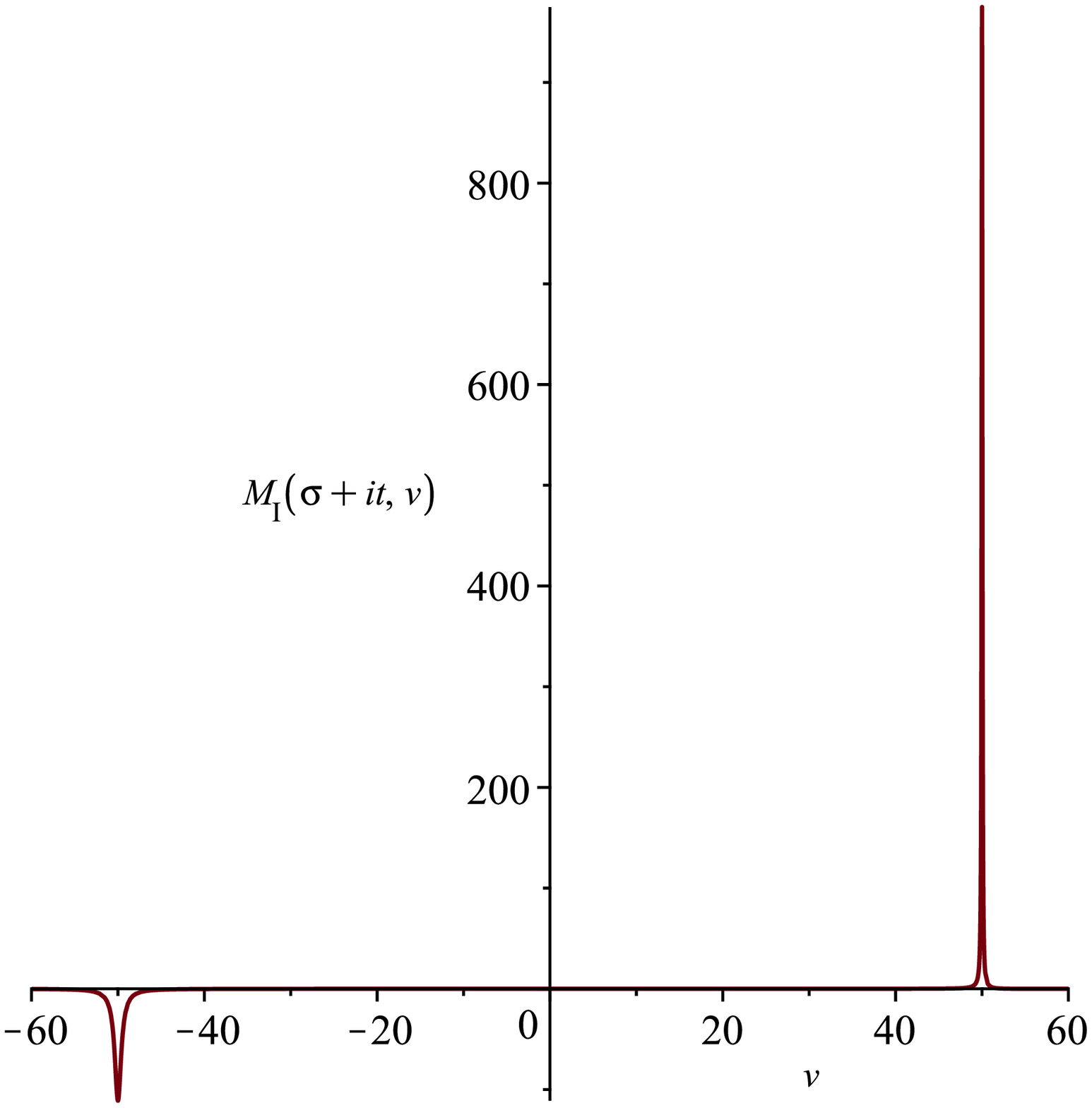}
}
\caption{This Figure summarizes the properties of the transfer function $M_{I}(\sigma+i{\it t},v)$ using $\sigma=0.1$ over a range of $v$ with ${\it t}=100$.}
\label{fig:FigMsOver} \leavevmode
\\
\subfloat []
{
\includegraphics [width=.45\textwidth]{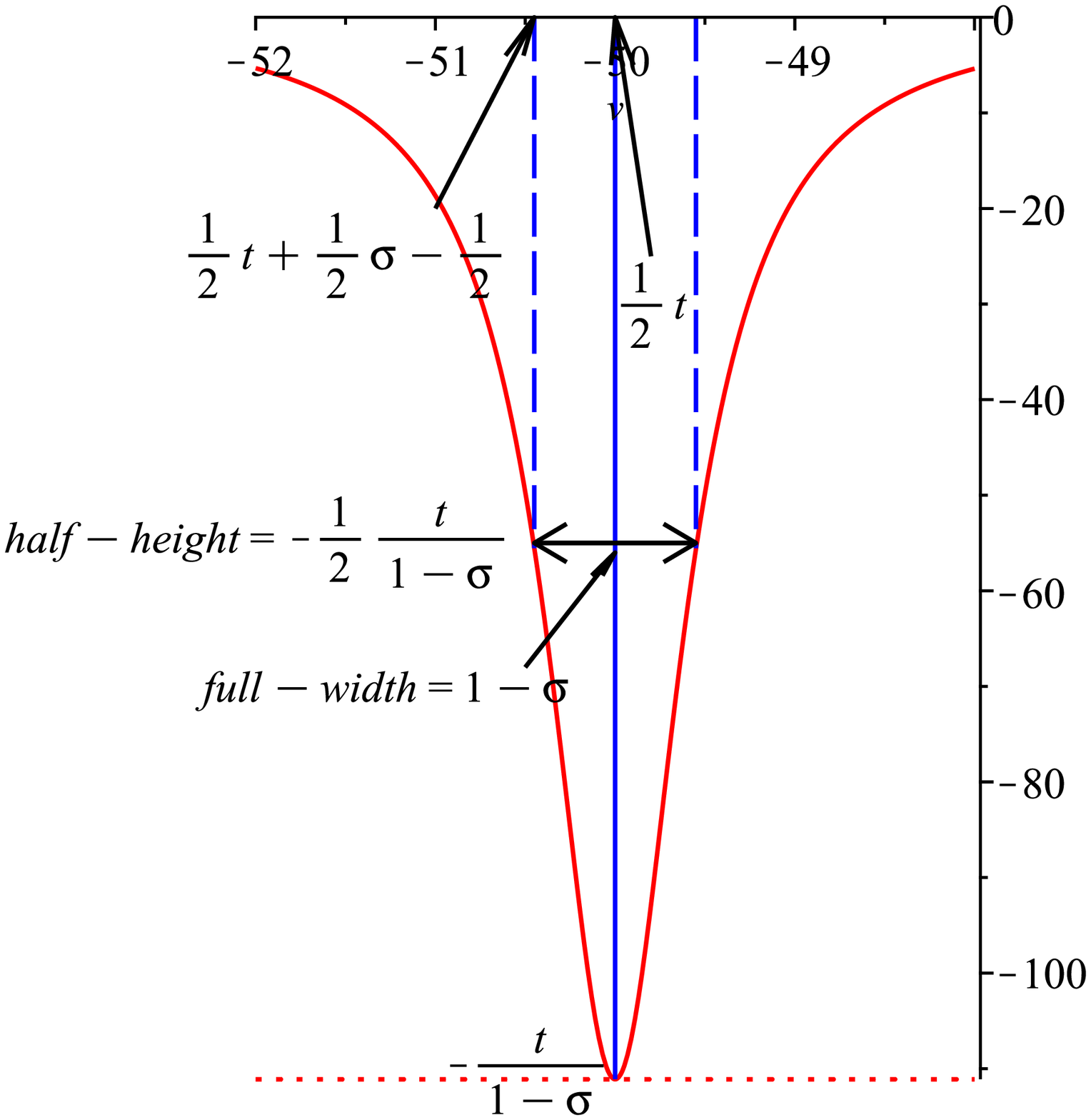} 
} 
\subfloat []
{
\includegraphics [width=.45\textwidth]{{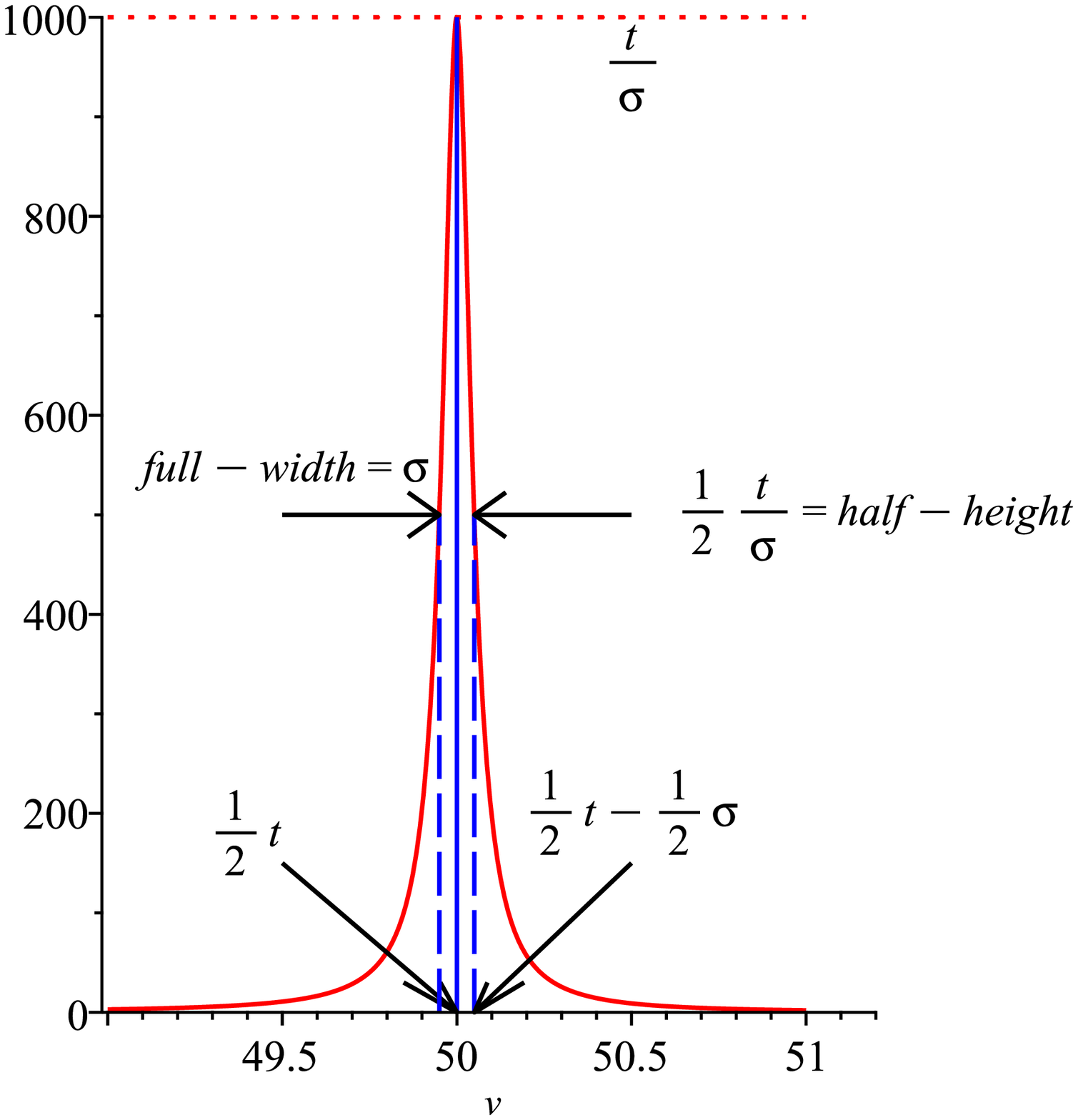}}
}
\caption{This Figure details the properties of $M_{I}(\sigma+i{\it t},v)$ when ${\it t}=100$, $\sigma=0.1$ and $v<0$ (left) and $v>0$ (right). Shown is the value of the extremum (dotted) as well as the location of one of its half-width intercepts (dashed) on the $v-$axis to first (asymptotic) order in ${\it t}^{-1}$. Note the difference in vertical scales for the same values of ${\it t}$ and $\sigma$ between these two Figures, as illustrated in Figure \ref{fig:FigMsOver} where both are shown at the same scale.}
\label{fig:FigMsHiLow}
\end{figure}

This appendix gives a detailed description of the properties of $M_I(c,s,v)$ in the case $c=-1$ and $t>0$. For simplicity, $c$ dependence is omitted from the notation. Because the interest is in the properties as a function of ${\it t}$, in this Appendix, I use $M_{I}(\sigma+i{\it t},v)\equiv M_I({\it t},v)$ and $M_R(\sigma+i{\it t},v)\equiv M_R({\it t},v)$ except in those cases where $\sigma$ dependence is considered. Because of symmetry properties about the line $c=-5/4$ (see Figure \ref{fig:Figb}), the cases $c=-1$ and $c=-3/2$ (see Section \ref{sec:C3half}) are related by the functional equation of $\zeta(s)$. With reference to its definition \eqref{ModF}, equating $\frac{\partial}{\partial\,v} M_{I}({\it t},v)=0$ and solving for $v$, identifies, as ${\it t}\rightarrow \infty$, the fact that $M_{I}({\it t},v)$ has only three extrema, those being
\begin{itemize}
\item
an extremum near $v=0$, consistent with the illustrations in Figures \ref{fig:FigMsOver} and \ref{fig:FigMsHiLow};
\item two maxima,  
\begin{equation}
v_{\pm}=\pm \frac{{\it t}}{2} \,.
\label{VsHi}
\end{equation}

\end{itemize}
At the positive maximum, we find
\begin{equation}
M_{I}({\it t}\rightarrow \infty,v_{+})=\frac{{\it t}}{\sigma}
\end{equation}
and for $v<0$ we find

\begin{equation}
M_{I}({\it t}\rightarrow \infty,v_{-})=-\frac{{\it t}}{(1-\sigma)}.
\end{equation}
These are the characteristics of two sharp ``pulses" - see Figure \ref{fig:FigMsOver} - in which case it is of further interest to determine the width. For $v>0$, solving $M_{I}({\it t}\rightarrow \infty,v)=M_{I}({\it t}\rightarrow \infty,v_{+})/2\,$ identifies the abscissa of $v$ at half-height, that being
\begin{equation}
v_h=\frac{{\it t}}{2}\,\pm\frac{\sigma}{2}.
\label{Vhalf}
\end{equation}

Thus the full-width at half-height of the function $M_{I}({\it t}\rightarrow\infty,v)$ is $ \sigma$, independent of ${\it t}$ in the asymptotic limit and thereby mimics, but does not approach, a representation of a Dirac delta function whose width vanishes. Figure \ref{fig:FigMsHiLow}, which illustrates these estimates in specific numerical form, suggests that the function values of any integrand containing $M_{I}$ (e.g. see \eqref{Beq0}) are to be evaluated at $v=\pm\,{\it t}/2$ corresponding to the sign of $v$. Embedded in this model are the factors $({\it t}/\sigma)$ and ${\it t}/(1-\sigma)$ reflecting the magnitudes of the respective peaks that depend on the sign of $v$. Included also is a width factor ${\it w_{\pm}}$ intended to capture the fact that non-zero width exists in the asymptotic limit. 

\section{Properties of $M_R(\sigma+i{\it t},v)$} \label{sec:PropsOfMR}

\begin{figure}[t!] 
\centering
\begin{subfloat} 
\centering
\includegraphics[width=.49\textwidth]{{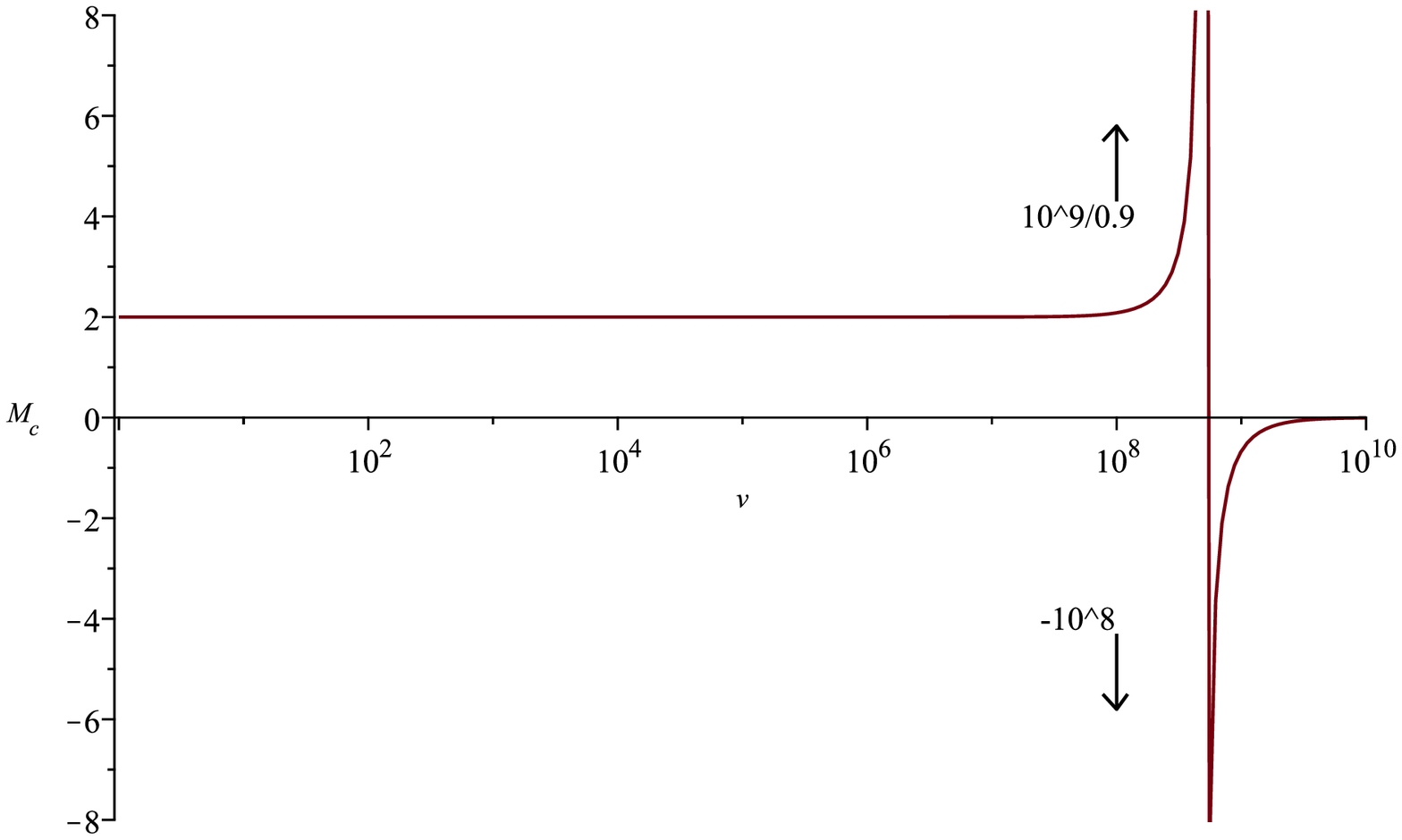}}
\caption{This Figure summarizes the properties of the transfer function $M_R(\sigma+i{\it t},v)$ using $\sigma=0.9$ over a range of $v$ with ${\it t}=10^{9}$.}
\label{fig:FigMcOver}
\end{subfloat}\leavevmode
\\
\subfloat [] 
{
\includegraphics [width=.45\textwidth]{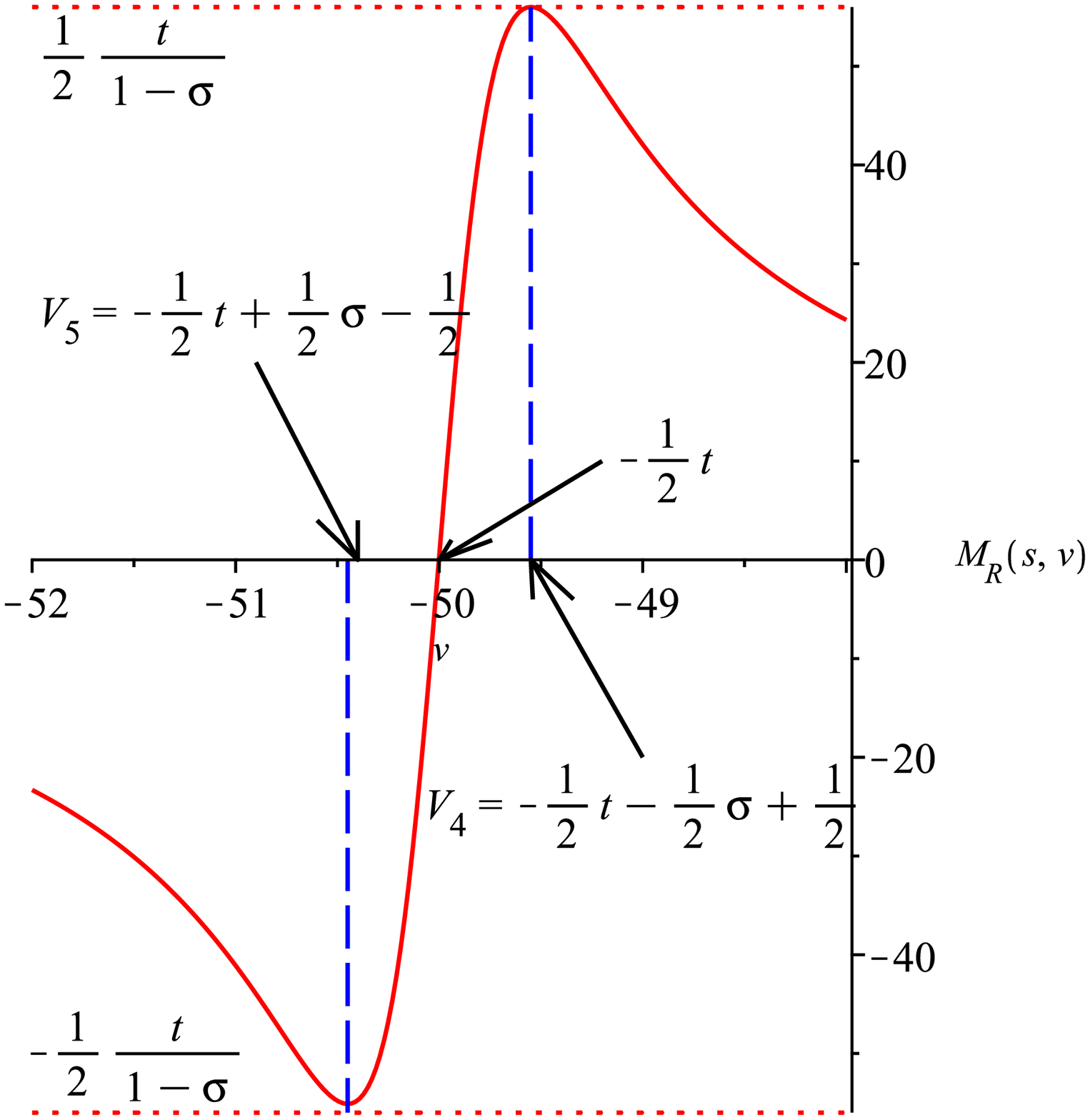} 
} 
\subfloat []
{
\includegraphics [width=.45\textwidth]{{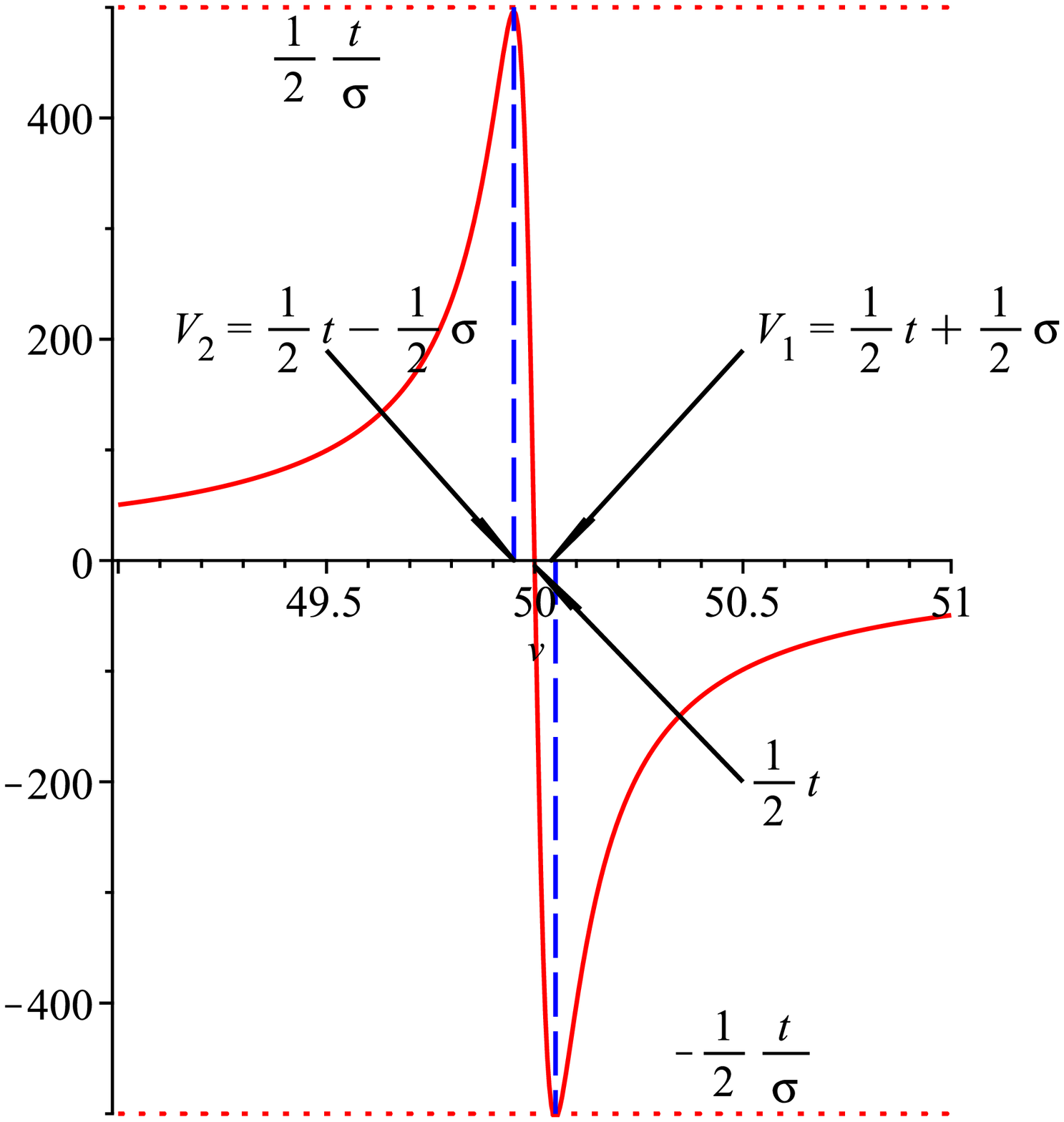}}
}
\caption{This Figure details the properties of the transfer function $M_R(\sigma+i{\it t},v)$ when ${\it t}=100$, $\sigma=0.1$ and $v<0$ (left) and $v>0$ (right). Shown are values of the extrema (dotted) as well as the location of their abscissa (dashed) on the $v-$axis to first (asymptotic) order in $1/t$. See \eqref{v1} to \eqref{v5}. Note the difference in vertical scales for the same value of ${\it t}$ and $\sigma$ between these two Figures.}
\label{fig:FigMcHiLow}
\end{figure}

The function $M_R(s,v)$, slightly more complicated, can be analyzed in a similar fashion, again for $t>0$. First of all, ignoring any extrema, in the asymptotic limit  ${\it t}\rightarrow \infty$, it is easily shown that $M_R\approx 2$ if $|v|<{\it t}/2$. An overview of $M_R(\sigma+i{\it t},v)$ can be seen in Figures \ref{fig:FigMcOver} and \ref{fig:FigMcHiLow}. Details follow below.\newline  

To locate the extrema, a solution of $\frac{\partial}{\partial\,v} M_R(\sigma+i{\it t},v)=0$ leads to a sixth order polynomial equation which does not appear to be amenable to direct analysis except in the case $\sigma=1/2$. To solve the general case $\sigma\neq 1/2$, we expect that any extremum of the function $M_R(\sigma+i{\it t},v)$ will approximately coincide with a zero of the derivative of the denominator, consistent with the presence of a nearby pole in the complex $v-$plane (see \eqref{Mx1}). Equate the derivative of this denominator to zero in the asymptotic limit ${\it t}\rightarrow \infty$, solve, then substitute the result into the expression for $\frac{\partial}{\partial\,v} M_R(\sigma+i{\it t},v)$. With a bit of experimentation using the result so-obtained as a guide, the location of three extrema can be found, those being to first order in ${\it t}^{-1}$, in order of decreasing $v$:


\begin{align} \label{v1}
&V_1=\frac{{\it t}}{2}+\frac{\sigma}{2}-\frac{\sigma(1-\sigma)}{2{\it t}}\\ \label{v2}
&V_2=\frac{{\it t}}{2}-\frac{\sigma}{2}-\frac{\sigma(1-\sigma)}{2{\it t}}\\ \label{v3}
&V_3 = \frac{(1/2-\sigma)}{4{\it t}}\\ \label{v4}
&V_4=-\frac{{\it t}}{2}+\frac{(1-\sigma)}{2}+\frac{\sigma(1-\sigma)}{2{\it t}}\\ \label{v5}
&V_5=-\frac{{\it t}}{2}-\frac{(1-\sigma)}{2}+\frac{\sigma(1-\sigma)}{2{\it t}}
\end{align}
Notice that the sets $\{V_1,V_5\}$ and $\{V_2,V_4\}$ are related by the interchange $\sigma\Leftrightarrow (1-\sigma)$ together with an overall change of sign. 
The result $V_3$ corresponds to an inflection near the origin.
\newline

The extreme values at each of these extremum locations can be found by substitution into $M_R(\sigma+i{\it t},v)$ in the limit  ${\it t}\rightarrow \infty$. For the case $v>0$ we obtain
\begin{equation}
M_R({\it t}\rightarrow\infty,\{V_1,V_2\})\approx\mp\frac{{\it t}}{2\sigma}-\frac{(1/2-\sigma)}{\sigma}\mp\frac{(2\sigma^2-2\sigma+1)}{4\sigma{\it t}}
\label{McMaxp}
\end{equation}
and for the case $v<0$ the extreme values are
\begin{equation}
M_R({\it t}\rightarrow\infty,\{V_4,V_5\})\approx\pm\frac{{\it t}}{2(1-\sigma)}+\frac{(1/2-\sigma)}{(1-\sigma)}\pm\frac{(2\sigma^2-2\sigma+1)}{4(1-\sigma){\it t}}
\label{McMaxm}
\end{equation}
all again in order, and, as in the previous case, the extrema are separated by a constant distance $\Delta v=\sigma $ or $\Delta v=1-\sigma$ according as $v>0$ or $v<0$. Using the same methods, the two midpoints can be more accurately obtained; in the asymptotic limit for $v>0$ and $v<0$ respectively, they are 
\begin{align} \label{Vmid1}
&V^{+}_{\rm mid} ={\it t}/2+3(\sigma-1)^2/(4{\it t})\\
&V^{-}_{\rm mid} =-{\it t}/2-3\sigma^2/(4{\it t})\,.
\label{Vmid2} 
\end{align}

Figure \ref{fig:FigMcOver} illustrates the above in general. $M_R({\it t},v)$ fundamentally differs from $M_{I}({\it t},v)$ because, in the asymptotic limit, it vanishes with slope $-{2{\it t}}/{(\sigma-1)^2}$ near either of the midpoints $V^{+}_{\rm mid}$ of the extrema rather than rising to an extreme value near that point, although it spans $2{\it t}$ between maximum and minimum ordinates as does $M_{I}({\it t},v)$. It can also be shown that, to a reasonable degree of approximation, $M_R({\it t},v)$ is nearly symmetric (in the asymptotic limit) about the midpoint. \newline

{\bf Remark:} In the physics lexicon, sharp peaks of the form shown in Figure \ref{fig:FigMsOver} are labelled ``Breit-Wigner" or Lorentzian. The shape is a reflection of the presence of a nearby pole in the complex $v-$plane, the closer, the sharper the peak. Choosing $c=-1$ with $0<\sigma<1$ guarantees that the peak shapes will be sharp, optimizing the approximation discussed in Section \ref{sec:Asympt}. As can be seen from Figure (\ref{fig:Figa}), when the two parametrized moving poles coalesce, they create a pole of one higher order, and this occurs only when $\sigma=1/2,t=0$. The locations of the poles can be read from the denominator parameters of \eqref{T1pa} and \eqref{T2pa}.

\section{Properties of the functions $T_{1,2,3,4}(s,v)$} \label{sec:Tprops}
The functions $T_{1,2,3,4}(s,v)$ were defined previously (see \eqref{T1Def} through \eqref{T4Def}). Since they are combinations of $M_{R,I}$ they share similar properties which can be visualized by examining Figure \ref{fig:P1234}. As before, the formal equations defining these functions, available from \eqref{Ms1}, \eqref{Mc1}, \eqref{T1Def} - \eqref{T4Def}
involve $6^{th}$ order rational polynomials which do not lend themselves to analytic solution in order to determine the explicit values of the extrema and their locations. However, in the limit ${\it t}\rightarrow\infty$, it is possible to obtain reasonable approximations by experimentation. It is likely that these estimates could be improved in future; they are listed below, employing to the following notation:
\begin{itemize}
\item $V_{j}(min/max)$ refers to the abscissa of an extremum projected onto the $v-$axis ; 
\item$T_{j}$ gives the magnitude (amplitude) of the extremum;
\item $Z_{j}$ locates any points where the function crosses the $v-$axis;
\item $H_{j}$ locates points corresponding to the half-height, and
\item $W_{j}$ gives the full-width at half-height,
\end{itemize}
all with reference to function $T_{j}(s,v)$, $j=1,2,3,4$. Note that in some cases, the minimum and maximum reverse or vanish according to whether $\sigma>1/2$, $\sigma<1/2$ or $\sigma=1/2$. 

\begin{align} \label{T1Min}
&\displaystyle T_{{1}} \left( {\it Min} \right) ={\frac {{\it t}}{\sigma\, \left( \sigma-1 \right) }}+\frac{3}{4{\it t}}\\ 
&\displaystyle H_{1}={\it t}/2\pm \,1/4\, \sqrt{-2-8\,{\sigma}^{2}+2\, \sqrt{4\,\sigma\, \left( \sigma-1 \right)  \left( 5\,{\sigma}^{2}-5\,\sigma+2 \right) +1}+8\,\sigma}
\label{H1}\\
&\displaystyle W_{1}  = 1/2\, \sqrt{-2-8\,{\sigma}^{2}+2\, \sqrt{4\,\sigma\, \left( \sigma-1 \right)  \left( 5\,{\sigma}^{2}-5\,\sigma+2 \right) +1}+8\,\sigma}
\label{W1} \\
&\displaystyle V_{{2}} \left( {\it Min}/ {\it Max} \right) ={\it t}/2\,\mp\left(1/4- \left( 1/2-\sigma \right) ^{2}\right)\label{V2MinMax} \\
%
%
&\displaystyle T_{{2}} \left( { {{\it Min}/}{{\it Max}}} \right) =\mp\,{\frac {{\it t}\left( -4\,\sigma\, \left( 1-\sigma \right)  \left( 2\,\sigma-1 \right) ^{2}+2 \right) }{\sigma\, \left( 1-\sigma \right)  \left( 4\,{\sigma}^{2}+1 \right)  \left( 4\, \left( 1-\sigma \right) ^{2}+1 \right) 
}} \label{T2MinMax} \\
&\displaystyle T_{{3}} \left( {\it Min}/{\it Max} \right) =\mp\,{\frac {2  {\it t}\left(\sigma -1/2 \right)}{\sigma\, \left( {\sigma}^{2}-2\,\sigma+2 \right) 
\mbox{} \left( 1-\sigma \right)  \left( {\sigma}^{2}+1 \right) }} \label{T3MinMax}\\
%
%
&\displaystyle V_{{3}} \left( {\it Min}/{\it Max} \right) ={\it t}/2\mp\,\sigma\, \left( 1-\sigma \right)/2 \label{V3MinMax} \\
%
%
&\displaystyle T_{{4}} \left( {\it Min} \right) =2\,{\it t}{\frac { \left(\sigma- 1/2 \right) }{\sigma\, \left( 1-\sigma \right) }} \label{T4Min}\\
&\displaystyle T_{{4}} \left( {\it Max} \right) =-2{\it t}\,{\frac { \left(\sigma-1/2 \right) }{1+2\, \sqrt{\sigma\, \left( 1-\sigma \right) }
}} \label{T4Max}\\
%
%
&\displaystyle V_{{4}} \left( {\it Min}/ {\it Max} \right) ={\it t}/2\,\mp\,1/2 \sqrt{\sigma\, \left( 1-\sigma \right) + \sqrt{\sigma\, \left( 1-\sigma \right) }} \label{V4MinMax} \\
&\displaystyle Z_{{4}} \left( {\it Min}/ {\it Max} \right) ={\it t}/2\mp1/2\, \sqrt{\sigma\, \left( 1-\sigma \right) }\label{Z4MinMax} \\
&\displaystyle {\it H_{4}}={\it t}/2\pm\,1/4\, \sqrt{-2+2\, \sqrt{4\, \left( \sigma-1 \right) ^{2}{\sigma}^{2}+1}}
\label{H4}\\
&\displaystyle {\it W_{4}}=1/2\, \sqrt{-2+2\, \sqrt{4\, \left( \sigma-1 \right) ^{2}{\sigma}^{2}+1}}
\label{W4}
\end{align}
%
\begin{figure}[htbp] 
\begin{subfloat} []
{
\includegraphics [width=.45\textwidth]{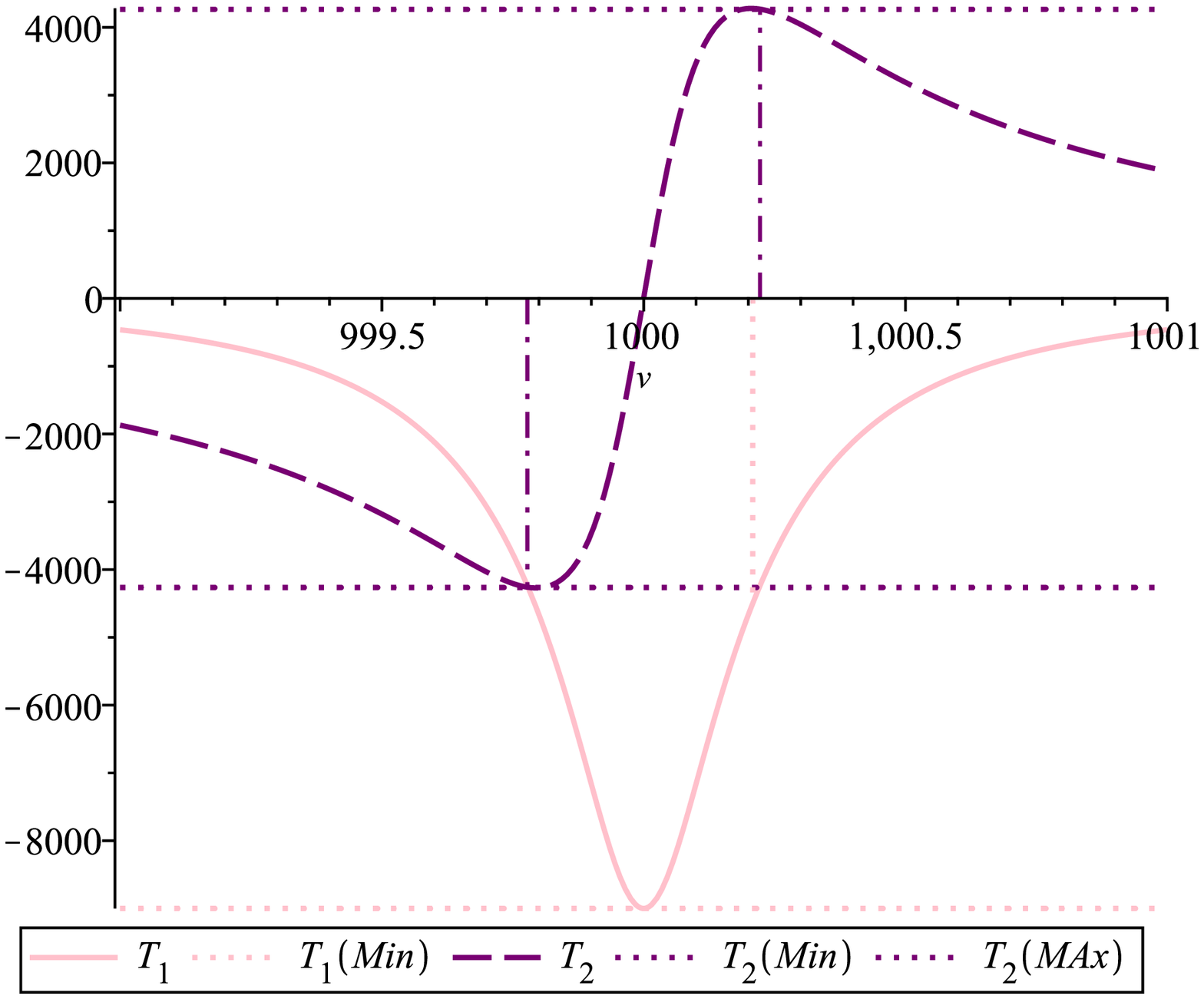} \label{fig:P1P2}
} 
\subfloat []
{
\includegraphics [width=.45\textwidth]{{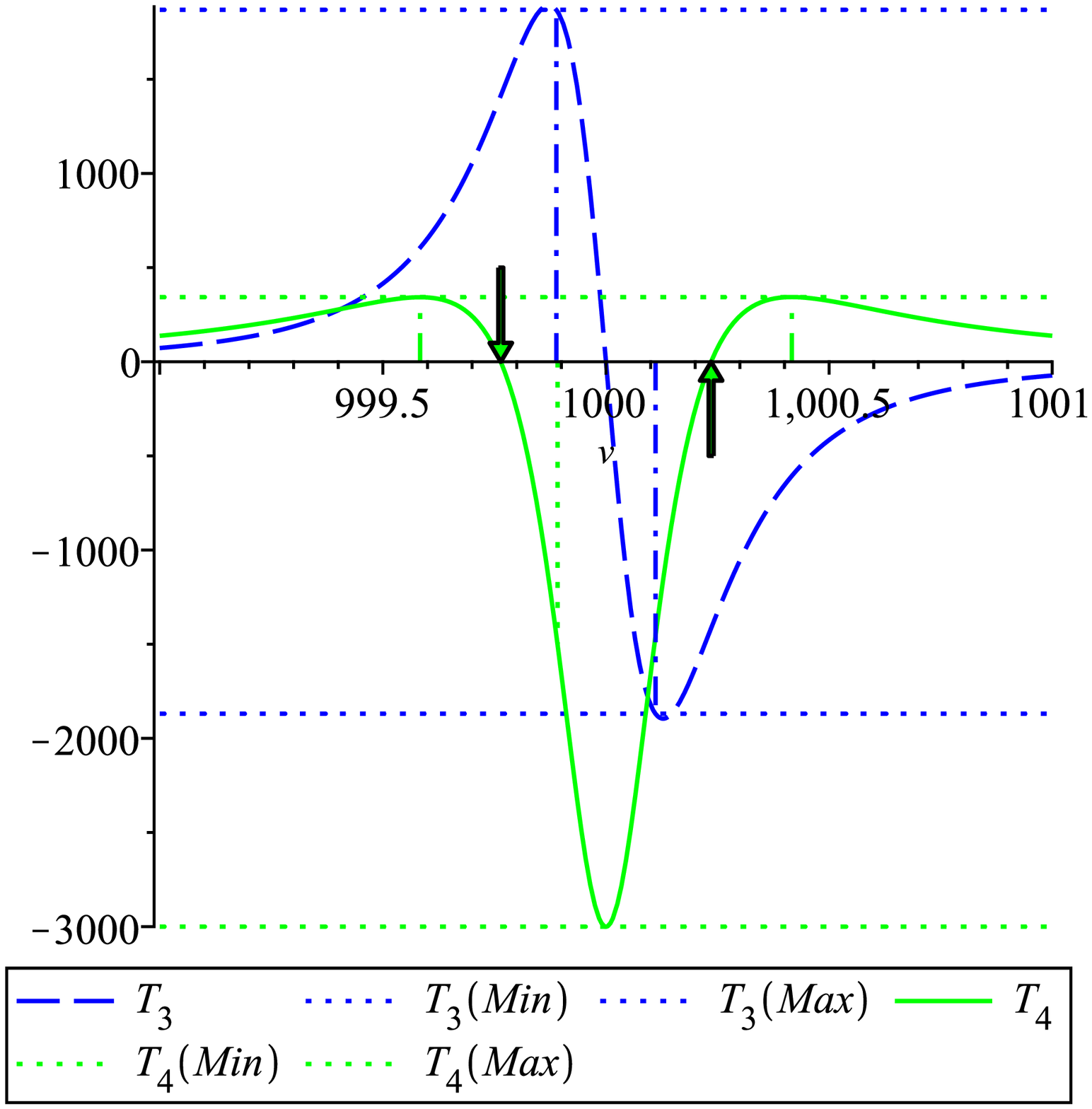}}
}
\label{fig:P3P4}
\end{subfloat}
\caption{This Figure illustrates the nature of the functions $T_1(s,v)$ to $T_4(s,v)$ using $\sigma=1/3, {\it t}=2000$. The extremum values are shown to demonstrate the accuracy of the estimates \eqref{T1Min}, \eqref{T2MinMax}, \eqref{T3MinMax}, \eqref{T4Min} and \eqref{T4Max}. The location points \eqref{V2MinMax}, \eqref{V3MinMax} and \eqref{V4MinMax} are also shown (dashdot). The arrows correspond to \eqref{Z4MinMax} and the vertical dotted lines show the location (\eqref{H1} and\eqref{H4}) of the half-height line.}
\label{fig:P1234}
\end{figure}

There are several important points to note about these functions:
\begin{itemize}
\item
To at least first order in ${\it t}$ each is more-or-less symmetric near the point $v={\it t}/2$;
\item 
Other than near $v={\it t}/2$ the functions otherwise vanish as $O({\it t}^{-2})$ or $O({\it t}^{-3})$ - for an overview see Figures \ref{fig:FigMsOver} and \ref{fig:FigMcOver}. Therefore the properties of $\xi(s)$ are determined to a large extent by the properties of the integrand of \eqref{EqmI} and \eqref{EqpR2}  near the point $v={\it t}/2$, although other non-vanishing, but ultimately cancelling terms obfuscate this result numerically;
\item
The function extrema approximate both a unit delta function and it's derivative in the limit ${\it t}\rightarrow\infty$;
\item
The magnitude estimates given here are all expected to be accurate to $O({\it t})$; the locations (and therefore widths) less so; in Section (\ref{sec:ZsApprox}) both are used to obtain an asymptotic estimate of $\zeta(s)$ which is therefore expected to be reasonably accurate, but not exact because the widths are only estimated. It is the magnitude(s) of the pole terms embedded in $T_1(s)$ and $T_4(s)$ that govern the leading order in ${\it t}$; since these are exact to leading order in ${\it t}$, so will be the accuracy of $\xi(s)$ (and hence $\zeta(s)$) to the leading order in ${\it t}$. The widths govern the accuracy of the approximation details shown in various Figures in the main text, but do not appear to affect the derivation of the location of the zeros. 
\end{itemize}

\section{Listing of Sums and Moments} \label{sec:RomSums}
The following sums were calculated by applying higher order derivative with respect to $x$ to Romik's Theorem 1 \cite{Romik} at $x=\pi$. They are needed to obtain the moments of $\Upsilon(2iv)$ below, for use in Section \ref{sec:Backg}. Note: $\Gamma(3/4)^n$ means $(\Gamma(3/4))^n$.

\begin{align} \label{MelSum1}
\displaystyle  \sum _{n=1}^{\infty }-{n}^{2}{{\rm e}^{-\pi\,{n}^{2}}}&=-{\frac {1}{{8\,\pi}^{3/4}\Gamma \left( 3/4 \right) }}
\\ \label{MelSum2}
\displaystyle \sum _{n=1}^{\infty }{n}^{4}{{\rm e}^{-\pi\,{n}^{2}}}&={\frac {3+{\frac {{\pi}^{4}}{2 \Gamma \left( 3/4 \right)^{8}}}}{32{\pi}^{7/4}\Gamma \left( 3/4 \right) 
\mbox{}}}
\\ \label{MelSum3}
\displaystyle \sum _{n=1}^{\infty }-{n}^{6}{{\rm e}^{-\pi\,{n}^{2}}}&=-{\frac {15}{128\,{\pi}^{11/4}\Gamma \left( 3/4 \right) }}-{\frac {15\,{\pi}^{5/4}}{256\, \Gamma \left( 3/4 \right)^{9}}}
\\ \label{MelSum4}
\displaystyle \sum _{n=1}^{\infty }{n}^{8}{{\rm e}^{-\pi\,{n}^{2}}}&={\frac {105}{512\,{\pi}^{{\frac{15}{4}}}\Gamma \left( 3/4 \right) }}
\mbox{}+{\frac {105\,\sqrt [4]{\pi}}{512\,\Gamma \left( 3/4 \right)^{9}}}-{\frac {{\pi}^{{\frac{17}{4}}}}{2048\,\Gamma \left( 3/4 \right) ^{17}}}
\\ \label{MelSum5}
\displaystyle \sum _{n=1}^{\infty }-{n}^{10}{{\rm e}^{-\pi\,{n}^{2}}}&=-{\frac {945}{2048\,{\pi}^{{\frac{19}{4}}}\Gamma \left( 3/4 \right) }}
\mbox{}-{\frac {1575}{2048\,{\pi}^{3/4}\Gamma \left( 3/4 \right)  ^{9}}}+{\frac {45\,{\pi}^{{\frac{13}{4}}}}{8192\,  \Gamma \left( 3/4 \right)  ^{17}}}
\\ \label{MelSum6}
\displaystyle \sum _{n=1}^{\infty }{n}^{12}{{\rm e}^{-\pi\,{n}^{2}}}&={\frac {10395}{8192\,{\pi}^{{\frac{23}{4}}}\Gamma \left( 3/4 \right) }}
\mbox{}+{\frac {51975}{16384\,{\pi}^{7/4}\Gamma \left( 3/4 \right)   ^{9}}}-{\frac {1485\,{\pi}^{9/4}}{32768\, \Gamma \left( 3/4 \right)  ^{17}}}\\  \nonumber
&\;\;\;\;+{\frac {51\,{\pi}^{{\frac{25}{4}}}}{65536\, \Gamma \left( 3/4 \right)^{25}}}
\\ \label{MelSum7}
\displaystyle \sum _{n=1}^{\infty }-{n}^{14}{{\rm e}^{-\pi\,{n}^{2}}}&=-{\frac {135135}{32768\,{\pi}^{{\frac{27}{4}}}\Gamma \left( 3/4 \right) 
\mbox{}}}-{\frac {945945}{65536\,{\pi}^{11/4}\Gamma \left( 3/4 \right)  ^{9}}}+{\frac {45045\,{\pi}^{5/4}}{131072\, \Gamma \left( 3/4 \right) ^{17}}}
\mbox{} \\ \nonumber
&\;\;\;\;-{\frac {4641\,{\pi}^{{\frac{21}{4}}}}{262144\, \Gamma \left( 3/4 \right) ^{25}}}
\end{align}

From the above sums, it is possible to obtain the moments of $\Upsilon(2iv)$ by evaluating the derivatives of \eqref{Rans0} at $x=\pi$ as follows:

\begin{align} \label{Mom1}
\displaystyle \int_{0}^{\infty }\!v{\it \Upsilon_{I}} \left( 2\,iv \right) \,{\rm d}v&=-\pi/4\, \left( {\frac {{\pi}^{1/4}}{2\Gamma \left( 3/4 \right) }}-1 \right) 
\\ \label{Mom2}
\displaystyle \int_{0}^{\infty }\!{v}^{2}{\it \Upsilon_{R}} \left( 2\,iv \right) \,{\rm d}v&=\pi/8+{\frac {{\pi}^{5/4}}{32\Gamma \left( 3/4 \right) }}
\mbox{}-{\frac {{\pi}^{{\frac{21}{4}}}}{64\, \Gamma \left( 3/4 \right)  ^{9}}}
\\ \label{Mom3}
\displaystyle \int_{0}^{\infty }\!{v}^{3}{\it \Upsilon_{I}} \left( 2\,iv \right) \,{\rm d}v&={\frac {3\,{\pi}^{{\frac{21}{4}}}
\mbox{}-10\,{\pi}^{5/4}  \Gamma \left( 3/4 \right)   ^{8}-16\,\pi\,  \Gamma \left( 3/4 \right)^{9}}{256\, \Gamma \left( 3/4 \right)^{9}}}
\\ \label{Mom4}
\displaystyle \int_{0}^{\infty }\!{v}^{4}{\it \Upsilon_{R}} \left( 2\,iv \right) \,{\rm d}v&={\frac {-{\pi}^{{\frac{37}{4}}}
\mbox{}-76\,{\pi}^{{\frac{21}{4}}}  \Gamma \left( 3/4 \right) ^{8}+68\,{\pi}^{5/4} \Gamma \left( 3/4 \right) ^{16}
\mbox{}-64\,\pi\,  \Gamma \left( 3/4 \right) ^{17}}{2048\,  \Gamma \left( 3/4 \right) ^{17}}}
\\ \label{Mom5}
\displaystyle \int_{0}^{\infty }\!{v}^{5}{\it \Upsilon_{I}} \left( 2\,iv \right) \,{\rm d}v&={\frac {5\,{\pi}^{{\frac{37}{4}}}
\mbox{}+420\,{\pi}^{{\frac{21}{4}}}\Gamma \left( 3/4 \right)^{8}-484\,{\pi}^{5/4} \Gamma \left( 3/4 \right)^{16}
\mbox{}+128\,\pi\,\Gamma \left( 3/4 \right)^{17}}{8192\, \Gamma \left( 3/4 \right)^{17}}}
\\ \nonumber
\displaystyle \int_{0}^{\infty }\!{v}^{6}{\it \Upsilon_{R}} \left( 2\,iv \right) \,{\rm d}v&=\frac{1}{65536\, \Gamma \left( 3/4 \right)^{25}}\left( -51\,{\pi}^{{\frac{53}{4}}}
\mbox{}-350\,{\pi}^{{\frac{37}{4}}}\Gamma \left( 3/4 \right)^{8}-11644\,{\pi}^{{\frac{21}{4}}} \Gamma \left( 3/4 \right)  ^{16}
\right.\\ \label{Mom6}
&\left. \hspace{3.5cm}
+5768\,{\pi}^{5/4} \Gamma \left( 3/4 \right)^{24}+512\,\pi\, \Gamma \left( 3/4 \right)^{25}\right)
\\ \nonumber
\displaystyle \int_{0}^{\infty }\!{v}^{7}{\it \Upsilon_{I}} \left( 2\,iv \right) \,{\rm d}v&=\frac{1}{262144\,\Gamma \left( 3/4 \right) ^{25
}}\times\left( 357\,{\pi}^{{\frac{53}{4}}}
\mbox{}+2590\,{\pi}^{{\frac{37}{4}}}  \Gamma \left( 3/4 \right) ^{8}+93492\,{\pi}^{{\frac{21}{4}}} \Gamma \left( 3/4 \right)  ^{16}
\right.\\ 
&\left. \hspace{3.5cm}
-54760\,{\pi}^{5/4} \Gamma \left( 3/4 \right)^{24}-1024\,\pi\, \Gamma \left( 3/4 \right)^{25}\right)
\label{Mom7}
\end{align}

\end{appendices}
\end{flushleft}
\end{maplegroup}
\end{document}